# Node polynomials for curves on surfaces
## by Thomas Dedieu

## 1 – Introduction and statement of results

The main theme of this text is the following set of formulae, first found by Vainsencher [38] (for $\delta \leqslant 6$) and shortly after improved by Kleiman and Piene [20, 21, 22, 23], counting $\delta$-nodal curves in a complete linear system on a smooth surface, if $\delta \leqslant 8$ and the corresponding line bundle is sufficiently positive.

**(1.1)** The formulae involve the *complete exponential Bell polynomials* $P_r \in \mathrm{Z}[X_1, \ldots, X_r]$, $r \in \mathbf{N}$, defined by the formal identity

$$\sum_{r \geqslant 0} P_r \frac{T^r}{r!} = \exp\left(\sum_{q \geqslant 1} X_q \frac{T^q}{q!}\right)$$

---

in $\mathbf{Q}(X_1, \ldots, )[[T]]$. Thus $P_0 = 1$, $P_1 = X_1$, $P_2 = X_1^2 + X_2$, $P_3 = X_1^3 + 3X_2X_1 + X_3$, and so forth. For all $r$, the polynomial $P_r$ is weighted homogeneous of degree $r$, if one assigns each $X_i$ the weight $i$. For convenience we also set $P_r = 0$ for all negative integers $r$.

**(1.2)** Let $S$ be a smooth irreducible complex projective surface, and $L$ be a line bundle on $S$. Consider the four Chern numbers

$$d = L^2 \; ; \quad k = L \cdot K_S \; ; \quad s = K_S^2 \; ; \quad x = c_2(S),$$

and set

$$\begin{aligned}
a_1 &= 3d + 2k + x \; ; \\
a_2 &= -42d - 39k - 6s - 7x \; ; \\
a_3 &= 1380d + 1576k + 376s + 138x \; ; \\
a_4 &= -72360d - 95670k - 28842s - 3888x \; ; \\
a_5 &= 5225472d + 7725168k + 2723400s + 84384x \; ; \\
a_6 &= -481239360d - 778065120k - 308078520s + 7918560x \; ; \\
a_7 &= 53917151040d + 93895251840k + 40747613760s - 2465471520x \; ; \\
a_8 &= -7118400139200d - 13206119880240k - 6179605765200s + 516524964480x.
\end{aligned}$$

**(1.3) Theorem** (Vainsencher, Kleiman–Piene). *Let $\delta \leqslant 8$ be an integer, and assume that $L$ may be written as $L = mA + B$ for an integer $m \geqslant 3\delta$ and line bundles $A$ and $B$, such that $A$ is very ample and $B$ is globally generated.[1] Fix $n - \delta$ general points $p_1, \ldots, p_{n-\delta} \in S$, where $n$ is the dimension of the complete linear system $|L|$. The number of $\delta$-nodal curves among members of $|L|$ passing through $p_1, \ldots, p_{n-\delta}$ is given by the polynomial expression in $d, k, s, x$:*

(1.3.1) $$N_\delta = \frac{1}{\delta!} P_\delta(a_1, \ldots, a_\delta).$$

These formulae are obtained as an application of a more general result, that I shall now describe.

**(1.4)** Let $\pi : F \to Y$ be a smooth projective family of surfaces, where the base $Y$ is equidimensional and Cohen–Macaulay. Let $D$ be a relative effective divisor on $F/Y$. Consider the three Chern classes

$$v = c_1\big(\mathcal{O}_F(D)\big) \quad \text{and} \quad w_j = c_j\big(\Omega^1_{F/Y}\big) \text{ for } j = 1, 2,$$

and set

$$\begin{aligned}
x_2 &= v^3 + v^2 w_1 + v w_2 \; ; \\
x_3 &= v^6 + 4v^5 w_1 + 5v^4(w_1^2 + w_2) + v^3(2w_1^3 + 11w_1 w_2) + v^2(6w_1^2 w_2 + 4w_2^2) + 4v w_1 w_2^2 \; ; \\
x_4 &= v^{10} + 10v^9 w_1 + v^8(40w_1^2 + 15w_2) + v^7(82w_1^3 + 111w_1 w_2) \\
&\quad + v^6(91w_1^4 + 315w_1^2 w_2 + 63w_2^2) + v^5(52w_1^5 + 29w_1^3 w_2 + 324w_1 w_2^2) \\
&\quad + v^4(12w_1^6 + 282w_1^4 w_2 + 593w_1^2 w_2^2 + 85w_2^3) + v^3(72w_1^5 w_2 + 464w_1^3 w_2^2 + 259w_1 w_2^3) \\
&\quad + v^2(132w_1^4 w_2^2 + 246w_1^2 w_2^3 + 36w_2^4) + v(72w_1^3 w_2^3 + 36w_1 w_2^4).
\end{aligned}$$

We shall use these polynomial expressions to define recursively eight Chern classes $b_1, \ldots, b_8$ on $F$. To do so we define, for all polynomial expressions $R(v, w_1, w_2)$ and integers $i$, a new polynomial expression $Q_{R,i}(v, w_1, w_2)$ as follows: (i) consider a new variable $e$, and form the polynomial expression

$$R(v - ie, \, w_1 + e, \, w_2 - e^2) \; ;$$

---

[1] if $\delta = 8$, one needs to add the technical assumption that among curves in $|L|$ with an ordinary 4-tuple point, the analytic type of the 4-tuple point is not constant, i.e., the cross-ratio of the four tangents varies.

(ii) perform the Euclidean division of this polynomial by $e^3 + w_1 e^2 + w_2 e$ ; (iii) consider the coefficient in $e^2$ of the remainder of this division, and define $Q_{R,i}(v, w_1, w_2)$ to be its opposite. The geometric meaning of these operations will be apparent later in the text. Then, we define for all $s = 0, \ldots, 7$:

$$
\begin{aligned}
(1.4.1) \quad b_{s+1} = {} & P_s\big(Q_{b_1,2}, \ldots, Q_{b_s,2}\big) \cdot x_2 \\
& - s(s-1)(s-2) \cdot P_{s-3}\big(Q_{b_1,3}, \ldots, Q_{b_{s-3},3}\big)\big) \cdot x_3 \\
& + 3281 \cdot 7! \cdot P_{s-7} \cdot x_4.
\end{aligned}
$$

In this definition, the polynomials $P_r$ are as defined in (1.1) above; note that for $s = 0, 1, 2$, only the first line is involved (as the other two are zero), and for $s = 3, 4, 5, 6$, only the first two lines are involved; for $s = 7$, $P_{s-7}$ is merely $P_0 = 1$.

**(1.5)** In order to state the general result leading to (1.3), we need to consider *(multi-)equisingularity classes*. Informally this is a collection of planar singularity types, e.g., $\delta A_1$ denotes a collection of $\delta$ nodes, $(\delta - 1)A_1 + A_2$ denotes a collection of $\delta - 1$ nodes and one cusp, and so on. Roughly speaking a singularity type is the same as the name of a singularity; for instance "ordinary quadruple point" is a singularity type, and all ordinary quadruple points are of this type, regardless of the cross-ratio of the tangents of their four local branches. A precise definition will be given in Section 4, where we will see in particular that equisingularity classes may be encoded by Enriques diagrams.

Let $S$ be a multi-equisingularity class. In the setup of (1.4), we can form the subset $Y(S) \subseteq Y$ consisting of those $y$ for which the curve $D_y$ in the surface $F_y$ has singularities of the type defined by $S$. It is a constructible set in $Y$, which has an expected codimension depending only on $S$, and that we shall denote by $\operatorname{expcod}(S)$; thus, $\operatorname{codim}_Y(Y(S))$ is smaller or equal than $\operatorname{expcod}(S)$, and if equality holds we say that $Y(S)$ has the expected codimension in $Y$. For convenience, we agree that the empty set has the expected codimension in $Y$.

Moreover, we denote by $Y(\infty)$ the subset of $Y$ consisting of those points $y$ such that the curve $D_y$ is non-reduced. It is a constructible set as well.

**(1.6) Theorem** (Vainsencher, Kleiman–Piene). *Consider the setup of (1.4). Let $\delta \leqslant 8$ be a non-negative integer, and assume that the two following conditions hold:*[2]
*(i) $\operatorname{codim}\big(Y(\infty)\big) \geqslant \delta + 1$;*
*(ii) for all multi-equisingularity classes $S$, $\operatorname{codim}_Y\big(Y(S)\big) \geqslant \min\big(\delta + 1, \operatorname{expcod}(S)\big)$.*

*Then $Y(\delta A_1)$, the locus of $\delta$-nodal curves in $Y$, has pure codimension $\delta$ (if non-empty); its closure $\overline{Y(\delta A_1)}$ is the support of a natural nonnegative cycle $U(D, \delta)$, and the following equality in the Chow ring of $Y$ holds:*

$$
(1.6.1) \qquad \big[U(D, \delta)\big] = \frac{1}{\delta!} P_\delta\big(\pi_* b_1, \ldots, \pi_* b_\delta\big),
$$

*where $P_\delta$ is the $\delta$-th Bell polynomial, and $b_1, \ldots, b_\delta$ are the rational equivalence classes on $F$ defined in (1.4).*

To save the reader from a possible moment of confusion, let me emphasize that the Bell polynomials appear twice in this story, first in the definition of the classes $b_1, \ldots, b_\delta$, and second in the expression of $[U(D, \delta)]$ in terms of these very $b_1, \ldots, b_\delta$.

In setup (1.4) one considers a family of curves in a possibly non-constant family of surfaces. This is indeed necessary for the approach of Vainsencher and Kleiman–Piene to Theorem (1.3), in which one considers, given a fixed surface $S$, the family $\pi : F \to Y = S$ such that for all $s \in S$, the surface $F_s$ is the blow-up of $S$ at the point $s$. This also gives the possibility of further applications, including the following important one.

---

[2] for $\delta = 8$, add the following assumption: (iii) for $S$ the equisingularity class of an ordinary quadruple point, the analytic type of $D_y$ at its quadruple point is non-constant along $Y(S)$.

**(1.7) Theorem** (Vainsencher, Kleiman–Piene). *Let $V$ be a general threefold in $\mathbf{P}^4$ of degree $m$. If $m \geq 4$, then $V$ contains precisely the following number of 6-nodal plane curves of degree $m$:*

$$(1.7.1) \quad \frac{5}{6!}m\left(m^{17} - 12m^{16} + 24m^{15} + 155m^{14} - 405m^{13} + 1082m^{12} - 18469m^{11} + 66446m^{10}\right.$$
$$- 192307m^9 + 1242535m^8 - 4049006m^7 + 11129818m^6 - 53664614m^5$$
$$\left. + 166756120m^4 - 415820104m^3 + 1293514896m^2 - 2517392160m + 1781049600\right).$$

For $m = 5$, Vainsencher and Kleiman–Piene have refined this formula by computing the number of irreducible 6-nodal plane quintics in $V$ (the refinement resides in the irreducibility), see Theorem (3.5). This computation has been important for the understanding of the enumeration of rational curves in Calabi–Yau threefolds, see [XI, Section **].

In the situation of Theorem (1.7), the family of surfaces $\pi : F \to Y$ to be considered is the universal family of planes in $\mathbf{P}^4$ over the Grassmannian $\mathbf{G}(2,4)$, and the relative divisor is cut out by the hypersurface $V$: for all $y = [\Pi] \in Y = \mathrm{Gr}(3,5)$, $F_y$ is the plane $\Pi \subseteq \mathbf{P}^4$ itself, and $D_y$ is the plane curve $V \cap \Pi$.

**(1.8)** One may of course expect that formulae analogous to (1.3.1) and (1.6.1) hold for all $\delta$. The existence of formulae analogous to (1.3.1) without the requirement that they be structured in terms of Bell polynomials, i.e., the existence of polynomials $T_\delta$ in the Chern classes $d, k, s, x$ which enumerate (under suitable ampleness assumptions) $\delta$-nodal curves in a complete linear system on a fixed surface was the object of the Göttsche conjecture. The polynomials $T_\delta$ are called node polynomials, or universal polynomials. The existence of formulae analogous to (1.6.1) without the Bell polynomials structure, i.e., the existence of universal polynomials in the Chern classes $v, w_1, w_2$, is a family version of the Göttsche conjecture. The Göttsche conjecture has been proved by Tzeng [37], and another proof by Kool–Shende–Thomas [26] has appeared shortly after. The family version of the conjecture has been proved by Laarakker [27] (then a student of Kool), by elaborating on the proof of [26]. We report on the Göttsche conjecture in [XIV]; in this text, we focus rather on giving explicit formulae.

It had been proved by Göttsche that if his conjecture holds true, then the node polynomials $T_\delta$ satisfy some structure conditions encoded in the so-called Göttsche–Yau–Zaslow formula, see [XIV, Theorem **]. Qviller [30] (then a student of Piene) has shown that this formula implies the Bell polynomial structure of the node polynomials $T_\delta$ as in Formula (1.3.1). Once established the existence of formulae analogous to (1.3.1) for all $\delta$, one may compute the polynomials $a_i(d, k, s, x)$ by interpolation from known cases. In fact, in [30, Algorithm 2.1] Qviller explains how to extract the $a_i$ from the Göttsche–Yau–Zaslow formula. He gives them up to $i = 15$ in [31, Table 2] (only up to $i = 11$ in [30]), which we reproduce below; following his notation, we give $\tilde{a}_i = \frac{1}{(i-1)!}a_i$:

$\tilde{a}_9 = 26842726680d + 52612204910k + 26239943207s - 2617350984x$ ;

$\tilde{a}_{10} = -513240952752d - 1055936555124k - 556487181661s + 62064807888x$ ;

$\tilde{a}_{11} = 9861407170992d + 21186861410508k + 11720114258490s - 1410986931936x$ ;

$\tilde{a}_{12} = -190244562607008d - 425029422316200k - 245491696730341s + 31230909182592x$ ;

$\tilde{a}_{13} = 3682665360521280d + 8525631885908256k + 5119580760611226s - 678769122880224x$ ;

$\tilde{a}_{14} = -71494333556133600d - 171005998538392560k - 106382292871378404s - 14560213534363728x$ ;

$\tilde{a}_{15} = 1391450779290676680d + 3429957097334083248k + 2203960837196658328s - 309288199242633956x$.

In the particular case of the projective plane, the Göttsche conjecture implies the existence of univariate polynomials $R_\delta$ such that for large enough $d$, $R_\delta(d)$ is the number of $\delta$-nodal, degree $d$ plane curves passing through the appropriate number of fixed general points: this was known as the Di Francesco–Itzykson conjecture, and has been proved by Fomin and Mikhalkin [12] using floor diagrams. Building on their work, Block [5] has developed an algorithm to compute the polynomials $R_\delta$; he gives them for all $\delta \leqslant 14$. In [6] he has extended these results to the relative setting, i.e., to the enumeration of plane curves with prescribed tangency pattern with a fixed line.

**(1.9)** The last section presents a brief account on Kazarian's Thom polynomials, see [18] (anecdotally, my last words during the final session of our seminar in Tor Vergata were for this very topic). In a differential framework, they give (explicit!) formulae enumerating curves with virtually any kind of singularities, thus providing a vast generalization of the formulae discussed in detail in this text. In particular, one can use Kazarian's Thom polynomials to recover virtually all the formulae for surfaces in $\mathbf{P}^3$ that have been given in [XII]. On the latter application, see also [33] by Sasajima and Ohmoto, which deals with threefolds in $\mathbf{P}^4$ as well; for considerations in the same circle of ideas by the same authors, see [34].

We also review applications of Kazarian's Thom polynomials to the enumeration of linear spaces with prescribed tangency pattern with a hypersurface.

**(1.10)** It is natural to wonder why the universal enumerative formulae for $\delta$-nodal curves may be structured in terms of Bell polynomials (besides their miraculous fitness for the recursive computation in Paragraph (2.14)). An explanation has been provided by Qviller [31], which resides in a computation by inclusion-exclusion with respect to the lattice of polydiagonals in some appropriate fiber product, and the fact that the coefficients of the Bell polynomials may be expressed in terms of numbers of partitions. We shall see this in Section 6.

**(1.11)** The Göttsche universal polynomials are enumerative as soon as the line bundle $L$ is $\delta$-very ample as is proved in [26], see [XIV, Theorem (1.2)]. It follows that so are the formulae analogous to (1.3.1) for all $\delta$; in particular, Theorem (1.3) holds under the weaker assumption that $L$ is $\delta$-very ample. See also [24] for even better validity conditions in case $S$ is a rational surface; when $S$ is the projective plane, this is the result described in Paragraph (5.11).

I have nevertheless chosen to give the original proof by Kleiman and Piene of the validity of Formula (1.3.1) for $\delta \leqslant 8$ and for $3\delta$-very ample line bundles, as I reckon it displays interesting and instructive techniques. The workaround to avoid Kleiman and Piene's arguments, which indeed necessitate the full assumptions of Theorem (1.3), would be to extract $a_i$ for $i = 1, \ldots, 8$ from the Göttsche–Yau–Zaslow formulae, and then observe that they coincide with the $a_i$ given in (1.2); then Theorem (1.3) with the weaker assumptions is a mere corollary of [XIV, Theorem (1.2)]. Note that while $\delta$-ampleness gives sufficient control on the codimensions in $|L|$ of the loci of bad curves (i.e., curves with singularities worse than nodes) in order to prove [XIV, Theorem (2.1)], it does not give the full condition (ii) of Theorem (1.6).

The organization of the text is as follows. We directly give the proof in Section 2 of the main Theorem (1.6); this uses the theory of equisingularity, although it is only treated in Section 4. In Section 3 we do the computations giving Formulae (1.3.1) and (1.7.1), leaving aside the question of their enumerativity; the latter question is studied in Section 5, where we thus complete the proof of Theorem (1.3). Section 6 is devoted to Qviller's study of the Bell polynomials structure of formulae (1.3.1) and (1.6.1), and Section 7 to Kazarian's Thom polynomials. There is an appendix providing the basic results on Bell polynomials.

**Thanks.** I am grateful to Ragni Piene for her observations on a preliminary version of this text.

## 2 – The inductive procedure

This section is dedicated to the proof of Theorem (1.6), following Kleiman and Piene [23]; we discuss the relation with Vainsencher's ideas [38] in the last Section 2.4.

The idea is the following. Consider the situation set up in (1.4). Using sheaves of principal parts, also known as jet bundles, it is possible to compute the rational equivalence in $F$ of the locus $X_2$ along which the relative divisor $D/Y$ has a double point; then, pushing the class $[X_2]$ down to $Y$ should give the class $[U(D, 1)]$ of the (closure of the) locus in $Y$ of points $y$ such that the divisor $D_y$ has one node. To add one more node, one considers the pull-back to $X_2$ of the family $\pi : F \to Y$, and modify it so that the surface over the point $x \in X_2$ is the blow-up of $F_{\pi(x)}$ at $x$; then one considers the proper transform of the relative divisor pulled-back from $D$, and repeat the procedure.

In section 2.1 we explain the construction of the pairs $(F_i/X_i, D_i)$ induced by the pair $(F/Y, D)$ we start with; this gives the general setup for the proof of Theorem (1.6). In section 2.2 we state and prove the recursive formula which is the engine powering Theorem (1.6). After that there only remains a (substantial) computation to be made in order to obtain Theorem (1.6): we do it in section 2.3.

### 2.1 – The induced pairs

**(2.1)** Let, as in (1.4), $\pi : F \to Y$ be a smooth projective family of surfaces, where $Y$ is equidimensional and Cohen–Macaulay, and let $D$ be a relative effective divisor on $F/Y$. Consider the product $F \times_Y F$, and denote by $\mathrm{pr}_1, \mathrm{pr}_2 : F \times_Y F \to F$ the two projections, by $\Delta \subseteq F \times_Y F$ the diagonal, and by $\mathcal{I}_\Delta$ its ideal. Then, for all non-negative integers $i$, the $i$-th *bundle of principle parts* of $D$ is
$$\mathcal{P}^i_{F/Y}(D) = \mathrm{pr}_{2*}\bigl(\mathrm{pr}_1^* \mathcal{O}_F(D)/\mathcal{I}_\Delta^{i+1}\bigr).$$
The fibre of $\mathcal{P}^i_{F/Y}(D)$ over a point $x \in F$ identifies with the space of global sections of the restriction of $\mathcal{O}_F(D)$ to the $i$-th order infinitesimal neighbourhood of $x$ in $F$. You may consult [9, Section 7.2] for more details.

The relative divisor $D$ is defined by a global section $s$ of $\mathcal{O}_F(D)$, and the latter induces a section $s_i$ of $\mathcal{P}^i_{F/Y}(D)$ for all $i$. For all $i \geqslant 1$, let $X_i = \mathrm{V}(s_{i-1}) \subseteq F$ be the scheme of zeroes of the section $s_{i-1}$. As a set, $X_i$ is the locus of points $x \in F$ such that the divisor $D_{\pi(x)}$ has a point of multiplicity at least $i$ at $x$. In particular, $X_1 = D$.

**(2.2)** Next, let $\varepsilon \colon F' \to F \times_Y F$ be the blow-up along the diagonal $\Delta$, and denote its exceptional divisor by $E$. Let $\pi'$ be the composed map $\mathrm{pr}_2\varepsilon$; then $\pi' \colon F' \to F$ is another smooth family of surfaces, with fibre over a point $x \in F$ the blow-up of the surface $F_{\pi(x)}$ at the point $x$.

For all $i \geqslant 2$, set $F_i := \pi'^{-1}(X_i)$, and denote by $\pi_i \colon F_i \to X_i$ the restriction of $\pi'$. Eventually, let $\rho_i$ and $\rho'_i$ be the closed immersions $X_i \hookrightarrow F$ and $F_i \hookrightarrow F'$ respectively. This gives the following commutative diagram.

(2.2.1)
$$\begin{array}{ccccccc}
F & \xleftarrow{\mathrm{pr}_1} & F \times_Y F & \xleftarrow{\varepsilon} & F' & \xhookleftarrow{\rho'_i} & F_i \\
{\scriptstyle\pi}\downarrow & & {\scriptstyle\mathrm{pr}_2}\downarrow & & {\scriptstyle\pi'}\downarrow & & \downarrow{\scriptstyle\pi_i} \\
Y & \xleftarrow{\pi} & F & = & F & \xhookleftarrow{\rho_i} & X_i
\end{array}$$

By the definition of $X_i$, the pull-back of $D$ to $F_i$ by the composed map $\rho'_i \varepsilon \operatorname{pr}_1$ contains the exceptional divisor $\rho'^*_i E$ with multiplicity $i$; we set

$$D_i = (\rho'_i \varepsilon \operatorname{pr}_1)^* D - i \cdot \rho'^*_i E.$$

It is a relative effective divisor on $F_i/X_i$. See [23, Lemma 2.2] for more information.

**(2.3)** The sheaf $\mathcal{P}^{i-1}_{F/Y}(D)$ fits into an exact sequence

(2.3.1) $$0 \longrightarrow \mathcal{O}_F(D) \otimes \operatorname{Sym}^{i-1}(\Omega^1_{F/Y}) \longrightarrow \mathcal{P}^{i-1}_{F/Y}(D) \longrightarrow \mathcal{P}^{i-2}_{F/Y}(D) \longrightarrow 0,$$

where the first term is the symmetric power of the sheaf of relative differentials twisted by $\mathcal{O}_F(D)$, see, e.g., [9, Theorem 7.2]. This implies in particular that $\mathcal{P}^{i-1}_{F/Y}(D)$ is locally free of rank $\binom{i+1}{2}$. Therefore, at each scheme point $x \in X_i$, the codimension of $X_i$ in $F$ is at most $\binom{i+1}{2}$, and if equality holds then $X_i$ is a local complete intersection in $F$ at $x$, hence equidimensional and Cohen–Macaulay, so that the assumptions of (1.4) hold for the new family $\pi_i : F_i \to X_i$.

Moreover, if $X_i$ has the expected codimension $\binom{i+1}{2}$ in $F$, then its class $[X_i]$ equals the top Chern class of the vector bundle $\mathcal{P}^{i-1}_{F/Y}(D)$. It is then an exercise in Chern classes computations to express $[X_i]$ in terms of

(2.3.2) $$v = c_1(\mathcal{O}_F(D)) \quad \text{and} \quad w_j = c_j(\Omega^1_{F/Y}) \text{ for } j = 1, 2,$$

using the above exact sequence for all $i' \leqslant i$ (see [9, Section 7.3.2] for details). We shall use (compare with (1.4)):

(2.3.3)
$$[X_2] = v^3 + w_1 v^2 + w_2 v;$$
$$[X_3] = v^6 + 4w_1 v^5 + (5w_1^2 + 5w_2)v^4 + (2w_1^3 + 11w_1 w_2)v^3 + (6w_2 w_1^2 + 4w_2^2)v^2 + 4vw_1 w_2^2;$$
$$[X_4] = v^{10} + 10w_1 v^9 + (15w_2 + 40w_1^2)v^8 + (82w_1^3 + 111w_1 w_2)v^7$$
$$+ (91w_1^4 + 315w_2 w_1^2 + 63w_2^2)v^6 + (52w_1^5 + 429w_2 w_1^3 + 324w_1 w_2^2)v^5$$
$$+ (12w_1^6 + 282w_2 w_1^4 + 593w_2^2 w_1^2 + 85w_2^3)v^4 + (72w_2 w_1^5 + 464w_2^2 w_1^3 + 259w_1 w_2^3)v^3$$
$$+ (132w_2^2 w_1^4 + 246w_2^3 w_1^2 + 36w_2^4)v^2 + (72w_2^3 w_1^3 + 36w_1 w_2^4)v.$$

**(2.4)** We will want to apply the above formulae to the pairs $(F_i/X_i, D_i)$. Let $e = [E]$ be the class of the exceptional divisor on $F'$. By a fairly common and convenient abuse of notation, we let $e, v, w_1, w_2$ also denote their own pull-backs by whichever map is appropriate. We have the following relation on $F'$:

(2.4.1) $$e^3 + w_1 e^2 + w_2 e = 0,$$

see [9, Theorem 13.14] or [13, Remark 3.2.4, p. 55].

By definition of $D_i$,

(2.4.2) $$c_1(\mathcal{O}_{F_i}(D_i)) = v - ie.$$

**(2.5) Lemma.** *The relative Chern classes of $F'/F$ are*

(2.5.1) $$c_1(\Omega^1_{F'/F}) = w_1 + e \quad \text{and} \quad c_2(\Omega^1_{F'/F}) = w_2 - e^2.$$

Pulling-back by $\rho_i$ and $\rho'_i$, we obtain

$$c_1(\Omega^1_{F_i/X_i}) = w_1 + e \quad \text{and} \quad c_2(\Omega^1_{F_i/X_i}) = w_2 - e^2$$

as well.

*Proof.* For this proof it is necessary to indicate explicitly some of the pull-backs, so I will indicate those which I feel are relevant. Besides, I will drop the superscript in '$\Omega$' to improve readability. The map $\pi' : F' \to F$ is $\mathrm{pr}_2\varepsilon$, so we have an exact sequence

$$0 \to \mathrm{pr}_2^*\Omega_F \to \Omega_{F'} \to \Omega_{F'/F} \to 0.$$

Notice for instance that I dropped the $\varepsilon^*$ in the leftmost term of the sequence; also, $\Omega_F$ and $\Omega_{F'}$ should be understood as $\Omega_{F/Y}$ and $\Omega_{F'/Y}$. From this sequence we get that

(2.5.2) $\qquad c_1(\Omega_{F'/F}) = c_1(\Omega_{F'}) - \mathrm{pr}_2^*c_1(\Omega_F)$

(2.5.3) $\qquad c_2(\Omega_{F'/F}) = c_2(\Omega_{F'}) - \mathrm{pr}_2^*c_2(\Omega_F) - \mathrm{pr}_2^*c_1(\Omega_F) \cdot c_1(\Omega_{F'/F}).$

On the other hand, from [13, Example 15.4.3], we have

(2.5.4) $\qquad c_1(\Omega_{F'}) = \varepsilon^*c_1(\Omega_{F^2}) + e$

(2.5.5) $\qquad c_2(\Omega_{F'}) = \varepsilon^*c_2(\Omega_{F^2}) + j_*\varepsilon_E^*c_1(\Omega_\Delta) - e^2,$

where $F^2$ stands for $F \times_Y F$, $\varepsilon_E : E \to \Delta$ is the restriction of $\varepsilon : F' \to F$, and $j$ is the inclusion map $E \hookrightarrow F'$. There is an isomorphism $\Omega_{F^2} \cong \mathrm{pr}_1^*\Omega_F \oplus \mathrm{pr}_2^*\Omega_F$, hence

(2.5.6) $\qquad c_1(\Omega_{F^2}) = \mathrm{pr}_1^*c_1(\Omega_F) + \mathrm{pr}_2^*c_1(\Omega_F)$

(2.5.7) $\qquad c_2(\Omega_{F^2}) = \mathrm{pr}_1^*c_2(\Omega_F) + \mathrm{pr}_2^*c_2(\Omega_F) + \mathrm{pr}_1^*c_1(\Omega_F) \cdot \mathrm{pr}_2^*c_1(\Omega_F).$

From (2.5.2), (2.5.4), and (2.5.6), one gets

(2.5.8) $\qquad c_1(\Omega_{F'/F}) = \mathrm{pr}_1^*c_2(\Omega_F) + e,$

which is the first formula we wanted to prove. Plugging it in (2.5.3), we get

(2.5.9) $\qquad c_2(\Omega_{F'/F}) = c_2(\Omega_{F'}) - \mathrm{pr}_2^*c_2(\Omega_F) - \mathrm{pr}_2^*c_1(\Omega_F) \cdot \mathrm{pr}_1^*c_2(\Omega_F) - \mathrm{pr}_2^*c_1(\Omega_F) \cdot e,$

hence, with (2.5.5) and (2.5.7),

$$\begin{aligned}c_2(\Omega_{F'/F}) &= \mathrm{pr}_1^*c_2(\Omega_F) + \mathrm{pr}_2^*c_2(\Omega_F) + \mathrm{pr}_1^*c_1(\Omega_F) \cdot \mathrm{pr}_2^*c_1(\Omega_F) \\ &\quad + j_*\varepsilon_E^*c_1(\Omega_\Delta) - e^2 \\ &\quad - \mathrm{pr}_2^*c_2(\Omega_F) - \mathrm{pr}_2^*c_1(\Omega_F) \cdot \mathrm{pr}_1^*c_2(\Omega_F) - \mathrm{pr}_2^*c_1(\Omega_F) \cdot e \\ &= \mathrm{pr}_1^*c_2(\Omega_F) - e^2 + j_*\varepsilon_E^*c_1(\Omega_\Delta) - \mathrm{pr}_2^*c_1(\Omega_F) \cdot e.\end{aligned}$$

Since $\Delta$ is the diagonal, the second projection $\mathrm{pr}_2$ realizes an isomorphism $\Delta \simeq F$, and thus the two terms $j_*\varepsilon_E^*c_1(\Omega_\Delta)$ and $\mathrm{pr}_2^*c_1(\Omega_F) \cdot e$ are equal (note that there is an omitted $\varepsilon^*$ in the latter term); this gives the second formula we wanted to prove. $\square$

## 2.2 – The recursive formula

**(2.6) Definition.** *Let $\pi : F \to Y$ be a smooth projective family of surfaces with a relative effective divisor $D$ as in (1.4). Let $\delta$ be a non-negative integer. If assumption (ii) of Theorem (1.6) holds, then we define $U(D,\delta)$ to be the closure of $Y(\delta\mathrm{A}_1)$ in $Y$, i.e., $U(D,\delta)$ is the closure of the locus of points $y \in Y$ such that the singularities of the curve $D_y$ are exactly $\delta$ ordinary nodes. Note that $U(D,0) = Y$. For negative $\delta$, we set $U(D,\delta) = 0$ by convention.*

If assumption (ii) of Theorem (1.6) does not hold, this definition is not the same as [23, Definition 5.1], but we won't need to consider this situation. Beware moreover that, as Kleiman and Piene point out [23, p. 14], the definition of $U(D,\delta)$ they give in [23, Definition 5.1] differs from that given in [20, Section 4]; the difference lies in the way non-reduced curves $D_y$ are considered, see, e.g., [23, Display (5.3)] and the discussion wrapped around it.

**(2.7)** The idea of the proof of Theorem (1.6) is to relate $U(D, \delta)$ with $U(D_2, \delta - 1)$, where $D_2$ is the relative divisor of the family $F_2/X_2$ defined in Section 2.1. As we will see, the pair $(F_2/X_2, D_2)$ inherits good properties from assumption (ii) for $(F/Y, D)$. Still, there are components of $U(D_2, \delta - 1)$ that do not come from $U(D, \delta)$ and therefore should be discarded.

Let us consider the case $\delta = 4$ to understand the issue: $U(D, 4)$ is the (closure of the) family of 4-nodal curves, and has pure codimension 4 in $Y$ by assumption (ii). In turn, $X_2$ is the locus in $F$ where $D$ has a point of multiplicity at least 2, and $D_2$ is $D - 2E$ (with our convention that pull-backs are implicit) in the family $F_2/X_2$ of surfaces blown-up at a point of $X_2$. For $y \in Y(4\mathrm{A}_1)$, the curve $D_y$ gives rises to four 3-nodal curves in the family $D_2$, one over each of the four points in $X_2$ corresponding to the four nodes of $D_y$. This gives a locus of pure codimension 3 in $X_2$ of 3-nodal curves, which is the one we care about.

It turns out however that there is another locus of pure codimension 3 of 3-nodal curves in $X_2$. Consider indeed the locus $Y(\mathrm{D}_4)$[3] of curves with one ordinary triple point in $Y$; it has pure codimension 4 by assumption (ii) as well; surely, a curve with one ordinary triple point and no further singularity is not a limit of 4-nodal curves in $Y$. For $y \in Y(\mathrm{D}_4)$, the curve $D_y$ gives rise to one 3-nodal curve in the family $D_2$, namely the fibre of $D_2$ over the triple point $x \in F_y$ of $D_y$, which belongs to $X_2$; since $D_2 = D - 2E$, the curve $(D_2)_y$ is the proper transform of $D_y$ in the blow-up of $D_y$ at $x$ plus the exceptional divisor $E_y$ with multiplicity 1. Since $Y(\mathrm{D}_4)$ has pure codimension 4 in $Y$ and is entirely contained in $X_2$, the corresponding component of $U(D_2, 3)$ has pure codimension 3 in $X_2$.

An illustrative way of phrasing this, although arguably unconventional, is that ordinary triple points may be considered as double points with three infinitely near double points.

In Proposition (2.9) below we use the inclusion-exlucsion principle in order to take this phenomenon into account, and come up with a recursive formula for the class $[U(D, \delta)]$. Before we can state it we need to discuss the properties inherited by the pairs $(F_i/X_i, D_i)$ from assumption (ii) of Theorem (1.6) on the initial pair $(F/Y, D)$.

Given $\delta \geqslant 1$ and $i \geqslant 2$, we set

$$(2.7.1) \qquad \delta_i = \delta - \binom{i+1}{2} + 2.$$

Note that $\binom{i+1}{2} - 2$ is the expected codimension for the equisingular locus in $Y$ of fibres of $D$ with an ordinary $i$-fold point (this is a standard dimension count, which is generalized in Proposition (4.23)).

**(2.8) Proposition.** *Let $\pi : F \to Y$ be a smooth projective family of surfaces with a relative effective divisor $D$ as in (1.4). Let $\delta$ be a non-negative integer, and assume that condition (ii) of Theorem (1.6) holds. Then, for all $i \geqslant 2$:*
*(i) condition (ii) of Theorem (1.6) holds for the pair $(F_i/X_i, D_i)$ in codimension $\delta_i$, i.e., for all multi-equisingularity class $\mathsf{S}$, $\mathrm{codim}_{X_i}\bigl(X_i(\mathsf{S})\bigr) \geqslant \min\bigl(\delta_i + 1, \mathrm{expcod}(\mathsf{S})\bigr);$*
*(ii) for all $s \leqslant \delta_i$, each component of $U(D_i, s)$ is naturally birational to a component of either the equisingular locus in $Y$ of curves with an ordinary $i$-fold point and $s$ nodes, or that of curves with an ordinary $(i+1)$-fold point and $s - i - 1$ nodes.*

In order to prove this proposition properly, it is best to have an accurate description of equisingularity classes in terms of arrangements of infinitely near multiple points. For expository reasons I have chosen to describe this later in the text (in Section 4), so that I could directly give the essential ideas underlying Theorem (1.6), in the hope that the reader would have an intuitive understanding of these matters. Persevering with this hope, I will prove the proposition now, but feel free to read Section 4 first.

---

[3] we use roman capitals to denote singularity classes, so that for instance $\mathrm{D}_4$ does not get confused with the divisor $D_4$ on the surface $F_4/X_4$.

*Proof.* Let $\mathsf{S}'$ be an equisingularity class. The equisingular locus $X_i(\mathsf{S}')$ is the locus of points $x \in X_i$ such that the curve $(D_i)_x = (D - iE)_x$ in $F'_x$, the blow-up of $F_{\pi(x)}$ at $x$, has singularities of type $\mathsf{S}'$. In particular, the curve $(D_i)_x$ is reduced, hence it either does not contain the exceptional divisor $E_x$, or has it as a multiplicity 1 component. Thus in the blow-down $F_{\pi(x)}$, the curve $D_{\pi(x)}$ has multiplicity either $i$ or $i+1$ at $x$. Let $\mathsf{S}$ be its equisingularity class $\mathsf{S}$; it is determined by the multiplicity at $x$, the equisingularity class $\mathsf{S}'$, and the intersection pattern of $(D_i)_x$ with the exceptional divisor $E_x$. In particular the following inequality holds,

$$(2.8.1) \qquad \operatorname{expcod}(\mathsf{S}) \geqslant \operatorname{expcod}(\mathsf{S}') + \binom{i+1}{2} - 2,$$

see Section 4.4. Moreover, for any irreducible component $Z$ of $X_i(\mathsf{S}')$ containing $x$, there is a dense subset $U \subseteq Z$ such that $\pi(U) \subseteq Y(\mathsf{S})$. Then,

$$\operatorname{codim}_Y\bigl(\pi(U)\bigr) \geqslant \min\bigl(\delta + 1, \operatorname{expcod}(\mathsf{S})\bigr)$$

by condition (ii) of Theorem (1.6) for the pair $(F/Y, D)$. Since $F$ has dimension 2 over $Y$,

$$\operatorname{codim}_F(Z) = \operatorname{codim}_Y\bigl(\pi(U)\bigr) + 2.$$

In addition, $X_i$ has codimension at most $\binom{i+1}{2}$ in $F$, see (2.3), so eventually

$$\operatorname{codim}_{X_i}(Z) \geqslant \min\bigl(\delta + 1, \operatorname{expcod}(\mathsf{S})\bigr) + 2 - \binom{i+1}{2},$$

and, taking (2.8.1) into account, assertion (i) of the proposition follows.

To prove assertion (ii) we need to show that if $x$ is general in $X_i(\mathsf{S}')$ (or rather in the irreducible component $Z$), then the curve $D_{\pi(x)}$ in the blow-down has an ordinary singularity at $x$. Consider first the case when $(D_i)_x$ does not contain $E_x$; then $D_{\pi(x)}$ has a point of multiplicity $i$ at $x$, and it is ordinary if and only $(D_i)_x$ intersects $E_x$ transversely. If the singularity at $x$ were not ordinary, then we would have a strict inequality in (2.8.1), and thus $Z$ would give by projection by $\pi$ a family of too large dimension in $Y(\mathsf{S})$, given that (ii) of Theorem (1.6) holds for $(F/Y, D)$; note in particular that this condition implies that $X_i$ has codimension exactly $\binom{i+1}{2}$ in $F$ (we are now tacitly assuming $\delta_i \geqslant 0$, for otherwise $U(D_i, s) = 0$ by definition). The proof in the case when $(D_i)_x$ contains $E_x$ is similar: then $D_{\pi(x)}$ has a point of multiplicity $i+1$ at $x$, the inequality

$$\operatorname{expcod}(\mathsf{S}) \geqslant \operatorname{expcod}(\mathsf{S}') + \binom{i+2}{2} - 2$$

holds instead of (2.8.1), and it is strict if and only if the singularity at $x$ is non-ordinary. The conclusion follows as in the first case, noting that $X_{i+1}$ has codimension $\binom{i+2}{2}$ in $Y$. □

A much more careful version of assertion (ii) of the previous proposition is given in [23, Section 3]. We are now ready to state the recursive formula for $[U(D, \delta)]$.

**(2.9) Proposition.** *Consider the situation of Theorem (1.6). Then the following formula holds:*

$$(2.9.1) \qquad [U(D, \delta)] = \frac{1}{\delta}\, \pi_* \biggl( \sum_{i \geqslant 2} (-1)^i \rho_{i*} [U(D_i, \delta_i)] \biggr),$$

*where $\rho_i \colon X_i \to F$ denotes the inclusion, see Diagram (2.2.1), and the integers $\delta_i$ are as defined in (2.7.1).*

*Proof.* The proof is along the lines of the discussion in paragraph (2.7) above. There is a $\delta : 1$ covering of $Y(\delta A_1)$ in $X_2((\delta - 1)A_1) = X_2(\delta_2 A_1)$, corresponding to the fact that each $\delta$-nodal curve of the family $D/Y$ has $\delta$ different partial normalisations appearing as members of $D_2/X_2$. By Proposition (2.8), the components of $X_2(\delta_2 A_1)$ are of two kinds: those of the first kind sum up as the total space of the above mentioned $\delta : 1$ covering, whereas those of the second kind amount to curves of the family $D/Y$ with a triple point and $\delta_2 - 3 = \delta_3$ nodes. In sum,

$$[U(D_2, \delta_2)] = \sum_{\substack{Z \subseteq X_2(\delta_2 A_1) \\ \text{type I}}} [Z] + \sum_{\substack{Z \subseteq X_2(\delta_2 A_1) \\ \text{type II}}} [Z] \quad \text{and} \quad \pi_* \left( \sum_{\substack{Z \subseteq X_2(\delta_2 A_1) \\ \text{type I}}} [Z] \right) = \delta [U(D, \delta)].$$

In turn, again by Proposition (2.8), the components of $X_3(\delta_3 A_1)$ are of two kinds: those of the first kind are birational to the components of type II of $X_2(\delta_2 A_1)$, whereas those of the second kind amount to curves of the family $D/Y$ with a triple point and $\delta_3 - 4 = \delta_4$ nodes. In sum,

$$[U(D_3, \delta_3)] = \sum_{\substack{Z \subseteq X_3(\delta_3 A_1) \\ \text{type I}}} [Z] + \sum_{\substack{Z \subseteq X_3(\delta_3 A_1) \\ \text{type II}}} [Z] \quad \text{and} \quad \pi_* \left( \sum_{\substack{Z \subseteq X_3(\delta_3 A_1) \\ \text{type I}}} [Z] \right) = \pi_* \left( \sum_{\substack{Z \subseteq X_2(\delta_2 A_1) \\ \text{type II}}} [Z] \right).$$

Thus, continuing until $i$ is large enough so that $\delta_i < 0$, we find

$$\delta [U(D, \delta)] = \sum_{\substack{i \geq 2 \\ \delta_i \geq 0}} (-1)^i [U(D_i, \delta_i)],$$

which is the asserted formula. $\square$

**(2.10)** Both Propositions (2.8) and (2.9) hold without any limitation on $\delta$. The reason why one assumes $\delta \leq 8$ in Theorem (1.6) is that, although condition (ii) of this very Theorem (1.6) propagates indefinitely to the induced pairs $(F_i/X_i, D_i)$ with $\delta$ replaced by $\delta_i$ (this is assertion (i) of Proposition (2.8)), this is not the case for condition (i) for $\delta \geq 8$.[4]

To explain the issue, set $\delta = 8$ and assume conditions (i) and (ii) of Theorem (1.6) hold for $(F/Y, D)$. The equisingular stratum $Y(8A_1)$ has codimension 8 in $Y$, and the idea of the theorem is to consider the pair $(F_2/X_2, D_2)$ and the equisingular stratum $X_2(7A_1)$, which has codimension 7 in $X_2$, in order to understand $Y(8A_1)$. Now $X_4$ (the scheme of points of multiplicity at least 4 in the divisor $D$) is contained in $X_2$, by definition of the $X_i$'s, and it too has codimension 7 in $X_2$, as $X_2$ and $X_4$ have codimension 3 and 10 in $F$ respectively. Since $D_2$ is the pull-back of $D$ minus twice the exceptional divisor of $\varepsilon$, the fibres of $D_2$ over those points of $X_2$ that lie in $X_4$ contain the exceptional divisor with multicipity 2. Thus $X_4 \subseteq X_2(\infty)$, and condition (i) for $\delta_2 = \delta - 1 = 7$ does not hold for $(F_2/X_2, D_2)$.

This creates trouble because if we form the "$X_2$" of the pair $(F_2/X_2, D_2)$, it will contain the whole exceptional divisor over $X_4$, i.e., $\pi_2^{-1} X_4 \cap E$, which has codimension 8 in $F_2$, whereas the singular points of the fibres over $X_2(7A_1)$ have codimension 9 in $F_2$. The upshot is that $\pi_2^{-1} X_4 \cap E$ will give an excess contribution when one computes the class of $X_2(7A_1)$, and this has to be controlled if one is willing to prove Theorem (1.6) for $\delta = 8$.

Before we proceed, let me emphasize that the reason why Proposition (2.9) holds regardless of this issue is that, by definition, the general point in any component of the support of $U(D_i, \delta_i)$ corresponds to a curve of the family $D_i/X_i$ which is $\delta_i$-nodal, i.e., its singularities are exactly $\delta_i$ ordinary double points, hence non-reduced curves are completely disregarded in Formula (2.9.1).

---

[4] In fact the issue arises for all $\delta \geq 8$; for $\delta = 8$, Kleiman and Piene are able to overcome the issue, see Paragraph (2.15).

## 2.3 – The computation — Proof of the main theorem

In this section we compute the class $[U(D, \delta)]$ by induction using formula (2.9.1), which will prove Theorem (1.6). I do it in (2.14) below, but before that, let me start by working out in a pedestrian way the formulae for the very first values of $\delta$.

**(2.11) Example.** For $\delta = 1$ we have $[U(D, 1)] = \pi_*[X_2]$ thanks to assumptions (i) and (ii) in Theorem (1.6), and thus the formula for $[X_2]$ given in (2.3.3) gives us the formula we are looking for. Explicitly, we set

$$b_1(v, w_1, w_2) = x_2(v, w_1, w_2)$$
$$= v^3 + w_1 v^2 + w_2 v$$

(compare with the definition in Display (1.4.1)), and then

(2.11.1) $$[U(D, 1)] = \pi_*\big(b_1(v, w_1, w_2)\big).$$

Next, for $\delta = 2$, we get

$$[U(D, 2)] = \frac{1}{2} \pi_* \rho_{2*} [U(D_2, 1)]$$

from the recursive formula (2.9.1), with the the notation as in Diagram (2.2.1). Applying formula (2.11.1) to the pair $(F_2/X_2, D_2)$, and using (2.4.2) and (2.5.1) to get an expression for the classes in $F_2$ analogous to $v, w_1, w_2$ in $F$, we get

(2.11.2) $$[U(D_2, 1)] = \pi_{2*}\big(b_1(v - 2e, w_1 + e, w_2 - e^2)\big).$$

We will need to push this class forward to $F$ via $\rho_2$ to derive the formula for $[U(D, 2)]$ as $\pi_*$ of some class in $F$. We do it in (2.13) below using the following general formula.

**(2.12) Lemma.** *Let $R$ be a polynomial expression in $v, w_1, w_2$. Consider an integer $i \geqslant 2$ and classes $v, w_1, w_2, e$ as in (2.3.2) and (2.4). Then*

$$\pi_{i*} \rho_i'^* R(v - ie, w_1 + e, w_2 - e^2) = \rho_i^* \big[\pi^* \pi_* R(v, w_1, w_2) + Q_{R,i}(v, w_1, w_2)\big],$$

*in the notation of Diagram (2.2.1), and with $Q_{R,i}$ as defined in (1.4).*

*Proof.* By commutativity of Diagram (2.2.1), one has

$$\pi_{i*} \rho_i'^* R(v - ie, w_1 + e, w_2 - e^2) = \rho_i^* \pi'_* R(v - ie, w_1 + e, w_2 - e^2)$$
$$= \rho_i^* \operatorname{pr}_{2*} \varepsilon_* R(v - ie, w_1 + e, w_2 - e^2).$$

We shall then use the two relations

(2.12.1) $$\varepsilon_* e = 0 \quad \text{and} \quad \varepsilon_* e^2 = -[\Delta],$$

where $\Delta$ is the diagonal in $F \times_Y F$, see [13, Corollary 4.2.2] or [9, Proposition 13.12]. We expand the polynomial expression $R(v - ie, w_1 + e, w_2 - e^2)$ and reduce it modulo the relation (2.4.1); expanding a notation shorthand, we obtain

$$R(v - ie, w_1 + e, w_2 - e^2) = \varepsilon^* R(v, w_1, w_2) + e \cdot \varepsilon^* S_{R,i}(v, w_1, w_2) - e^2 \cdot \varepsilon^* Q_{R,i}(v, w_1, w_2)$$

for some polynomial $S_{R,i}$, and with $Q_{R,i}$ as defined in (1.4) (the minus sign is a feature, not a bug). Thus, with the relations (2.12.1), we get the equality of classes in $F \times_Y F$,

$$\varepsilon_* R(v - ie, w_1 + e, w_2 - e^2) = R(v, w_1, w_2) + Q_{R,i}(v, w_1, w_2) \cdot [\Delta].$$

In the above, $R(v, w_1, w_2)$ is a shorthand for $\mathrm{pr}_1^* R(v, w_1, w_2)$ so, since $\mathrm{pr}_{2*}\mathrm{pr}_1^* = \pi^*\pi_*$ by the commutativity of Diagram (2.2.1),

$$\mathrm{pr}_{2*}\varepsilon_* \mathrm{pr}_1^* R(v, w_1, w_2) = \pi^*\pi_* R(v, w_1, w_2).$$

On the other hand,

$$\mathrm{pr}_{2*}\bigl(\mathrm{pr}_1^* Q_{R,i}(v, w_1, w_2) \cdot [\Delta]\bigr) = Q_{R,i}(v, w_1, w_2)$$

since $\Delta$ is the diagonal. Summing up,

$$\mathrm{pr}_{2*}\varepsilon_* R(v - ie, w_1 + e, w_2 - e^2) = \pi^*\pi_* R(v, w_1, w_2) + Q_{R,i}(v, w_1, w_2),$$

which proves the lemma. $\square$

**(2.13) Example (2.11) continued.** Expanding a notation shorthand in (2.11.2), we have

$$\begin{aligned}[U(D_2, 1)] &= \pi_{2*}\rho_2'^{\,*}\bigl(b_1(v - 2e, w_1 + e, w_2 - e^2)\bigr) \\ &= \rho_2^*\bigl[\pi^*\pi_* b_1(v, w_1, w_2) + Q_{b_1, 2}(v, w_1, w_2)\bigr]\end{aligned}$$

as an application of Lemma (2.12). Then, as $\rho_2$ is the closed immersion $X_2 \hookrightarrow F$,

$$\begin{aligned}\rho_{2*}[U(D_2, 1)] &= \rho_{2*}\rho_2^*\bigl[\pi^*\pi_* b_1(v, w_1, w_2) + Q_{b_1, 2}(v, w_1, w_2)\bigr] \\ &= \bigl[\pi^*\pi_* b_1(v, w_1, w_2) + Q_{b_1, 2}(v, w_1, w_2)\bigr] \cdot [X_2],\end{aligned}$$

after what

$$\begin{aligned}\pi_*\rho_{2*}[U(D_2, 1)] &= \pi_*\Bigl(\bigl[\pi^*\pi_* b_1(v, w_1, w_2) + Q_{b_1, 2}(v, w_1, w_2)\bigr] \cdot [X_2]\Bigr) \\ &= \pi_* b_1(v, w_1, w_2) \cdot \pi_*[X_2] + \pi_*\bigl(Q_{b_1, 2}(v, w_1, w_2) \cdot [X_2]\bigr)\end{aligned}$$

where we used the projection formula for the first summand. Finally, recalling that

$$[X_2] = b_1(v, w_1, w_2) = x_2(v, w_1, w_2),$$

we get

$$\pi_*\rho_{2*}[U(D_2, 1)] = \pi_* b_1(v, w_1, w_2) \cdot \pi_* b_1(v, w_1, w_2) + \pi_*\bigl(Q_{b_1, 2}(v, w_1, w_2) \cdot x_2(v, w_1, w_2)\bigr);$$

we set

$$b_2(v, w_1, w_2) = Q_{b_1, 2}(v, w_1, w_2) \cdot x_2(v, w_1, w_2),$$

(compare with Display (1.4.1)), and then

$$\begin{aligned}\pi_*\rho_{2*}[U(D_2, 1)] &= (\pi_* b_1)^2 + \pi_* b_2 \\ &= P_2\left(\pi_* b_1, b_2\right),\end{aligned}$$

as we wanted to show.

**(2.14) Proof of Theorem (1.6) when $\delta \leqslant 7$.** The argument is by induction on $\delta$. For $\delta = 1$, assumptions (i) and (ii) imply that $X_2$ has pure codimension 3 in $F$, and the general point in any component of $X_2$ is the only singular point of a 1-nodal curve in the family $D$. Indeed, on the one hand $X_2$ has codimension at most 3 in $F$ by its description in Section 2.1, and on the other hand the locus of 1-nodal curves in $D$ has pure codimension 1 in $Y$, whereas all other equisingular loci as well as $Y(\infty)$ have codimension strictly greater; thus only the families of 1-nodal curves in $D$ are large enough for the corresponding singular points to fill up a dense open subset in a component of $X_2$, and these very components have codimension exactly 3 in $F$. It follows that $[U(D, 1)] = \pi_*[X_2]$, which proves the theorem for $\delta = 1$.

Assume next that $\delta \geqslant 2$ and that the theorem holds for all families verifying the hypotheses of the theorem with $\delta$ replaced by $\delta' < \delta$. The recursive formula (2.9.1) gives

(2.14.1) $$\delta[U(D, \delta)] = \pi_*\Big(\rho_{2*}[U(D_2, \delta - 1)] - \rho_{3*}[U(D_3, \delta - 4)]\Big)$$

since $\delta_i = 0$ for all $i > 3$, because $\delta \leqslant 7$. We want to apply the theorem by induction to the induced pairs $(F_i/X_i, D_i)$ for $i = 2, 3$, in order to compute $[U(D, \delta_i)]$. That condition (ii) in the theorem holds for these two pairs, for the integers $\delta_2$ and $\delta_3$ instead of $\delta$, is assertion (i) of Proposition (2.8). Condition (i) on the other hand is not granted a priori, and we will need to remove some negligible loci in $Y$ to make it hold. By the definition of $D_2$ and $D_3$, $X_2(\infty)$ and $X_3(\infty)$ are

$$\big(\pi^{-1}Y(\infty) \cap X_2\big) \cup X_4 \quad \text{and} \quad \big(\pi^{-1}Y(\infty) \cap X_3\big) \cup X_5,$$

respectively. By condition (i) for the inital pair $(F/Y, D)$, $Y(\infty)$ has codimension at least $\delta + 1$ in $Y$, hence we may replace $Y$ by $Y - Y(\infty)$. Moreover, since $\delta \leqslant 7$, by condition (ii) the locus $Y(\text{4-tuple})$ (over which the fibres of $D$ have a point of multiplicity at least 4) has codimension at least $\delta + 1$ as well, so we may further replace $Y$ by $Y - Y(\text{4-tuple})$. The upshot is that condition (i) holds for the pairs $(F_i/X_i, D_i)$, $i = 2, 3$, as well. Finally, note that $X_2$ and $X_3$ are local complete intersections in $F$, as they have the expected codimension everywhere; it follows that they are Cohen–Macaulay, since $F$ is, as it is smooth over the Cohen–Macaulay $Y$. Thus, by induction, formula (2.14.1) gives
(2.14.2)
$$\delta[U(D,\delta)] = \pi_*\Big(\rho_{2*}\frac{1}{(\delta-1)!}P_{\delta-1}\big(\pi_{2*}b_1(D_2),\dots\big) - \rho_{3*}\frac{1}{(\delta-4)!}P_{\delta-4}\big(\pi_{3*}b_1(D_3),\dots\big)\Big).$$

Then, by Lemma (2.12), and arguing as in paragraph (2.13), we get for $i = 2, 3$:

$$\begin{aligned}\rho_{i*}P_{\delta_i}\big(\pi_{i*}b_1(D_i),\dots\big) &= \rho_{i*}\rho_i^*\big[P_{\delta_i}\big(\pi^*\pi_*b_1(D) + Q_{b_1,i}(D),\,\dots\big)\big] \\ &= \big[P_{\delta_i}\big(\pi^*\pi_*b_1(D) + Q_{b_1,i}(D),\,\dots\big)\big] \cdot x_i(D).\end{aligned}$$

In turn, using identity (A.7.2) about Bell polynomials, this gives:

$$\rho_{i*}P_{\delta_i}\big(\pi_{i*}b_1(D_i),\dots\big) = \sum_{s=0}^{\delta_i}\binom{\delta_i}{s}P_s\big(\pi^*\pi_*b_1(D),\,\dots\big)P_{\delta_i-s}\big(Q_{b_1,i}(D),\,\dots\big) \cdot x_i(D),$$

hence, by the projection formula,

$$\pi_*\rho_{i*}P_{\delta_i}\big(\pi_{i*}b_1(D_i),\dots\big) = \sum_{s=0}^{\delta_i}\binom{\delta_i}{s}P_s\big(\pi_*b_1(D),\,\dots\big) \cdot \pi_*\Big(P_{\delta_i-s}\big(Q_{b_1,i}(D),\,\dots\big) \cdot x_i(D)\Big).$$

Thus, returning to Display (2.14.2),

$$
\begin{aligned}
\delta!\,[U(D,\delta)] &= \sum_{s=0}^{\delta-1}\binom{\delta-1}{s}P_s\bigl(\pi_*b_1(D),\dots\bigr)\cdot \pi_*\Bigl(P_{\delta-1-s}\bigl(Q_{b_1,2}(D),\dots\bigr)\cdot x_2(D)\Bigr) \\
&\quad -\frac{(\delta-1)!}{(\delta-4)!}\sum_{s=0}^{\delta-4}\binom{\delta-4}{s}P_s\bigl(\pi_*b_1(D),\dots\bigr)\cdot \pi_*\Bigl(P_{\delta-4-s}\bigl(Q_{b_1,3}(D),\dots\bigr)\cdot x_3(D)\Bigr) \\
&= \sum_{s=0}^{\delta-1}\binom{\delta-1}{s}P_s\bigl(\pi_*b_1(D),\dots\bigr)\cdot \pi_*\Bigl(P_{\delta-1-s}\bigl(Q_{b_1,2}(D),\dots\bigr)\cdot x_2(D) \\
&\quad -(\delta-1-s)(\delta-2-s)(\delta-3-s)P_{\delta-4-s}\bigl(Q_{b_1,3}(D),\dots\bigr)\cdot x_3(D)\Bigr) \\
&= \sum_{s=0}^{\delta-1}\binom{\delta-1}{s}P_s\bigl(\pi_*b_1(D),\dots\bigr)\cdot \pi_*b_{\delta-s}(D),
\end{aligned}
$$

where the last equality follows directly from the definition of $b_{\delta-s}(D)$, see Display (1.4.1). By identity (A.7.1) for Bell polynomials, this reads

$$\delta!\,[U(D,\delta)] = P_\delta\bigl(\pi_*b_1(D),\dots\bigr),$$

which is the formula we wanted to prove.  $\square$

**(2.15) Discussion of Theorem (1.6) when $\delta = 8$.** For $\delta \geqslant 8$ non-reduced curves inevitably enter the picture, and with them excess contributions, see Paragraph (2.10). For $\delta = 8$ Kleiman and Piene are able to make the appropriate correction, along the following lines. First apply Proposition (2.9) as before, to the effect that

$$8\,[U(D,8)] = \pi_*\Bigl(\rho_{2*}[U(D_2,7)] - \rho_{3*}[U(D_3,4)] + \rho_{4*}[U(D_4,0)]\Bigr).$$

For $i = 3, 4$, the induced pairs $(F_i/X_i, D_i)$ do satisfy the conditions of the theorem for $\delta = \delta_i$ (shrinking $Y$ appropriately), hence the two last summands may be computed using the theorem for $\delta' < 8$. Namely, we have

$$[U(D_3,4)] = \frac{1}{4!}P_4\bigl(\pi_{3*}b_1(D_3),\dots,\pi_{3*}b_4(D_3)\bigr) \quad \text{and} \quad [U(D_4,0)] = [X_4].$$

For $i = 2$ however, condition (i) does not hold as $X_2(\infty) = X_4$ has codimension 7 in $X_2$ as we have seen in Paragraph (2.10). Still, if we remove the troublesome $X_2(\infty)$, then Formula (1.6.1) does give the class of the locus of 7-nodal curves in $X_2 \setminus X_4$; the upshot is that

$$\frac{1}{7!}P_7\bigl(\pi_{2*}b_1(D_2),\dots,\pi_{2*}b_7(D_2)\bigr)$$

decomposes as the sum of $[U(D_2,7)]$ and a term supported on $X_4$.

Kleiman and Piene [23] show that in fact, provided the additional condition in Footnote 2 holds, the latter term supported on $X_4$ is actually $C.[X_4]$ for some integer $C$ independent of the family $(F/Y, D)$. It is then possible to compute the constant $C$ if one is able to enumerate the number of 8-nodal curves in one single example; Kleiman and Piene [21, Example 3.8] compute by hand the number of 8-nodal plane quintics passing through 12 general points (one could also have used the recursive formula of Caporaso and Harris), and thus find that $C = -3280$. The upshot is that

$$8\,[U(D,8)]\\ = \pi_*\Big(\rho_{2*}\Big(\frac{1}{7!}P_7\big(\pi_{2*}b_1(D_2),\,\ldots\big) + 3280\,x_4(D)\Big) - \rho_{3*}\frac{1}{4!}P_4\big(\pi_{3*}b_1(D_3),\,\ldots\big) + \rho_{4*}x_4(D)\Big),$$

from which Formula (1.6.1) follows by the same computations as before. □

The proof that the correction term is indeed $C.[X_4]$ for some integral constant $C$ is fairly elaborate, and I will not report on it. Please see [23, Sections 6–7].

## 2.4 – Comparison with Vainsencher's computation

The idea of Vainsencher's proof of his version of Theorem (1.6), see [38], is essentially equivalent to that of Kleiman and Piene which I have been describing above. Yet I believe it is interesting to compare the two, as it sheds some light on the recursive formula (2.9.1).

**(2.16)** Vainsencher makes a direct iterated use of the formulae for the class $[X_2]$: he first considers the $X_2$ of the pair $(F/Y, D)$ and the induced pair $(F_2/X_2, D_2)$ as in Section 2.1, then he considers "$X_2$" of $(F_2/X_2, D_2)$, $X_{2,2}$ say, with the corresponding induced pair $(F_{2,2}/X_{2,2}, D_{2,2})$, and so on.

He coined the kind of singularities one directly enumerates with this procedure, with the possibility of using any $X_i$ instead of $X_2$ at any stage: given a divisor $D$ in a smooth variety $V$, a *singularity sequence of type* $(m_1, \ldots, m_r)$ on $D$ is an ordered sequence of (possibly infinitely near) points $(x_1, \ldots, x_r)$ on $D$ such that:
  (i) $x_1$ is a point of multiplicity at least $m_1$ of $D$, and
 (ii) in the blow-up $\varepsilon_1 : V_1 \to V$ at $x_1$, with exceptional divisor $E_1$, the sequence $(x_2, \ldots, x_r)$ is a singularity sequence of type $(m_2, \ldots, m_r)$ of the divisor $\varepsilon_1^* D - m_1 E_1$.

Then, given a pair $(F/Y, D)$ of arbitrary relative dimension $n$, satisfying conditions analogous to those of Theorem (1.6), the locus in $Y$ representing members of the family $D$ with a singularity sequence of type $(m_1, \ldots, m_r)$ is closed in $Y$ of the expected codimension $\sum_{i=1}^{r} \binom{n+m_i-1}{n} - r\,n$, and it is the support of a natural nonnegative cycle, namely the push-forward by the composition of natural maps

$$X_{m_1,\ldots,m_r} \subseteq F_{m_1,\ldots,m_{r-1}} \longrightarrow \cdots \longrightarrow X_{m_1,m_2} \subseteq F_{m_1} \longrightarrow X_{m_1} \subseteq F \longrightarrow Y$$

of the fundamental class $[X_{m_1,\ldots,m_r}]$, which may be computed as in Section 2.1.

**(2.17)** Denote the type of singularity sequence $(m, \ldots, m)$ with $r$ repetitions by $(m^r)$ for short. Given a $\delta$-nodal curve $C$, each ordering of its nodes yields a singularity sequence of type $(2^\delta)$. However not all singularity sequences of type $(2^\delta)$ are of this kind: one may of course have any kind of double points instead of nodes, but these are not a problem as non-ordinary double points have expected codimension larger than that of ordinary double points; the non-nodal singularity sequences of type $(2^\delta)$ one needs to be careful about are those arising from a point of multiplicity larger than 2: for instance, an ordinary triple point has expected codimension 4 and gives rise to 3! singularity sequences of type $(2^4)$, as explained in Paragraph (2.7).

Thus, the enumeration of curves with a singularity sequence of type $(2^\delta)$ in a family $D$ as in Paragraph (2.16) above amounts to the enumeration of curves with any of the (finitely many) singularity types that give rise to singularity sequences of type $(2^\delta)$ and have expected codimension $\delta$, with weights corresponding to the possible orderings. Thus, after one has enumerated curves with a singularity sequence of type $(2^\delta)$, one is left with the task of isolating the contribution of $\delta$-nodal curves in this enumeration.

| $\delta$ | type of singularity | number of sequences |
|---|---|---|
| $\leqslant 3$ | $\delta A_1$ | $\delta!$ |
| 4 | $4A_1$ | $4!$ |
|   | $D_4$ | $3!$ |
| 5 | $5A_1$ | $5!$ |
|   | $D_4 + A_1$ | 30 |
| 6 | $6A_1$ | $6!$ |
|   | $D_4 + 2A_1$ | 180 |
|   | $D_6$ | 30 |
| 7 | $7A_1$ | $7!$ |
|   | $D_4 + 3A_1$ | 1260 |
|   | $D_6 + A_1$ | 210 |
|   | $E_7$ | 30 |

Table 1: Singularity types giving rise to singularity sequences of type $(2^r)$, and with expected codimension $r$

It is worth noting that the recursive formula (2.9.1) of Kleiman and Piene provides an automated process for selecting only the contribution of $\delta$-nodal curves among all curves with a singularity sequence of type $(2^\delta)$.

**(2.18)** For $\delta \leqslant 7$, the list of the singularity types giving rise to singularity sequences of type $(2^\delta)$ and of the appropriate expected codimesion $\delta$ has been worked out by Vainsencher [38] and Kleiman–Piene [20]. To do so it is convenient to use Enriques diagrams, which are presented in Section 4; this is the approach adopted in [20]. The list is given in Table 1, see Section 4 for the (standard) notation.

The cases $\delta \leqslant 4$ have already been discussed above (note that $D_4$ is the equisingularity class of ordinary triple points). For $\delta = 5$ nothing really new happens, but let me enumerate the number of sequences arising from an ordinary triple point and a node: one may either let the genuine node come first, and then we have 3! orderings for the rest as we know from the $\delta = 4$ case, or let the triple point come first, after what we have four nodes in the divisor "$\varepsilon^* D - 2E$" (three of wich at the intersection of the proper transform with the exceptional divisor), which gives 4! possible orderings; the total number of singularity sequences is therefore $3! + 4! = 30$. We leave it as a fun exercise to the reader to count the number of singularity sequences of type $(2^6)$ and $(2^7)$ arising from singularities $D_4 + 2A_1$ and $D_4 + 3A_1$, respectively.

The singularity $D_6$ is a triple point with three local branches, two of which are tangent, or in more classical terms, a triple point with an infinitely near double point (see Section 4 for more details). Thus, after blowing-up the ambient surface at such a singularity, when taking "$\varepsilon^* D - 2E$", one arrives at one ordinary double point and an ordinary triple point, both at the intersection of the proper transform with the exceptional divisor; this is a collection of singularities of type $D_4 + A_1$, hence the number 30 of singularity sequences of type $(2^6)$ arising from a $D_6$ singularity.

The singularity $E_7$ is a triple point with two local branches, one ordinary cusp and one smooth branch, with their tangent cones supported on the same line; in other words, it is a triple point with an infinitely near tacnode. After blowing-up and taking "$\varepsilon^* D - 2E$", one arrives at a singularity of type $D_6$, with the two tangent local branches given by the exceptional divisor and the proper transform of the cuspidal branch respectively. Hence the number of singularity sequences of type $(2^7)$ arising from a $E_7$ singularity is yet again 30. □

From the classification given in Table 1 one deduces, under the assumptions of Theorem (1.6), the following decompositions of the cycles $[U(D, 2^r)]$ enumerating curves with a singularity sequence of type $(2^r)$ in $D$, obtained as the push-forward of $[X_{2^r}]$ as in Paragraph (2.16) above: $[U(D, 2^\delta)] = \delta! [U(D, \delta A_1)]$ for all $\delta \leqslant 3$, and

(2.18.1)
$$\begin{aligned}
\left[U(D, 2^4)\right] &= 4! \left[U(D, 4A_1)\right] + 3! \left[U(D, D_4)\right] \\
\left[U(D, 2^5)\right] &= 5! \left[U(D, 5A_1)\right] + 30 \left[U(D, D_4 + A_1)\right] \\
\left[U(D, 2^6)\right] &= 6! \left[U(D, 6A_1)\right] + 180 \left[U(D, D_4 + 2A_1)\right] + 30 \left[U(D, D_6)\right] \\
\left[U(D, 2^7)\right] &= 7! \left[U(D, 7A_1)\right] + 1260 \left[U(D, D_4 + 3A_1)\right] + 210 \left[U(D, D_6 + A_1)\right] + 30 \left[U(D, E_7)\right].
\end{aligned}$$

In these formulae, for all (multi)equisingularity class $S$, $U(D, S)$ denotes the closure in $Y$ of the constructible set $Y(S)$ of all points $y \in Y$ such that the singularities of the curve $D_y$ are of the type $S$; thus, $U(D, \delta A_1)$ is what was denoted by $U(D, \delta)$ before.

**(2.19)** To conclude this section let me outline, following [20], how to compute, under the assumptions of Theorem (1.6) for $\delta \leqslant 7$, the non-nodal terms in the decompositions of Display (2.18.1) (it would be outrageous to call them parasitic), so that it is possible to obtain $[U(D, \delta)]$ from $[U(D, 2^\delta)]$.

First of all, for $\delta \leqslant 3$, it is straightforward to obtain $[U(D, \delta)]$ from $[U(D, 2^\delta)]$, hence to compute $[U(D, \delta)]$ from the procedure in Paragraph (2.16) (correspondingly, there is only one non-trivial summand in the recursive formula (2.9.1) in this case).

Then, it is possible to perform the same computations for the pair $(F_3/X_3, D_3)$; note that the conditions of Theorem (1.6) still hold for $\delta \leqslant 3$ for this pair. Moreover, thanks to these very conditions, the members of the family $D_3$ generically correspond to the members of the initial family $D$ with an ordinary triple point (i.e., with a singularity of type $D_4$: once again I apologize for this regrettable notation). Thus, it is possible to compute $[U(D, D_4 + \delta A_1)]$ for all $\delta = 0, 1, 2, 3$:
$$\left[U(D, D_4 + \delta A_1)\right] = \pi_* \left[U(D_3, \delta A_1)\right].$$

We will see now how to compute $[U(D, D_6)]$, $[U(D, D_6 + A_1)]$, and $[U(D, E_7)]$. When this is done we will know all cycle classes in the formulae of Display (2.18.1) except $[U(D, \delta A_1)]$ for $\delta = 4, 5, 6, 7$, so that it is possible to invert the system of linear relations to find them. The idea is once again to apply the formulae we know so far to an aptly chosen induced pair.

As we have seen in Paragraph (2.18) above, a singularity $D_6$ in $D$ gives a $D_4$ and a node in $D_2$, with the two latter singularity on the exceptional divisor. Of course it is also possible to obtain $D_4 + A_1$ in $D_2$ from a curve in $D$ with singularities of type $D_4 + 2A_1$, in two ways as one may blow-up at either one of the two nodes. In fact, the following relation holds:
$$\pi_* \left[U(D_2, D_4 + A_1)\right] = 2 \left[U(D, D_4 + 2A_1)\right] + \left[U(D, D_6)\right]$$

(see [20, p. 233] for more details). Since we already know $[U(D', \delta A_1)]$ and $[U(D', D_4 + \delta A_1)]$ for $\delta \leqslant 3$ for any family $D'$, we obtain a general formula for $[U(D, D_6)]$. In a similar fashion, we have the relation:
$$\pi_* \left[U(D_2, D_4 + 2A_1)\right] = 3 \left[U(D, D_4 + 3A_1)\right] + \left[U(D, D_6 + A_1)\right],$$

which gives a general formula for $[U(D, D_6 + A_1)]$.

Finally, since a singularity $E_7$ gives a $D_6$ in the family $D_2$, see Paragraph (2.18), we have the analogous relation:
$$\pi_*\big[U(D_2, D_6)\big] = \big[U(D, D_6 + A_1)\big] + \big[U(D, E_7)\big].$$

It gives a general formula for $[U(D, E_7)]$, and thus we have completed our program.

## 3 – Applications

In this section we spell out the computations for various applications of Theorem (1.6), essentially Theorems (1.3) and (1.7). The question of their validity is left aside for the moment, and will be considered later in Section 5, after the necessary material in equisingularity theory is treated.

### 3.1 – Curves in a linear system on a fixed surface

In this section I show how to get Formula (1.3.1) from (1.6.1). The setup is thus as in Theorem (1.3).

**(3.1)** Let $S$ be a smooth irreducible complex projective surface, and $L$ be a line bundle on $S$; set
$$d = L^2 \; ; \quad k = L \cdot K_S \; ; \quad s = K_S^2 \; ; \quad x = c_2(S).$$

In order to apply Theorem (1.6), we let $Y$ be the complete linear system $|L|$, and set $F = S \times Y$; we take $\pi : F \to Y$ and $p : F \to S$ to be the second and first projections, respectively. We let the relative divisor $D$ in $F/Y$ be the universal family of the linear system $|L|$.

Explicitly, let $s_0, \ldots, s_n$ be a basis of $H^0(S, L)$, and consider the corresponding homogeneous coordinates $(a_0 : \ldots : a_n)$ on the projective space $|L|$. Then $D$ is defined in $F = S \times Y$ by the equation
$$a_0 s_0 + \cdots + a_n s_n = 0.$$

Thus we see that $\mathcal{O}_F(D) = p^*L \otimes \pi^*\mathcal{O}_Y(1)$, hence
$$v = c_1\big(\mathcal{O}_F(D)\big) = p^*l + \pi^*h, \quad \text{where} \quad l = c_1(L) \quad \text{and} \quad h = c_1\big(\mathcal{O}_Y(1)\big).$$

Also, $\Omega^1_{F/Y}$ is simply the pull-back $p^*\Omega^1_S$, hence
$$w_j = c_j\big(\Omega^1_{F/Y}\big) = p^*c_j\big(\Omega^1_S\big) \quad \text{for} \quad j = 1, 2.$$

**(3.2) Lemma.** *Let $\alpha, \beta_1, \beta_2$ be non-negative integers such that $\alpha + \beta_1 + 2\beta_2 = \delta + 2$. Then:*
$$\pi_*\big(v^\alpha w_1^{\beta_1} w_2^{\beta_2}\big) = \begin{cases} \binom{\delta+2}{2} d.h^\delta & \text{if } (\alpha, \beta_1, \beta_2) = (\delta + 2, 0, 0) \\ (\delta + 1) k.h^\delta & \text{if } (\alpha, \beta_1, \beta_2) = (\delta + 1, 1, 0) \\ s.h^\delta & \text{if } (\alpha, \beta_1, \beta_2) = (\delta, 2, 0) \\ x.h^\delta & \text{if } (\alpha, \beta_1, \beta_2) = (\delta, 0, 1) \\ 0 & \text{otherwise,} \end{cases}$$

*with $d, k, s, x$ as in (3.1) above.*

*Proof.* For dimension reasons, one has $w_1^{\beta_1} = 0$ if $\beta_1 > 2$ and $w_2^{\beta_2} = 0$ if $\beta_2 > 1$; similarly, when expanding $v^\alpha = (p^*l + \pi^*h)^\alpha$, all terms involving $p^*l$ to a power striclty larger than 2 vanish. Moreover, for integers $\alpha', \alpha'', \beta_1, \beta_2$,

$$\pi_*\big(\pi^*h^{\alpha'} \cdot p^*(l^{\alpha''}w_1^{\beta_1}w_2^{\beta_2})\big) = h^{\alpha'} \cdot \pi_* p^*(l^{\alpha''}w_1^{\beta_1}w_2^{\beta_2})$$
$$= \begin{cases} l^{\alpha''}w_1^{\beta_1}w_2^{\beta_2}.h^{\alpha'} & \text{if } \alpha'' + \beta_1 + 2\beta_2 = 2 \\ 0 & \text{otherwise.} \end{cases}$$

Then the result follows from standard algebraic manipulations. □

Finally, applying the recipes in (1.4) with the above incarnations of $v, w_1, w_2$ to compute

$$a_i = \pi_* b_i(D, F/Y) \quad \text{for} \quad i = 1, \ldots, 8,$$

one gets the formulae in (1.2), and thus (1.3.1), as desired.

## 3.2 – Plane curves in a solid in $\mathbf{P}^4$

In this section I show how to get Formula (1.7.1) from (1.6.1), and then how to derive the number of irreducible 6-nodal plane quintics lying in a general quintic threefold.

**(3.3) Computations for Formula** (1.7.1)**.** Let $Y = \mathrm{Gr}(3,5)$ be the Grassmannian of 2-planes $\Pi \subseteq \mathbf{P}^4$, and $F$ be the corresponding universal family: thus

$$F = \mathbf{P}(\mathcal{Q}) \subseteq Y \times \mathbf{P}^4,$$

where $\mathcal{Q}$ is the tautological rank 3 quotient bundle of $\mathcal{O}_Y^5$. Let $\pi : F \to Y$ and $p : F \to \mathbf{P}^4$ be the two projections.

Let $V \subseteq \mathbf{P}^4$ be a general threefold of degree $m$, and define

$$D = F \cap (Y \times V).$$

If $m \geqslant 3$, then $V$ cannot contain a plane, and thus for all $y \in Y$, the fibre $D_y$, which identifies with the intersection $F_y \cap V$ in $\mathbf{P}^4$, is a degree $m$ plane curve; the upshot is that $D$ is a relative effective divisor on $F/Y$.

We may then apply Theorem (1.6) to derive the class in $Y$ of the locus of $\delta$-nodal curves in the family $D$. The Grassmannian $Y$ has dimension 6, so this computation makes sense only if $\delta \leqslant 6$. For all $\delta = 1, \ldots, 6$, Kleiman and Piene show in [21, Section 4] that assumptions (i) and (ii) of Theorem (1.6) hold in our situation if $m \geqslant 4$. In fact, this follows without too much difficulty from the verification of these very assumptions in the case of plane curves, i.e., the situation of Section 3.1 in the particular case when $S = \mathbf{P}^2$ and $L = \mathcal{O}_{\mathbf{P}^2}(m)$, carried out in [21, Section 3]. I will report on the latter analysis in Section 5, but skip the details of the adaptation to our current situation.

I will now sketch the computations needed to apply the recipes of (1.4) in the case under consideration, leaving aside the particulars of intersection theory on Grassmannians. Since the divisor $D$ is defined in $Y \times \mathbf{P}^4$ by the equation of $V \subseteq \mathbf{P}^4$, constant along $Y$, one has $\mathcal{O}_F(D) = p^*\mathcal{O}_{\mathbf{P}^4}(m)$, hence

$$v = c_1\big(\mathcal{O}_F(D)\big) = m.p^*h, \quad \text{where} \quad h = c_1\big(\mathcal{O}_{\mathbf{P}^4}(1)\big).$$

To compute the Chern classes $w_j = c_j(\Omega^1_{F/Y})$, one uses the Euler exact sequence

$$0 \longrightarrow \Omega^1_{F/Y}(1) \longrightarrow \pi^*\mathcal{Q} \longrightarrow \mathcal{O}_F(1) \longrightarrow 0,$$

and finds

$$w_1 = \pi^* q_1 - 3p^*h \quad \text{and} \quad w_2 = \pi^* q_2 - 2p^*h \cdot \pi^* q_1, \quad \text{where} \quad q_j = c_j(\mathcal{Q}) \quad \text{for} \quad j = 1, 2.$$

To compute with these classes note that, for dimension reasons, $h^j = 0$ for $j > 4$, and $\pi_* p^* h^j = 0$ for $j = 0, 1$; moreover,

$$\pi_* p^* h^2 = [Y], \qquad \pi_* p^* h^3 = q_1, \qquad \pi_* p^* h^4 = q_1^2 - q_2.$$

Then, it is algorithmic to derive the $a_i = \pi_* b_i(D, F/Y)$, $i = 1, \ldots, 6$, with the notation as in (1.4). After that, Formula (1.6.1) gives an expression of the classes $[U(D, \delta)]$, $\delta = 1, \ldots, 6$, in the Chow ring of the Grassmannian $Y = \mathrm{Gr}(3, 5)$, as polynomials in the integer $m$ and the classes $q_1$ and $q_2$, if $m \geqslant \delta/2 + 1$.

For $\delta = 6$, this class has degree $6 = \dim(Y)$, and one finds Formula (1.7.1) by applying the identities describing the intersection product in the Chow ring of the $Y = \mathrm{Gr}(3, 5)$, namely

$$q_1^6 = 5, \qquad q_1^4 q_2 = 3, \qquad q_1^2 q_2^2 = 2, \qquad q_2^3 = 1.$$

For lower values of $\delta$, one may use the expression of $[U(D, \delta)]$ to compute other enumerative numbers. Kleiman and Piene give the example of

$$[U(D, 3)] \cdot q_1^3$$
$$= \frac{m}{6} \left( 5\,m^8 - 30\,m^7 + 33\,m^6 + 23\,m^5 + 102\,m^4 + 359\,m^3 - 2330\,m^2 + 2048\,m + 240 \right),$$

which is the number of 3-nodal plane sections of $V$ whose span meets three given lines in general position in $\mathbf{P}^4$.

**(3.4) Example.** For $\delta = 6$ and $m = 4$, one finds that in a general quartic solid there are 5600 6-nodal plane quartic curves; each of those is the sum of four coplanar lines, since plane quartics have arithmetic genus 3.

**(3.5) Theorem** (Katz, Vainsencher, Kleiman–Piene)**.** *Let $V$ be a general quintic threefold in $\mathbf{P}^4$; it contains exactly:*
*(i)* 2875 *lines;*
*(ii)* 609 250 *smooth conics;*
*(iii)* 17 601 000 *irreducible* 6-*nodal plane quintics.*

For the first two numbers, see [17, p. 175] and, for more details on the number of lines, [9, Section 6.5]. I will now explain how to compute the third from the first two and Formula (1.7.1), following Vainsencher and Kleiman–Piene; again, I will skip the verification of the validity of the computation, for which I refer to [21, Proof of Theorem 4.3].

**(3.6)** Maintain the notation of Paragraph (3.3), and assume $m = 5$. Applying Formula (1.7.1) one finds the number of 6-nodal plane quintics in $V$, 21 617 125. Among them the reducible ones are either the sum of a line and a 2-nodal quartic curve, or the sum of a conic and cubic.

Curves of the latter type are the easier to enumerate. Let $C$ be a smooth conic contained in $V$. Then $C$ spans a plane, which intersects $V$ residually in a cubic $C'$; for general $V$ it is proven in [21] that for all smooth conics $C$ in $V$ the cubic $C'$ is smooth and intersects $C$ transversely. Thus the number of 6-nodal quintics sum of a conic and a cubic in $V$ equals the number of smooth conics in $V$, namely 609 250.

Now let $L$ be a line in $V$. Planes in $\mathbf{P}^4$ containing $L$ form a projective plane in the Grassmannian $Y = \mathrm{Gr}(3, 5)$, and each of those cuts out $V$ residually in a plane quartic. We will compute the number of 2-nodal curves in this family by applying Theorem (1.6) yet again (taking for granted that the quartics intersect $L$ transversely). Then we will need to subtract 2875 times this number from the number of 6-nodal plane quintics in $V$.

Denote by $Y_L \subseteq Y$ the variety of planes $\Pi \subseteq \mathbf{P}^4$ containing $L$, and set $F_L = F \times_Y Y_L$, so that $\pi_L = \pi|_{F_L} : F_L \to Y_L$ is the two-dimensional family of planes we want to consider. Then we set
$$D_L = F_L \cap (Y_L \times V) \quad \text{and} \quad D'_L = D_L - (Y_L \times L);$$
thus $D_L$ is the restriction to $F_L$ of the universal family $D$ of plane quintics in $V$, and $D'_L$ is the universal family of plane quartics residual to $L$. We want to apply Theorem (1.6) to the pair $(D'_L, F_L/Y_L)$.

The two Chern classes $w'_j = c_j(\Omega^1_{F_L/Y_L})$ are simply the restrictions of the two classes $w_j = c_j(\Omega^1_{F/Y})$, as $F_L = \pi^{-1}(F_L)$. Now let us compute the class of $D'_L$. The class of $D_L$ is the restriction of $v = c_1(\mathcal{O}_F(D)) = 5p^*h$ to $F_L$, and we want to subtract the class of the divisor $Y_L \times L$. The inclusion of this divisor in $F_L$ corresponds to a quotient map $\mathcal{Q}|_{Y_L} \twoheadrightarrow \mathcal{O}^{\oplus 2}_{Y_L}$, and the sheaf of ideals of $Y_L \times L$ in $F_L$ is $\pi_L^* \mathcal{K}_L(-1)$, where $\pi : F_L \to Y_L$ is the projection and $\mathcal{K}_L$ is the locally free $\mathcal{O}_{Y_L}$-module fitting in the exact sequence
$$0 \longrightarrow \mathcal{K}_L \longrightarrow \mathcal{Q}|_{Y_L} \longrightarrow \mathcal{O}^{\oplus 2}_{Y_L} \longrightarrow 0.$$
Thus,
$$v' = c_1\big(\mathcal{O}_{F_L}(Y_L \times L)\big) = -c_1\big(\pi_L^* \mathcal{K}_L(-1)\big) = (p^*h - \pi^*q_1)|_{F_L}.$$
Finally we compute, in the notation of Paragraph (1.4),
$$\begin{aligned}
a_i &= \pi_{L*} b_i (v', w'_1, w'_2) \\
&= \pi_{L*} b_i \big({\rho'_L}^*(p^*h - \pi^*q_1),\, {\rho'_L}^* w_1,\, {\rho'_L}^* w_2\big) \\
&= \rho_L^* \pi_* b_i \big(p^*h - \pi^*q_1,\, w_1,\, w_2\big) \\
&= [Y_L] \cdot \pi_* b_i \big(p^*h - \pi^*q_1,\, w_1,\, w_2\big),
\end{aligned}$$
where $\rho_L$ and $\rho'_L$ are the two closed immersions in the following Cartesian diagram.

$$\begin{array}{ccc} F_L & \xrightarrow{\rho'_L} & F \\ \pi_L \downarrow & & \downarrow \pi \\ Y_L & \xrightarrow{\rho_L} & Y \end{array}$$

Using the formula
$$[Y_L] = (q_1^2 - q_2)^2$$
for the class of the Schubert cell $Y_L$ in the Grassmannian $Y$, one completes the computation and finds that the number of 2-nodal quartics in the family $D'_L$ is 1185 by Formula (1.6.1).

The upshot is that the number of irreducible 6-nodal quintics in $V$ is
$$21617125 - 609250 - 2875 \times 1185 = 17601000$$
as asserted. $\square$

# 4 – Equisingularity classes

This section is devoted to the definition and study of equisingularity classes. The first two sections 4.1 and 4.2 contain preliminaries on infinitely near points and their configurations, while the main definitions (including those of equisingularity classes and Enriques diagrams) are given in section 4.3. In the last section 4.4 we discuss the expected codimension of the locus of curves belonging to a prescribed equisingularity class in a given family of curves.
In this section, unless explicit indication of the contrary, curves are assumed to be reduced.

## 4.1 – Infinitely near points and proximity

**(4.1) Definition** (infinitely near points). *Let $p$ be a (closed) point on a smooth surface $S$. The* points in the first infinitesimal neighbourhood *of $p$ are the (closed) points of the exceptional divisor of the blow-up of $S$ at $p$. Then by induction, for all $i > 1$, the* points in the $i$-th infinitesimal neighbourhood *of $p$ are the points in the first infinitesimal neighbourhood of some point in the $(i-1)$-th infinitesimal neighbourhood of $p$.*

Points infinitely near *to $p$ are those points lying in the $i$-th infinitesimal neighbourhood of $p$ for some $i > 0$. By contrast, points lying on $S$ itself are called* proper points. *An infinitely near point of $S$ is a point infinitely near to some proper point of $S$.*

**(4.2) Definition** (precedence). *Let $p, q$ be two points, proper or infinitely near, on $S$. We say that $p$* precedes *$q$ if $q$ is infinitely near to $p$. When this is the case, we write $p < q$. We will also use the usual variants, like $p \leqslant q$ and so on.*

**(4.3) Definition** (proximity). *Let $p, q$ be two points, proper or infinitely near, on $S$. We say that $q$ is* proximate *to $p$ if it belongs to the exceptional divisor $E_p$ of the blow-up of $S$ at $p$, or to the proper transform of $E_p$ in some further blow-up of $S$. When this is the case, we write $p \prec q$. We will also use the usual variants, like $p \preccurlyeq q$ and so on.*

For example, the points in the first and second infinitesimal neighbourhood of the origin of the affine plane lying on the proper transforms of the curve $V(y^2 - x^3)$ are both proximate to the origin.

Observe that a point proximate to $p$ lies either in the first infinitesimal neighbourhood of $p$ or in the first infinitesimal neighbourhood of a point proximate to $p$ (this follows directly from the definition; checking this is a good test of your understanding of the definition).

**(4.4) Lemma.** *Let $p, q$ be two points on $S$, proper or infinitely near, and assume that $p < q$. Then $q$ is proximate to either one or two points infinitely near or equal to $p$.*

*Proof.* Let $\varepsilon : S' \to S$ be a composition of blow-ups

$$S' = S_n \xrightarrow{\varepsilon_n} S_{n-1} \xrightarrow{\varepsilon_{n-1}} \cdots \xrightarrow{\varepsilon_2} S_1 \xrightarrow{\varepsilon_1} S_0 = S$$

such that $q$ is a proper point of $S'$; we may assume that each $\varepsilon_i$ is the blow-up at some (proper) point $p_i \in S_{i-1}$. The exceptional divisor of $\varepsilon$ is a disjoint union of trees of smooth rational curves, the irreducible components of which are the proper transforms $E_1, \ldots, E_n$ of the exceptional divisors of $\varepsilon_1, \ldots, \varepsilon_n$ respectively.

Thus the point $q$ lies on either one or two of the curves $E_1, \ldots, E_n$, and correspondingly it is proximate to either one or two points infinitely near or equal to $p$. □

**(4.5) Definition** (free and satellite points). *Let $p$ be a proper point of $S$, and $q$ a point infinitely near to $p$. The point $q$ is* free *if it is proximate to exactly one point infinitely near or equal to $p$; otherwise it is a* satellite. *One says that a satellite point is a* satellite of *the last free point that precedes it.*

Thus, in the setup of the proof of Lemma (4.4) above, the point $q$ is free if it is a smooth point of the exceptional divisor of $\varepsilon$, and it is satellite otherwise. Since the exceptional divisor is a union of trees, a satellite point is proximate to exactly two points infinitely near or equal to $p$.

All points in the first infinitesimal neighbourhood of a proper point are free. The first infinitesimal neighbourhood of a free (infinitely near) point contains exactly one satellite point, whereas that of a satellite point contains exactly two satellite points.

## 4.2 – Multiplicities and proximity equalities

**(4.6)** Let $C$ be a curve on a smooth surface $S$, and $q$ be an infinitely near point on $S$. One defines the multiplicity of $C$ at $q$ as follows. Let $\varepsilon : S' \to S$ be a composition of blow-ups such that $q$ is a proper point of $S'$, and let $C_q$ be the proper transform of $C$ in $S'$. The multiplicity of $C$ at $q$ is the multiplicity of $C_q$ at $q$. We denote it by $\mathrm{mult}_q(C)$.

We say that the point $q$ lies on $C$ if $\mathrm{mult}_q(C)$ is positive; it is a multiple point of $C$ if $\mathrm{mult}_q(C) > 1$, and a simple point if $\mathrm{mult}_q(C) = 1$.

The following result will arguably feel very natural. I will not prove it; if necessary, you can find it as [7, Theorem 3.3.1].

**(4.7) Proposition** (Noether's formula)**.** *Let $A$ and $B$ be two curves on a smooth surface $S$, and $p$ be a proper point on $S$. The intersection multiplicity of $A$ and $B$ at $p$ (possibly infinite) is*

$$(4.7.1) \qquad (A \cdot B)_p = \sum_{p \leqslant q} \mathrm{mult}_q(A)\, \mathrm{mult}_q(B)$$

*(the sum is on all points infinitely near or equal to $p$).*

The sequence of multiplicities of a curve at infinitely near points satisfies the following relations.

**(4.8) Proposition** (Proximity equalities)**.** *Let $C$ be a curve on a smooth surface $S$, and $p$ be a point on $S$, proper or infinitely near. Then,*

$$(4.8.1) \qquad \mathrm{mult}_p(C) = \sum_{p \prec q} \mathrm{mult}_q(C)$$

*(the sum is on all points $q$ proximate to $p$).*

*Proof.* Let $\varepsilon : S' \to S$ be a composition of blow-ups such that $p$ is a proper point of $S'$, and let $\varepsilon' : S'' \to S'$ be the blow-up at $p$, with exceptional divisor $E_p$. Let $C_p$ and $C'_p$ be the proper transforms of $C$ in $S'$ and $S''$ respectively. Then,

$$\mathrm{mult}_p(C) = \mathrm{mult}_p(C_p) = \sum_{q \in E_p} \left(E_p \cdot C'_p\right)_q,$$

where the first equality is by definition. Then one obtains the relation by applying Noether's formula (4.7.1) to each summand, noting that $E_p$ is smooth hence has multiplicity one at all its points, proper or infinitely near. □

**(4.9) Remark.** As a consequence of the proximity equalities, if some satellite point $q$ lies on a curve $C$, then it is proximate to a multiple point of $C$. Indeed, let $p$ be the immediate predecessor of $q$, and $p_0$ be the other point $q$ is proximate to. Then both $p$ and $q$ are proximate to $p_0$ hence, by the proximity equalities,

$$\mathrm{mult}_{p_0}(C) \geqslant \mathrm{mult}_p(C) + \mathrm{mult}_q(C),$$

and the sum is greater than one since both $p$ and $q$ lie on $C$.

## 4.3 – Equisingularity and Enriques diagrams

**(4.10) Definition** (singular points of a curve)**.** *Let $C$ be a curve on a smooth surface $S$. A point $p$ of $C$, proper or infinitely near, is a* singular point *of $C$ if one of the following conditions holds:*
*(a) $p$ is a multiple point of $C$;*
*(b) $p$ is a satellite point;*
*(c) $p$ precedes a satellite point on $C$.*

Note that there may be simple and free points on a curve $C$ preceding a satellite point still on $C$, and these are singular points of $C$. For example, this is the case of the point in the first infinitesimal neighbourhood of the origin on the affine curve $\mathrm{V}(y^2 - x^3)$. Nevertheless, a proper point of $C$ may precede a satellite point on $C$ only if it is multiple by Remark (4.9), hence the proper singular points of $C$ are the proper multiple points of $C$.

The following alternative characterization of singular points is helpful in grasping the concept (so let me leave the proof as an exercise). Let $p$ be an infinitely near point on $C$, and consider $\pi : S' \to S$ the minimal series of blow-ups giving rise to $p$ i.e., the minimal series of blow-ups such that $p$ is a proper point of $S'$. Call $E$ the exceptional divisor of $\pi$; it is a tree of smooth rational curves. Call $C'$ the proper transform of $C$ in $S'$. Then $p$ is a singular point of $C$ if and only if the intersection multiplicity of $C'$ with $E$ at $p$ is larger than 1. Thus, if $p$ is a simple point of $C$ (equivalently, of $C'$), either $C'$ is tangent to $E$ at $p$, or $p$ is a double point of $E$.

**(4.11)** In order to define equisingularity of two curves with planar singularities, we need the following definitions. Let $C$ be a reduced curve on a smooth surface $S$. The *cluster of singular points* of $C$ is the set of all singular points of $C$, proper or infinitely near, endowed with the two order relations of precedence and proximity; we denote it by $\mathrm{Sing}(C)$. The *augmented cluster of singular points* of $C$ is the set of all singular points of $C$, together with those non-singular points of $C$ that lie in the first infinitesimal neighbourhood of a singular point, endowed with the two order relations of precedence and proximity; we denote it by $\mathrm{Sing}^+(C)$. It should be noted that, by resolution of singularities for curves, both sets $\mathrm{Sing}(C)$ and $\mathrm{Sing}^+(C)$ are finite.

The set underlying the augmented cluster of singular points may be defined alternatively as follows. Let $p_1, \ldots, p_r$ be the singular points of $C$ which are proper on $S$. For all $i = 1, \ldots, r$, let $\gamma_{i,1}, \ldots, \gamma_{i,s_i}$ be the local branches of $C$ at $p_i$. For all $i = 1, \ldots, r$ and $j = 1, \ldots, s_i$, let $q_{i,j}$ be the first point (for the precedence order relation) on the branch $\gamma_{i,j}$ which is non-singular in the sense of (4.10). Then the set underlying $\mathrm{Sing}^+(C)$ is the set of all proper or infinitely near points on $C$ preceding or equal to a point $q_{i,j}$.

**(4.12) Definition** (equisingularity)**.** *Let $C \subseteq S$ and $D \subseteq T$ be two reduced curves on two smooth surfaces respectively. The two curves $C$ and $D$ are* equisingular *if there exists a bijection $\varphi : \mathrm{Sing}^+(C) \to \mathrm{Sing}^+(D)$ such that both $\varphi$ and $\varphi^{-1}$ preserve the two order relations of precedence and proximity.*

We do not require anything on the multiplicities of $C$ and $D$ at their singular points in the above definition. In fact, it comes automatically that $\varphi$ preserves multiplicities.

**(4.13) Proposition.** *Let $\varphi : \mathrm{Sing}^+(C) \to \mathrm{Sing}^+(D)$ be a bijection as in Definition (4.12), realizing the equisingularity of $C$ and $D$. Then, for all points $p \in \mathrm{Sing}^+(C)$:*

$$\mathrm{mult}_p(C) = \mathrm{mult}_{\varphi(p)}(D).$$

*Proof.* The proof is by decreasing induction with respect to the precedence order on $\mathrm{Sing}^+(C)$. If $p$ is a singular point of $C$ maximal in $\mathrm{Sing}^+(C)$ for precedence, then so is $\varphi(p)$, and thus $p$ and $\varphi(p)$ are simple points of $C$ and $D$ respectively, in particular they both have multiplicity 1. Now if $p$ is non-maximal in $\mathrm{Sing}^+(C)$, then

$$\mathrm{mult}_p(C) = \sum_{p \prec q} \mathrm{mult}_q(C)$$

by the proximity equality (4.8.1), hence

$$\mathrm{mult}_p(C) = \sum_{\varphi(p) \prec \varphi(q)} \mathrm{mult}_{\varphi(q)}(D)$$

by induction, and thus $\mathrm{mult}_p(C) = \mathrm{mult}_{\varphi(p)}(D)$ by the proximity equality. $\square$

Let me remark that the above proposition is false if one considers plain clusters of singular points instead of augmented clusters. Indeed the cluster associated to an ordinary singularity of any multiplicity $m$ consists of a single point. So with Sing instead of $\mathrm{Sing}^+$ in Definition (4.12), all ordinary singularities would be equisingular, regardless of the multiplicity, which is certainly not the definition we want.

If one prefers, it is possible to define equisingularity by considering the clusters Sing and adding the condition that the bijection $\varphi$ preserves multiplicities.

**(4.14)** It follows from Definition (4.12) that equisingularity is an equivalence relation; we call *equisingularity classes* the equivalence classes under this equivalence relation. My generic notation for equisingularity classes is S. By definition, an equisingularity class S may be represented by a finite set $\mathrm{Sing}^+$ equipped with the two order relations of precedence and proximity; as it will not cause any trouble that this set is defined only up to bijection as in Definition (4.12), I will abuse notation and denote it by $\mathrm{Sing}^+(\mathsf{S})$. By removing the maximal points of $\mathrm{Sing}^+(\mathsf{S})$ with respect to precedence, one finds the set $\mathrm{Sing}(\mathsf{S})$. Moreover, by Proposition (4.13), all points $p \in \mathrm{Sing}^+(\mathsf{S})$ come with an assigned multiplicity $\mathrm{mult}_p(\mathsf{S})$.

I will say loosely that a reduced curve $C$ in a smooth surface $S$, or merely a reduced curve $C$ with planar singularities, *has singularities of type* S if it belongs to the equisingularity class S.

One may prefer to reserve the term equisingularity class for equivalence classes of germs of curves under equisingularity, so that an equisingularity class is the datum of one singularity only. From this perspective, it would be better to name the notion of equisingularity class used in the present text multi-equisingularity class, as it is equivalent to the data of fintely many equisingular classes of germs.

**(4.15) Other definitions of equisingularity.** Equinsigularity had originally been defined by Zariski, in different terms, but with a similar flavour. The definition given here is equivalent to Zariski's, see [7, Section 3.8]. Also, see [8, Section 1.1] for a discussion of Teissier's version of Zariski's definition and more references.

Equisingularity is also equivalent to topological equivalence, where two germs of isolated planar curve singularities $(C_1, 0) \subseteq (\mathbf{C}^2, 0)$ and $(C_2, 0) \subseteq (\mathbf{C}^2, 0)$ are *topologically equivalent* if there exists a homeomorphism $(\mathbf{C}^2, 0) \to (\mathbf{C}^2, 0)$ mapping $(C_1, 0)$ to $(C_2, 0)$. For references about this statement, again see [7, Section 3.8] and [8, Section 1.1].

**(4.16) Enriques diagrams.** I now describe a graphic representation of equisingularity classes, designed to show the singular points and the precedence and proximity orders between them. It has been introduced by Enriques, see [10, Libro Quarto, Capitolo I (to be found in volume II)], and these representations are now called *Enriques diagrams*. Let S be an equisingularity class.

To materialize the precedence relation, we see $\mathrm{Sing}^+(\mathsf{S})$ as the set of vertices of an oriented graph, with the following definition: for all $p, q \in \mathrm{Sing}(\mathsf{S}^+)$, there is an edge from $p$ to $q$ if $q$ is in the first infinitesimal neighbourhood of $p$. Thus, $p$ precedes $q$ if and only if there is an oriented chain of edges joining $p$ to $q$. This makes $\mathrm{Sing}(\mathsf{S}^+)$ a disjoint union of (oriented) trees. For a tree, an orientation is equivalent to the choice of an origin; thus we will indicate the orientation simply by indicating the respective origins of the various trees composing $\mathrm{Sing}(\mathsf{S}^+)$, which are the vertices corresponding to proper singular points.

The proximity relation on the other hand is represented by drawing the edges in two possible ways, either straight or curved, according to the following rules:

(i) if $q$ is in the first infinitesimal neighbourhood of $p$, and $q$ is free (equivalently, $q$ is proximate only to $p$), then the edge from $p$ to $q$ is a smooth curved edge with (if $p$ is not a proper point) the same tangent at $p$ as the edge ending at $p$ (there is a unique such edge since the graph is a union of trees);

(ii) if $q$ is in the first infinitesimal neighbourhood of $p$, then all points infinitely near to $q$ and proximate to $p$ (all such points are satellite), as well as the edges joining them, are represented on a straight half-line starting at $q$ and orthogonal to the edge from $p$ to $q$.

To avoid self-intersections in the diagram, half-lines as in rule (ii) above are alternately oriented to the right and to the left of the preceding edge.

Although this is not necessary (by Proposition (4.13)), it is often helpful to indicate at each vertex of the diagram the multiplicity of the corresponding point. Kleiman and Piene have chosen to represent only $\mathrm{Sing}(\mathsf{S})$ in their version of Enriques diagrams, which they call minimal Enriques diagrams; this causes no trouble if the multiplicities are indicated, as we have seen in the discussion after the proof of Proposition (4.13). Their choice is more economic, but representing the whole $\mathrm{Sing}^+(\mathsf{S})$ makes the diagram more suggestive: it has the advantage of immediately showing all local branches, and it lets the diagram do a better job in representing the geometric difference between free and satellite points. An enlightening discussion of the nice geometric features of Enriques diagrams is included in [7, p. 99–100].

*Graphic convention.* The origins of the trees (i.e., proper singular points) are represented by square black dots. Infinitely near singular points are represented by round black dots. Non-singular points (i.e., points of $\mathrm{Sing}^+(\mathsf{S}) \setminus \mathrm{Sing}(\mathsf{S})$) are represented by simple points (thus, they are just the endpoints of the edge pointing toward them, without any further inking). The multiplicities of points in $\mathrm{Sing}(\mathsf{S})$ are indicated, but not those of points in $\mathrm{Sing}^+(\mathsf{S}) \setminus \mathrm{Sing}(\mathsf{S})$ (because the latter are all simple).

**(4.17) Example.** The six diagrams below are the respective Enriques diagrams of affine plane curves, defined by the indicated equation. Thus the first three are the diagrams of a node, an ordinary tacnode, and an ordinary oscnode, respectively; the last three are the diagrams of an ordinary cusp of order two, and two different cusps of order three.

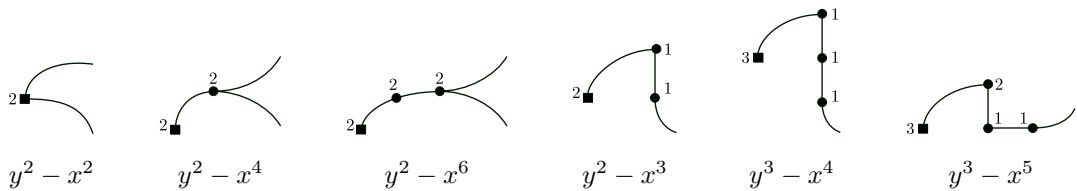

$y^2 - x^2$ $\quad\quad$ $y^2 - x^4$ $\quad\quad$ $y^2 - x^6$ $\quad\quad$ $y^2 - x^3$ $\quad\quad$ $y^3 - x^4$ $\quad\quad$ $y^3 - x^5$

For the two cusps of order three, name $p_0, p_1, p_2, p_3$ the points in S, following the order of precedence; then, for $V(y^3 - x^4)$, the two points $p_2, p_3$ are proximate to $p_0$ (in addition to their respective immediate predecessors), and for $V(y^3 - x^5)$, the point $p_2$ is proximate to $p_0$ whereas the point $p_3$ is proximate to $p_1$.

You may want to check the proximity equalities of (4.8) against these examples. More examples are given in Paragraph (4.28) below.

**(4.18) Notation.** In addition to the data mentioned in Paragraph (4.14), an equisingularity class $\mathsf{S}$ carries the following intrinsic characteristics:
— $r(\mathsf{S})$, the number of points in $\mathrm{Sing}^+(\mathsf{S}) \setminus \mathrm{Sing}(\mathsf{S})$;
— $\mathrm{roots}(\mathsf{S})$, the number of roots in $\mathrm{Sing}^+(\mathsf{S})$, i.e., singular points with no predecessor;
— $\mathrm{free}(\mathsf{S})$, the number of free singular points, i.e., free points in $\mathrm{Sing}(\mathsf{S})$ (by our Definition (4.5), free points are always infinitely near; beware that this is not the definition adopted by Kleiman and Piene).

Thus, $\mathrm{roots}(\mathsf{S})$ is the number of trees of which $\mathrm{Sing}^+(\mathsf{S})$ is the disjoint union; this is the number of proper singular points of the equisingularity class $\mathsf{S}$. On the other hand, $r(\mathsf{S})$ is the total number of local branches of curves in $\mathsf{S}$ (i.e., the sum of the number of local branches at the various proper singular points), see Paragraph (4.11).

**(4.19) Proposition.** *Let $C$ be a curve in the equisingularity class $\mathsf{S}$. Then the difference $p_a(C) - g(C)$ depends on $\mathsf{S}$ only, and equals*

$$(4.19.1) \qquad \delta(\mathsf{S}) = \sum_{p \in \mathrm{Sing}(\mathsf{S})} \binom{\mathrm{mult}_p}{2}.$$

*Proof.* The proof is by induction on the number of singular points. If there is none, the result is trivial. Implicitly, the curve $C$ lies on a smooth surface $S$. Let $p$ be a proper singular point of $C$, and consider the blow-up $\varepsilon : S' \to S$ at $p$, with exceptional divisor $E$. Then, the proper transform of $C$ in $S'$ is $C' = \varepsilon^*C - \mathrm{mult}_p.E$. By comparing the applications of the adjunction formula to $C \subseteq S$ and $C' \subseteq S'$ respectively, one finds

$$p_a(C') = p_a(C) - \binom{\mathrm{mult}_p}{2}.$$

The result then follows by induction, since the equisingularity class of $C'$ is obtained by removing $p$ from $\mathsf{S}$. □

Given a reduced curve $C$ in a smooth surface $S$, we define its *total Milnor number* as

$$\mu(C) = h^0(S, \mathcal{O}_S / \mathcal{J}_{C,S}),$$

where $\mathcal{J}_{C,S}$ is the Jacobian ideal of $C$ in $S$, locally generated by the partial derivatives of the equation of $C$ in $S$, see [II, Section 1]. Then, together with the celebrated Milnor–Jung formula, see for instance [15, Chapter I, Proposition 3.35], the above proposition gives the following. Beware that the Milnor–Jung formula is not always true in positive characteristic, see [3] for a recent account.

**(4.20) Corollary.** *The total Milnor number of $C$ depends on $\mathsf{S}$ only, and equals*

$$(4.20.1) \qquad \mu(\mathsf{S}) = 2\delta(\mathsf{S}) - r(\mathsf{S}) + \mathrm{roots}(\mathsf{S}).$$

**(4.21) Remark.** The definition of $r(\mathsf{S})$ in [20] is different from ours, but the two are equivalent. Indeed, the total number of local branches may be computed as follows:

$$r(\mathsf{S}) = \sum_{p \in \mathrm{Sing}(\mathsf{S})} \Big( \mathrm{mult}_p - \sum_{\substack{q \in \mathrm{Sing}(\mathsf{S}) \\ p \prec q}} \mathrm{mult}_q \Big).$$

This follows from the proximity equalities, and I leave the proof as an exercise.

### 4.4 – Equisingular loci

In this section we study, given a family $(F/Y, D)$ of curves on smooth surfaces as in Paragraph (1.4), the locus $Y(\mathsf{S}) \subseteq Y$ of curves belonging to a given equisingularity class $\mathsf{S}$. The results are taken from Kleiman and Piene [20, 22]. In particular, we show that it has an expected codimension, and that its actual codimension is everywhere smaller or equal to the expected codimension. The main statement is Proposition (4.23). We shall need the following numerical characters attached to an equisingularity class, in addition to those introduced in (4.18) above.

**(4.22)** Let $\mathsf{S}$ be an equisingularity class, as in (4.14). We set:

$$\dim(\mathsf{S}) = 2\,\mathrm{roots}(\mathsf{S}) + \mathrm{free}(\mathsf{S});$$
$$\deg(\mathsf{S}) = \sum_{p \in \mathrm{Sing}(\mathsf{S})} \binom{\mathrm{mult}_p + 1}{2};$$
$$\mathrm{expcod}(\mathsf{S}) = \deg(\mathsf{S}) - \dim(\mathsf{S}).$$

Now let $F/Y$ be a smooth projective family of surfaces over a Noetherian base, and $D$ be a relative effective divisor $D$ on $F$. Denote by $Y(\infty)$ the locus of points $y \in Y$ such that the curve $D_y$, i.e., the fibre of $D$ over $y$, is non-reduced. Moreover, for all equisingularity classes $\mathsf{S}$, denote by $Y(\mathsf{S})$ the locus of points $y \in Y$ such that the curve $D_y$ belongs to $\mathsf{S}$.

**(4.23) Proposition** (Kleiman–Piene). *In the above setup, the following hold:*
*(i) The set $Y(\infty)$ is closed in $Y$.*
*(ii) There are only finitely many equisingularity classes $\mathsf{S}$ such that $Y(\mathsf{S})$ is non-empty.*
*(iii) For all equisingularity classes $\mathsf{S}$, the set $Y(\mathsf{S})$ is constructible in $Y$, and at all point $y \in Y(\mathsf{S})$ we have*

$$\mathrm{codim}_Y\big(Y(\mathsf{S})\big)_y \leqslant \mathrm{expcod}(\mathsf{S}).$$

The idea of the proof is to elaborate on the following classical dimension count for nodal curves on a fixed surface $S$: for all $p \in S$, imposing that curves in some linear system $|L|$ have multiplicity at least 2 at $p$ amounts to 3 linear conditions on $|L|$; since $p$ moves in 2 dimensions along $S$, the locus of curves with one node has codimension at most one (if non-empty). The definition of equisingularity in the previous section makes it clear that any type of singularity may be understood as a collection of infinitely near multiple points, with appropriate precedence and proximity relations among them. So we will apply a similar dimension count to such collections.

**(4.24)** The first step is the construction of the locus $H(\mathsf{S})$ in the Hilbert scheme of $F/Y$ which parametrizes collections of infinitely near points with the configuration determined by the equisingularity class $\mathsf{S}$.

In order to understand this construction, let us see how to construct one such configuration on a given surface $S$. First choose pairwise distinct points on $S$ in correspondence with the roots of $\mathsf{S}$, and blow-up $S$ at these points. Next, for each root $p$, choose pairwise distinct points on the irreducible component $E_p$ of the exceptional divisor over $p$, in correspondence with the immediate successors of $p$ in $\mathrm{Sing}(\mathsf{S})$, and blow-up at all these newly chosen points. Then we want to consider those points in $\mathrm{Sing}(\mathsf{S})$ the second infinitesimal neighbourhood of the roots; the novelty is that from now on there may be satellite points. Let $q$ be one of the points chosen in the previous step; we choose pairwise distinct points on the irreducible component $E_q$ of the exceptional divisor over $q$, in correspondence with the immediate successors of $q$ in $\mathrm{Sing}(\mathsf{S})$, with the following new rule: for all immediate successors $r$ of $q$:

— if $r$ is free, the corresponding point on $E_q$ must be chosen off the intersection with the proper transforms of the previous exceptional divisors;
— if $r$ is proximate to some $p < q$, then the corresponding point must be the intersection point of $E_q$ with the proper transform of $E_p$.

We then repeat this last step until all points in $\mathrm{Sing}(\mathsf{S})$ have been chosen an incarnation as an infinitely near point on $S$.

In this construction, representatives of roots may be chosen virtually anywhere on $S$, representatives of free singular points may be chosen each virtually anywhere on an irreducible rational curve, and there is only one possible choice for satellite singular points. Thus the number of degrees of freedom in the construction is $\dim(\mathsf{S})$.

Now we must pass from a configuration of infinitely near points to a subscheme of $S$. This goes as follows. Let $\varepsilon : S' \to S$ be the sequence of all blow-ups considered above; include for convenience also the blow-ups of the last points of $\mathrm{Sing}(\mathsf{S})$ (with respect to precedence). For all point $p \in \mathrm{Sing}(\mathsf{S})$, which we identify with its incarnation, let $E_p$ be the irreducible exceptional divisor over $p$ (it dwells in some surface between $S$ and $S'$), and denote by abuse of notation by $\varepsilon^* E_p$ the pull-back of $E_p$ to $S'$ by those blow-ups coming after the blow-up of $q$; thus $\varepsilon^* E_p$ is a possibly reducible exceptional divisor of $\varepsilon$, with self-intersection $-1$. Then we set

$$\mathcal{I} = \varepsilon_* \mathcal{O}_{S'}\Big(-\sum_{p \in \mathrm{Sing}(\mathsf{S})} \mathrm{mult}_p.\, \varepsilon^* E_p\Big).$$

It defines a 0-dimensional subscheme of $S$, supported at the points we have chosen to incarnate the roots of $\mathsf{S}$, and of degree $\deg(\mathsf{S})$: the latter fact is a generalization of the fact for $p \in S$ with maximal ideal $\mathfrak{m}_p$, the ideal $\mathfrak{m}_p^{\mathrm{mult}_p}$ defines a 0-dimensional subscheme of length $\binom{\mathrm{mult}_p+1}{2}$. For a complete proof, using the framework of complete ideals, see [7, Section 8.3].

In sum, the above construction provides a locally closed subset $H(\mathsf{S})$ of dimension $\dim(\mathsf{S})$ of the Hilbert scheme of 0-dimensional, degree $\deg(\mathsf{S})$ subschemes of $S$, parametrizing collections of infinitely near points in the configuration determined by $\mathsf{S}$. In general, the following holds.

**(4.25) Proposition.** *Let $\mathsf{S}$ be an equisingularity class, and $F/Y$ be a smooth projective family of surfaces. The family of all collections of infinitely near points in the members of $F/Y$ in the configuration determined by $\mathsf{S}$ exists as a locally closed subset $H_{F/Y}(\mathsf{S})$ of $\mathrm{Hilb}_{F/Y}^{\deg(\mathsf{S})}$, the relative Hilbert scheme of $0$-dimensional subschemes of degree $\deg(\mathsf{S})$. Moreover, $H_{F/Y}(\mathsf{S})$ is smooth over $Y$, of pure (relative) dimension $\dim(\mathsf{S})$.*

The complete proof of this result is the object of [22], see Corollaries 5.5 and 5.8 in that paper. In fact Kleiman and Piene prove a more precise result, valid in arbitrary characteristic; they study infinitely near points from a functorial point of view, prove the existence of a so-to-speak moduli space $Q(\mathsf{S})$ of infinitely near points on $F/Y$ in the configuration determined by $\mathsf{S}$ and of a universal injection $\Psi : Q(\mathsf{S}) \to \mathrm{Hilb}_{F/Y}^{\deg(\mathsf{S})}$, which in characteristic 0 is an embedding.

Please see the introduction of [22] for references to previous works by various authors in the direction of Proposition (4.25) above.

**(4.26)** Consider now, in addition to the objects in Proposition (4.25), a relative effective divisor $D$ on $F/Y$. To relate $H_{F/Y}(\mathsf{S})$ to the study of members of the family $D$ in the equisingularity class $\mathsf{S}$, we form the intersection[5]

$$Z(\mathsf{S}) = H_{F/Y}(\mathsf{S}) \cap \mathrm{Hilb}_{D/Y}^{\deg(\mathsf{S})}.$$

The following is an important point for the proof of Proposition (4.23): by [1, Proposition 4], the Hilbert scheme $\mathrm{Hilb}_{D/Y}^{\deg(\mathsf{S})}$ is a subscheme of $\mathrm{Hilb}_{F/Y}^{\deg(\mathsf{S})}$ locally cut out by $\deg(\mathsf{S})$ equations. Thus, for all $z \in Z(\mathsf{S})$,

(4.26.1) $$\mathrm{codim}_{H_{F/Y}(\mathsf{S})}\bigl(Z(\mathsf{S})\bigr)_z \leqslant \deg(\mathsf{S}).$$

**(4.27)** We now prove Proposition (4.23).

*Proof of assertion* (i). In the notation of Section 2.1, the set $Y(\infty)$ is the image in $Y$ of the set of those $x \in X_2$ at which the fibre of $X_2/Y$ has positive dimension. The latter set is closed in $X_2$, which itself is closed in $F$ as explained in Section 2.1. The result follows since $F/Y$ is projective. □

*Proof of assertion* (ii). By flatness of $D/Y$, the curves $D_y$ have the same arithmetic genus for all $y \in Y$. Moreover, the number of components of their normalizations $\bar{D}_y$ is bounded from above (again by flatness, the degree of the $D_y$'s with respect to some relative polarization is constant along $Y$), so the geometric genus of the $D_y$'s is bounded from below. This implies that the $\delta$-invariants of the $D_y$'s, $\delta(D_y) = p_a(D_y) - g(D_y)$ is bounded from above. By Proposition (4.19), this bounds from above the number of multiple points, proper or infinitely near, of each $D_y$, as well as their multiplicities.

It remains to show that for each of these finitely many possible distributions of multiple points, there may be only finitely many simple singular points. Arguing local branch by local branch, the possible number of which is bounded from above by the previous considerations, this follows from [7, Theorem 3.5.8]. □

*Proof of assertion* (iii). Let $\mathsf{S}$ be an equisingularity class. By construction, $Y(\mathsf{S})$ is contained in $f(Z(\mathsf{S}))$, where $Z(\mathsf{S})$ is as in Paragraph (4.26), and $f$ denotes the structure morphism $\mathrm{Hilb}_{D/Y}^{\deg(\mathsf{S})} \to Y$.[6] The rest of $f(Z(\mathsf{S}))$ consists of points $y \in Y$ such that either $D_y$ is non-reduced, or it is reduced with singularities worse than $\mathsf{S}$. In the latter case, let $\mathsf{S}'$ be the corresponding equisingularity class. Then, again by construction, there is an injection

$$j : \mathrm{Sing}(\mathsf{S}) \hookrightarrow \mathrm{Sing}(\mathsf{S}')$$

such that for all $p \in \mathrm{Sing}(\mathsf{S})$, $j(p)$ has multiplicity at least that of $p$. Then, $\deg(\mathsf{S}')$ is larger than $\deg(\mathsf{S})$. The upshot is that

$$Y(\mathsf{S}) = f\bigl(Z(\mathsf{S})\bigr) \setminus \Bigl(Y(\infty) \cup \bigcup\nolimits_{\deg(\mathsf{S}')>\deg(\mathsf{S})} f\bigl(Z(\mathsf{S}')\bigr)\Bigr),$$

---

[5]if one is to work with the moduli space $Q(\mathsf{S})$, this has to be substituted with the fibre product $G(\mathsf{S}) = Q(\mathsf{S}) \times_{\mathrm{Hilb}_{F/Y}^{\deg(\mathsf{S})}} \mathrm{Hilb}_{D/Y}^{\deg(\mathsf{S})}$, see [23, Paragraph 2.5].

[6]for reference, I point out that in the series of articles of Kleiman and Piene, $U(\mathsf{S})$ is defined as the closure in $Y$ of $f(Z(\mathsf{S}))$, or a variant thereof: the definition in [20, 21] is exactly this, see [20, p. 227] and [21, p. 74], but that in [23, Definition 5.1] is different: there one removes the part supported on $Y(\infty)$.

| type | | normal form | | type | normal form | |
|---|---|---|---|---|---|---|
| $A_k$ | $(1 \leq k \leq 10)$ | $y^2 + x^{k+1}$ | | $J_{2,0}$ | $x^3 + ax^2y^2 + y^6$ | $(4a^3 + 27 \neq 0)$ |
| $D_k$ | $(4 \leq k \leq 10)$ | $x^2y + y^{k-1}$ | | $X_{1,1}$ | $x^4 + x^2y^2 + ay^5$ | $(a \neq 0)$ |
| $E_6$ | | $x^3 + y^4$ | | $J_{2,1}$ | $x^3 + x^2y^2 + ay^7$ | $(a \neq 0)$ |
| $E_7$ | | $x^3 + xy^3$ | | $X_{1,2}$ | $x^4 + x^2y^2 + ay^6$ | $(a \neq 0)$ |
| $E_8$ | | $x^3 + y^5$ | | $Z_{11}$ | $x^3y + y^5 + axy^4$ | |
| $X_{1,0}$ | | $x^4 + ax^2y^2 + y^4$ $(a^2 \neq 4)$ | | $Y_{1,1}$ | $x^5 + ax^2y^2 + y^5$ | $(a \neq 0)$ |

Table 2: Plane curve singularities with expected codimension at most 10

where the union ranges over the finitely many $\mathsf{S}'$ such that $Y(\mathsf{S}')$ is non-empty and $\deg(\mathsf{S}') > \deg(\mathsf{S})$. Since $Z(\mathsf{S})$ and the $Z(\mathsf{S}')$ are locally closed by Proposition (4.25), this proves that $Y(\mathsf{S})$ is constructible.

It remains to prove the bound on the codimension. Let $y \in Y(\mathsf{S})$. By construction there is a unique point $z \in Z(\mathsf{S})$ over $y$. Since $H_{F/Y}(\mathsf{S})$ is smooth of dimension $\dim(\mathsf{S})$ over $Y$ by Proposition (4.25), one has

$$\operatorname{codim}_Y\bigl(Y(\mathsf{S})\bigr)_y = \operatorname{codim}_{H_{F/Y}(\mathsf{S})}\bigl(Z(\mathsf{S})\bigr)_z - \dim(\mathsf{S}),$$

and by (4.26.1), this is smaller than or equal to $\deg(\mathsf{S}) - \dim(\mathsf{S}) = \operatorname{expcod}(\mathsf{S})$. □

**(4.28) Classification of singularities of small codimension.** In [20, Section 2], Kleiman and Piene show how one can use Enriques diagrams in order to classify plane curve singularities using combinatorial methods. In particular, they give the complete list of plane curve singularities with expected codimension at most 10 (as they point out, this may also be deduced from Arnold's classification results in [2]). As a sideremark, note that these are the singularities one expects to find on hyperplane sections of a sufficiently general surface $S \subseteq \mathbf{P}^N$ for $N \leq 10$, hence those to be considered in order to study the corresponding stratification of the dual $S^\vee \subseteq \check{\mathbf{P}}^N$, as we have done for surfaces in $\mathbf{P}^3$ in [XII].

We give their list in Table 2, in increasing order of expected codimension. It is limited to equisingularity classes with only one root, for any other class is the combination of several of those. Note that $X_{1,0}$ is the ordinary quadruple point, and has expected codimension 8. Singularities with expected codimension at most 8 are those in the first column (except $A_k$ and $D_k$ for $k > 8$).

Below, we give the Enriques diagrams for a sample of them (you may find all of them, in Kleiman–Piene's minimal form, in [20]). Note that we have already given the diagrams for singularities $A_1$, $A_3$, $A_5$, $A_2$, $E_6$, and $E_8$ (in that order) in Example (4.17).

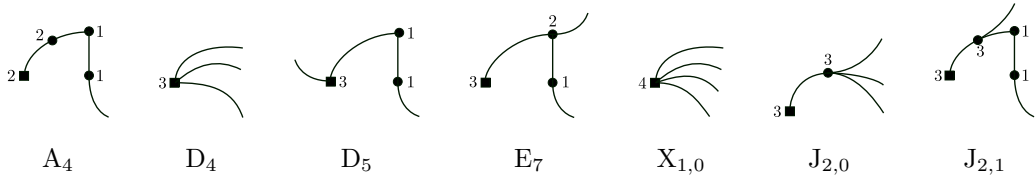

$A_4$  $D_4$  $D_5$  $E_7$  $X_{1,0}$  $J_{2,0}$  $J_{2,1}$

You may observe that $J_{2,0}$ is a triple point followed by $D_4$ (i.e., two infinitely near triple points), and $J_{2,1}$ is a triple point followed by $D_5$.

## 5 – Regularity conditions for curves in a fixed linear system

In this section we concentrate on the situation of Theorem (1.3), with the goal of providing sufficient conditions for the enumerative validity of the formulae in that theorem. Thus we will consider a smooth irreducible complex projective surface $S$, a line bundle $L$ on $S$, and curves in the fixed linear system $|L|$, which we will denote by $Y$. We denote by $D$ the universal divisor in $F = S \times Y$, see Section 3.1.

Before we start, let me state the following little lemma (see [20, Lemma (4.7)]), which is often overlooked. It asserts that imposing the passing through $c$ general points on $S$ defines a sufficiently general (for our enumerative purposes) linear subsystem of $|L|$ of codimension $c$. Set $n = \dim Y$, and let $r$ be an integer such that $0 \leqslant r \leqslant n$. For all $\mathbf{p} = (p_1, \ldots, p_{n-r}) \in S^{n-r}$, denote by $\Lambda(\mathbf{p})$ the $r$-dimensional linear subsystem of $Y$ of members of $Y$ passing through $p_1, \ldots, p_{n-r}$.

**(5.1) Lemma.** *Let $U$ be a reduced subscheme of $Y$ of pure codimension $r$, and $\partial U$ its boundary, i.e., $\partial U = \overline{U} \setminus U$. Assume that for all $y \in U$, the corresponding curve $D_y$ is reduced. Then there exists a non-empty open subset $\mathcal{S}^\circ$ of $S^{n-r}$ with the following property: for all $\mathbf{p} \in \mathcal{S}^\circ$, $\Lambda(\mathbf{p}) \cap U$ is finite and reduced, and $\Lambda(\mathbf{p}) \cap \partial U$ is empty.*

*Proof.* Let $D^\circ$ be the smooth locus in $D$ of the projection $D \to Y$, and form the fibered product $\mathcal{D}^\circ = D^\circ \times_Y \cdots \times_Y D^\circ$; it is smooth over $Y$. Hence $\mathcal{D}^\circ \times_Y U$ is reduced. Moreover, it is dense in $D^{n-r} \times_Y U$, where $D^{n-r}$ is the fibered product $D \times_Y \cdots \times_Y D$. Also, $D^{n-r} \times_Y U$ is of pure dimension $2(n-r)$. Consider the natural map,

$$D^{n-r} \times_Y U \to S^{n-r},$$

where $S^{n-r}$ is the plain product $S \times \cdots \times S$. By Sard's theorem, its fibers are finite and reduced over a non-empty open subset $\mathcal{S}_0$ of $S^{n-r}$. These fibers are simply the $\Lambda(\mathbf{p}) \cap U$.

On the other hand, $D^{n-r} \times_Y \partial U$ is of dimension at most $2(n-r) - 1$, so the map,

$$D^{n-r} \times_Y \partial U \to S^{n-r}$$

cannot be surjective. Hence, its fibers are empty over a nonempty subset $\mathcal{S}_0$ of $\mathcal{S}_1$. These fibers are the $\Lambda(\mathbf{p}) \cap \partial U$. Thus we may take $\mathcal{S}^\circ = \mathcal{S}_1$ and the lemma is proved. □

We shall need the following definition throughout.

**(5.2) Definition.** *Let $k$ be a non-negative integer. The line bundle $L$ is $k$-very ample if for all $0$-dimensional, length $k+1$, subscheme $Z$ of $S$, the restriction map*

$$H^0(S, L) \to H^0(Z, L|_Z)$$

*is surjective.*

### 5.1 – Application of equisingular deformation theory

**(5.3) Proposition.** *Let $\mathsf{S}$ be an equisingularity class such that $Y(\mathsf{S})$ is non-empty. If $L$ is $(\mu(\mathsf{S}) - 1)$-very ample, then $Y(\mathsf{S})$ is locally closed in $Y$, smooth, and of pure codimension*

$$\operatorname{codim}_Y\bigl(Y(\mathsf{S})\bigr) = \operatorname{expcod}(\mathsf{S}).$$

**(5.4)** In order to prove the proposition, we shall consider the deformation theory of the singularities of the members of $Y(\mathsf{S})$. We refer to [II] on this subject, and recall the following necessary facts.

Let $y \in Y(\mathsf{S})$, and consider the corresponding curve $D_y$. Then $D_y$ has isolated singularities, and we may consider a semi-universal deformation space $B$ for these singularities. The space $B$ is smooth, with tangent space at the origin

$$H^0(S, \mathcal{O}_S/\mathcal{T}_{D_y}),$$

where $\mathcal{T}_{D_y}$ is the Tjurina ideal of $D_y$ in $S$, locally generated by the equation of $D_y$ and its partial derivatives. Thus $\mathcal{T}_{D_y}$ contains the Jacobian ideal of $D_y$, hence its degree is at most the total Milnor number $\mu(D_y) = \mu(\mathsf{S})$, see Corollary (4.20) for this equality.

There is a natural map $\phi : U \to B$ defined on an analytic neighbourhood of $y$ in $Y$, and its differential at $y$ is the (quotient of the) restriction map

$$H^0(S, L)/H^0(S, \mathcal{O}_S) \longrightarrow H^0\big(S, L/\mathcal{T}_{D_y}L\big).$$

Since $L$ is $(\mu(\mathsf{S}) - 1)$-very ample, this differential is surjective, hence the map $\phi : U \to B$ is smooth (up to shrinking $U$ if necessary).

**(5.5)** In the deformation space $B$, there is a smooth closed subset $B^{\mathrm{es}}$, whose points represent the equisingular deformations of the singularities of $D_y$; the codimension of $B^{\mathrm{es}}$ in $B$ is

(5.5.1) $$\operatorname{codim}_B(B^{\mathrm{es}}) = \operatorname{expcod}(\mathsf{S}).$$

I give some indications about this below, but before that let us conclude with the proposition.

*Proof of Proposition (5.3).* Via the map $\phi$ defined in Paragraph (5.4) above,

$$Y(\mathsf{S}) \cap U = \phi^{-1}(B^{\mathrm{es}}).$$

Thus, $Y(\mathsf{S})$ is locally closed and smooth since $B^{\mathrm{es}}$ is closed and smooth in $B$, and

$$\operatorname{codim}_Y\big(Y(\mathsf{S})\big)_y = \operatorname{codim}_B(B^{\mathrm{es}})_0 = \operatorname{expcod}(\mathsf{S}).$$

$\square$

The two books [15] and [16] provide a fairly thorough presentation of equisingular deformation theory. See [16, Remark 2.2.44.1, (ii), p. 196] for a concise discussion of the smoothness of $B^{\mathrm{es}}$ and its codimension in $B$; Equality (5.5.1) then amounts to [16, Corollary 1.1.64]. The bottom line is that $B^{\mathrm{es}}$ is equal to the stratum $B^\mu$ in $B$ of deformations with constant Milnor number [15, Chapter II, Corollary 2.68], and the latter is smooth by [15, Chapter II, Theorem 2.61] and [15, Chapter II, Theorem 2.38].

The expression of the codimension of $B^{\mathrm{es}}$ in terms of $\operatorname{expcod}(\mathsf{S})$ is attributed to Shustin [35] by Kleiman and Piene; to the best of my knowledge, it follows from the theory only in a rather roundabout way. The argument is given in [20, p. 222–223], and goes as follows. First, for large enough $d$, there exists plane curves of degree $d$ in the equisingularity class $\mathsf{S}$, see [16, Theorem 4.5.13]. Next, if $d$ is large enough, the stratum $V_d^{\mathsf{S}} \subseteq |\mathcal{O}_{\mathbf{P}^2}(d)|$ parametrizing degree $d$ plane curves in the equisingularity class $\mathsf{S}$ is smooth of the expected dimension by [16, Theorem 2.2.40] or [16, Theorem 4.3.2]; in particular it has codimension $\operatorname{codim}_B(B^{\mathrm{es}})$ in $|\mathcal{O}_{\mathbf{P}^2}(d)|$. Finally, one computes as in the next Section 5.2 (see also the proof of [16, Corollary 1.1.64]) that the codimension of $V_d^{\mathsf{S}}$ in $|\mathcal{O}_{\mathbf{P}^2}(d)|$ equals $\operatorname{expcod}(\mathsf{S})$, for $d$ large enough. This altogether proves Equality (5.5.1).

## 5.2 – Application of Bertini's theorem

Let $\mathsf{S}$ be an equisingularity class, and recall the following notation from Section 4.4 (in the setup of the present section, $F = S \times Y$ as indicated above):

— $H_{F/Y}(\mathsf{S})$ is the parameter space in $\mathrm{Hilb}_{F/Y}^{\deg(\mathsf{S})}$ of infinitely near multiple points in the configuration determined by $\mathsf{S}$;

— $Z(\mathsf{S})$ is the intersection $H_{F/Y}(\mathsf{S}) \cap \mathrm{Hilb}_{D/Y}^{\deg(\mathsf{S})}$

— $f$ is the structure map $\mathrm{Hilb}_{D/Y}^{\deg(\mathsf{S})} \to Y$.

By construction, $Y(\mathsf{S})$ is contained in $f(Z(\mathsf{S}))$. Here we consider moreover

$$Z_0(\mathsf{S}) = \bigl(f^{-1}Y(\mathsf{S})\bigr) \cap Z(\mathsf{S}).$$

Then $f$ induces a bijection $Z_0(\mathsf{S}) \to Y(\mathsf{S})$, since by definition, for $y \in Y(\mathsf{S})$ the singularities of the curve $D_y$ gives exactly one configuration of infinitely near multiple points as determined by $\mathsf{S}$.

In the present situation, being $F = S \times Y$, $H_{F/Y}(\mathsf{S})$ is the product $H_S(\mathsf{S}) \times Y$, where $H_S(\mathsf{S}) \subseteq \mathrm{Hilb}_S^{\deg(\mathsf{S})}$ is the parameter space of configurations of infinitely near points in the fixed surface $S$. The projection $p : F \to S$ induces a map $Z(\mathsf{S}) \to H_S(\mathsf{S})$: despite appearances it is simple, indeed it maps a divisor containing a configuration of multiple infinitely near points to that configuration in the Hilbert scheme of $S$.

**(5.6) Proposition.** *Assume that $L$ is $\deg(\mathsf{S})$-very ample. Then the following hold:*
*(i) $Z(\mathsf{S})/H_S(\mathsf{S})$ is a bundle of projective spaces of codimension $\deg(\mathsf{S})$ in $|L|$;*
*(ii) $Z_0(\mathsf{S})$ is a dense open subset of $Z(\mathsf{S})$.*

*Proof.* Let us consider a point of $H_S(\mathsf{S})$: this is a configuration of infinitely near points on $S$ as determined by $\mathsf{S}$, materialized as a 0-dimensional subscheme of length $\deg(\mathsf{S})$ of $S$; denote by $\mathcal{I}$ the ideal of this subscheme. Then the fibre of $Z(\mathsf{S})$ over this point of $H_S(\mathsf{S})$ is merely the linear subsystem $\mathbf{P}H^0(S, \mathcal{I}L)$ of $Y$, parametrizing those divisors containing the chosen configuration infinitely near points.

Since $L$ is $(\deg(\mathsf{S}) - 1)$-very ample, the restriction map $H^0(S, L) \to H^0(S, L/\mathcal{I}L)$ is surjective, hence $H^0(S, \mathcal{I}L)$ has constant codimension $\deg(\mathsf{S})$ in $H^0(S, L)$ as $\mathcal{I}$ varies. Thus the fibres of $Z(\mathsf{S})/H_S(\mathsf{S})$ are all projective spaces of the same dimension, since $H^0(S, L/\mathcal{I}L) = \deg(\mathsf{S})$ for all $\mathcal{I}$, and assertion (i) follows.

Consider next the minimal blow-up $\varepsilon : S' \to S$ on which our infinitely near points are realized as proper points. Then on $S'$, $\varepsilon^{-1}\mathcal{I}\cdot \varepsilon^*L$ is an invertible sheaf which we will denote by $L'$, and the two linear systems $\mathbf{P}H^0(S, \mathcal{I}L)$ and $\mathbf{P}H^0(S', L')$ are the same. By definition, members of this linear system that are smooth and transverse to every component of the exceptional divisor of $\varepsilon$ map on $S$ to curves in the equisingularity class $\mathsf{S}$. Now, since $L$ is $\deg(\mathsf{S})$-very ample, the linear system $\mathbf{P}H^0(S', L')$ is free from base points and thus, by Bertini's theorem, smooth members transverse to the exceptional divisor form a dense open subset of $\mathbf{P}H^0(S', L')$. Assertion (ii) follows. □

**(5.7) Corollary.** *If $L$ is $\deg(\mathsf{S})$-very ample, then $Y(\mathsf{S})$ is a non-empty, dense, open subset of $f(Z(\mathsf{S}))$, and it has pure codimension $\mathrm{expcod}(\mathsf{S})$ in $Y$.*

*Proof.* By the proposition, the closures of the images $f(Z_0(\mathsf{S}))$ and $f(Z(\mathsf{S}))$ in $Y$ are equal. Since $Z_0(\mathsf{S})$ maps bijectively to $Y(\mathsf{S})$, $Y(\mathsf{S})$ is dense in $f(Z(\mathsf{S}))$. Now $Z(\mathsf{S})$ is a bundle of projective spaces of codimension $\deg(\mathsf{S})$ in $|L| = Y$ over $H_S(\mathsf{S})$, and the latter is smooth of dimension $\dim(\mathsf{S})$ by Proposition (4.25); the assertion on the codimension follows. □

### 5.3 – Conclusion

In this section we still consider a line bundle $L$ on a fixed surface $S$, and the corresponding pair $(F/Y, D)$ where $Y$ is the complete linear system $|L|$, $F = S \times Y$, and $D$ is the universal divisor. We shall use the results of the two preceding Sections 5.1 and 5.2 to prove the following.

**(5.8) Proposition.** *Let $1 \leqslant r \leqslant 8$, and assume that $L$ is $3r$-very ample.*
*(i) Let $\mathsf{S}$ be an equisingularity class such that $\operatorname{expcod}(\mathsf{S}) \leqslant r$. Then $Y(\mathsf{S})$ is locally closed in $Y$, smooth, of pure codimension $\operatorname{expcod}(\mathsf{S})$; moreover, $Y(\mathsf{S})$ is open and dense in $f(Z(\mathsf{S}))$.*
*(ii) Conditions (i) and (ii) of Theorem (1.6) hold.*

Let us first prove the first assertion, which is a direct application of the results in the two preceding sections.

*Proof of assertion* (i). By looking at the classification[7] in Paragraph (4.28) above, one finds the two following inequalities hold:

$$\mu(\mathsf{S}) \leqslant \operatorname{expcod}(\mathsf{S}) + 1 \quad \text{and} \quad \deg(\mathsf{S}) \leqslant 3\operatorname{expcod}(\mathsf{S}).$$

Then, we may apply Propositions (5.3) and (5.6), and the assertion follows. □

To prove assertion (ii) however, we need to prove that $\operatorname{codim}_Y\bigl(Y(\mathsf{S})\bigr) \geqslant \min\bigl(r+1, \operatorname{expcod}(\mathsf{S})\bigr)$ for all equisingularity class $\mathsf{S}$, and in general our positivity assumption on $L$ is not sufficient in order to apply the results of the previous sections to all possible $\mathsf{S}$ (even though those $\mathsf{S}$ for which $Y(\mathsf{S})$ is non-empty are finitely many). The trick to get control on all $Y(\mathsf{S})$ is to find for all $\mathsf{S}$ an equisingularity subclass to which one can apply the previous results.

Recall from Section 4.3 that an equisingularity class may be seen as a finite set $\operatorname{Sing}(\mathsf{S})$ endowed with the two order relations of precedence and proximity, and with a multilicity (or weight) attached to each element.

**(5.9) Definition.** *An equisingularity subclass $\mathsf{S}'$ of $\mathsf{S}$ is an equisingularity class corresponding to a subset of $\operatorname{Sing}(\mathsf{S})$, endowed with the induced relations of precedence and proximity, such that the multiplicities for $\mathsf{S}'$ are lower or equal to those for $\mathsf{S}$.*

Thus for instance, $A_{2i-1}$ is a subclass of $A_{2i}$ and of $A_{2i+1}$, but $A_{2i}$ is not a subclass of $A_{2i+1}$; also, $D_{2i}$ is a subclass of $D_{2i+1}$ and of $D_{2i+2}$, but $D_{2i+1}$ is not a subclass of $D_{2i+2}$. The result we shall use is the following.

**(5.10) Lemma.** *Let $1 \leqslant r \leqslant 8$, and $\mathsf{S}$ an equisingularity class. If $\operatorname{expcod}(\mathsf{S}) \geqslant r + 1$, then $\mathsf{S}$ contains a subclass $\mathsf{S}'$ such that $\operatorname{expcod}(\mathsf{S}') \geqslant r$ and $\deg(\mathsf{S}') \leqslant 3r$.*

This is [20, Lemma (4.4)]: the proof given there is combinatorial and I will skip it, although I genuinely don't take combinatorial to imply either easy or boring. The authors there announce a proof valid without any restriction on $r$ in a forthcoming paper, but I could not find it, see footnote 7 again.

*Proof of assertion* (ii). Let us first prove that condition (ii) of Theorem (1.6) holds. Let $\mathsf{S}$ be an equisingularity class with $\operatorname{expcod}(\mathsf{S}) > r$; we need to show that $Y(\mathsf{S})$ has codimension at least $r + 1$ at every point. To this end, let us consider $\mathsf{S}'$ as in Lemma (5.10) above. Being $\mathsf{S}'$ a subclass of $\mathsf{S}$, we have $Z(\mathsf{S}) \subseteq Z(\mathsf{S}')$, and thus $Y(\mathsf{S}) \subseteq f(Z(\mathsf{S}'))$. Since $\deg(\mathsf{S}') \leqslant 3r$, the class

---

[7] Kleiman and Piene [20, p. 228] state that the second inequality holds without the assumption $\operatorname{expcod}(\mathsf{S}) \leqslant 8$, and announce a proof in a forthcoming (at the time of [20]) paper; I have not been able to find it in their series of papers however, my (second) best guess being (the first one being my own short-sightedness) that it got lost in the flow of versions the whole body of their work has come through.

S$'$ is liable to Proposition (5.6) and Corollary (5.7); in particular, $Y(\mathsf{S}')$ is dense in the image of $Z(\mathsf{S}')$ in $Y$, and has pure codimension expcod($\mathsf{S}'$) in $Y$. If expcod($\mathsf{S}'$) $> r$, then $Y(\mathsf{S}')$ is already small enough for the inclusion $Y(\mathsf{S}) \subseteq f(Z(\mathsf{S}'))$ to grant $\mathrm{codim}_Y(Y(\mathsf{S}))$ at all points of $Y(\mathsf{S})$.

It may happen however that expcod($\mathsf{S}'$) $= r$. Then $\mathsf{S}'$ is a strict subclass of $\mathsf{S}$, as we assume expcod($\mathsf{S}$) $> r$, and we will argue that $Y(\mathsf{S})$ is contained in a proper closed subset of the closure of $Y(\mathsf{S}')$, which implies that $Y(\mathsf{S})$ indeed has codimension strictly greater than $r$ at all points. Since expcod($\mathsf{S}'$) $= r$, the class $\mathsf{S}'$ is liable to assertion (i), and thus $Y(\mathsf{S}')$ is locally closed in $Y$. This implies that $Y(\mathsf{S}')$ is a dense open subset of the closure of $f(Z(\mathsf{S}'))$. Since also $Y(\mathsf{S})$ is contained in $f(Z(\mathsf{S}'))$, it must be contained in a proper closed subset of $f(Z(\mathsf{S}'))$, and our claim is proved.

We will now prove that condition (i) of Theorem (1.6) holds with a similar argument. Consider the equisingularity class $\mathsf{S} = r\mathrm{A}_1$ of $r$ ordinary double points. Then the locus of non-reduced curves, $Y(\infty)$, lies in $f(Z(r\mathrm{A}_1))$ because we can find $r$ distinct double points schematically contained in a non-reduced curve. As above, $Y(r\mathrm{A}_1)$ is locally closed and dense in the closure $f(Z(r\mathrm{A}_1))$, of the expected codimension $r$; $Y(\infty)$ is then contained in a proper closed subset of the closure of $f(Z(r\mathrm{A}_1))$, hence it has codimension strictly greater than $r$ as we wanted to show. $\square$

We have now proved the two assertions of Proposition (5.8). This completes the proof of Theorem (1.3): the proposition says that the hypotheses of Theorem (1.6) are verified in the conditions of Theorem (1.3); we may thus apply Theorem (1.6) as explained in Section 3.1, which gives the wanted result, thanks also to the little Lemma (5.1).

### 5.4 – Particular cases

**(5.11) Plane curves.** For plane curves of degree $d$ with $\delta$ nodes, Kleiman and Piene prove that Formula (1.3.1) is valid if

$$d \geqslant \tfrac{1}{2}\delta + 1,$$

see [21, Theorem 3.1] and in particular [21, Lemma 3.3]. To do so, they use in a crucial way the following result of Greuel and Lossen [14].

Let $Y$ be the complete linear system of plane curves of degree $d$, and let $\mathsf{S}$ be an equisingularity class. Let $C$ be a plane curve with singularities of type $\mathsf{S}$. The equisingular deformations of $C$ in the plane are governed by the so-called *equisingular ideal* which we denote by $I^{\mathrm{es}}$, see [II, Section 6]. Call $\tau^{\mathrm{es}}$ the co-length of this ideal; this is the codimension of $B^{\mathrm{es}}$ in $B$ in the notation of Section 5.1 above. Then, by [14, Corollary 3.9], if $C$ is not the union of $d$ lines through a point, and if

$$4d > 4 + \tau^{\mathrm{es}},$$

then at the point $[C]$, $Y(\mathsf{S})$ is smooth and has codimension $\tau^{\mathrm{es}}$ in $Y$.

**(5.12) Curves on $K3$ surfaces.** If $S$ is a $K3$ or Enriques surface, the results of Knutsen in [25] give necessary and sufficient conditions for $k$-very ampleness of the line bundle $L$, hence sufficient conditions for the validity of Formula (1.3.1).

In particular, if $S$ is a $K3$ surface with Picard number 1, denoting by $L_1$ the positive generator of the Picard group, then for all $k \geqslant 0$:
— $L = L_1$ is $k$-very ample if $k \leqslant \tfrac{1}{4}L^2$;
— for all $m > 1$, $L = mL_1$ is $k$-very ample if $k + 1 < \tfrac{m-1}{m^2}L^2$.

**(5.13) Curves on Abelian surfaces.** In this case, Kleiman and Piene provide the following sufficient conditions, in [21, Lemma 5.3]. Let $S$ be an Abelian surface with Picard number 1, and denote by $L_1$ the positive generator of the group of divisors modulo numerical equivalence. Let $L$ be an effective line bundle on $S$. Then for all $k \geqslant 0$:
— if $L \equiv L_1$, then $L$ is $k$-very ample if $k+1 \leqslant \frac{1}{4}L^2$;
— if $L \equiv mL_1$ for some $m > 1$, then $L_1$ is $k$-very ample if $k+1 < \frac{m-1}{m^2}L^2$.

These conditions are derived by Kleiman and Piene directly from the Beltrametti–Sommese theorem [4]; this is also the basic idea underlying [25], but Knutsen's results are much more precise, valid without restriction on the Picard number.[8]

In [21, Section 5], Kleiman and Piene carry out an enumeration of curves on Abelian surfaces that goes beyond the situation of Theorem (1.3) and of the present section, as it is about curves in a given homology class instead of a given linear equivalence class. This can be done thanks to the flexibility of the setup of Theorem (1.6): then, one takes the parameter space $Y$ to be a projective bundle over the dual Abelian surface $\widehat{A}$, with fibre over the point $[L] \in \widehat{A}$ the complete linear system $|L|$.

# 6 – Inclusion-exclusion and the structure of node polynomials

In this section I explain, following Qviller [31], the structure of node polynomials in terms of the Bell polynomials. The main point is to give an intersection theoretic definition of, say, the numbers $a_i$ appearing in Formula (1.3.1), but we need to set things up before we can do so.

## 6.1 – The principle and an example

**(6.1) Introduction.** Let $\pi : F \to Y$ be a smooth projective family of surfaces, and let $D$ be a relative effective divisor on $F/Y$. The starting point is, as in the Kleiman–Piene approach, to consider the critical locus $X$ of the map $D \to Y$ induced by $\pi$: this is the $X_2$ from Section 2.1, but since here we will not consider any $X_i$ with $i > 2$, we denote it simply by $X$. The next steps, however, are different.

Whereas Kleiman–Piene consider an induced family over $X$, namely $(F_2/X_2, D_2)$ in the notation of Section 2.1, the critical locus of this new family, and then repeat this by induction, Qviller seeks to consider the self-intersection of $X$ in the sense of topologists, also known as the double locus in classical algebraic geometric terms. Concretely, he accesses it by means of the intersection
$$p_1^{-1}X \cap p_2^{-1}X \quad \text{in} \quad F \times_Y F$$
where $p_1$ and $p_2$ denote the two projections of $F \times_Y F$ to its factors; equivalently, this is the fibre product $X \times_Y X$. Of course, one considers the double locus when one is interested in 2-nodal curves, in general one ought to consider finer self-intersections, by means of
$$p_1^{-1}X \cap \cdots \cap p_\delta^{-1}X \quad \text{in} \quad F \times_Y \cdots \times_Y F,$$
or, equivalently, $X \times_Y \cdots \times_Y X$.

There are however contributions in this intersection that we need to get rid of in order to isolate that of $\delta$-nodal curves in the family $D$: there are the various polydiagonals in the product $X \times_Y \cdots \times_Y X$, which give excess contributions, as well as other contributions supported on various proper closed subsets of the polydiagonals; Example (6.2) below should give a good idea of what is going on.

---

[8] Beware that the application of Knutsen's theorem to Enriques surfaces given in [21, Remark 5.5] is pointless, as all Enriques surfaces have Picard number 10. Similarly, [21, Lemma 5.3] does not apply to bielliptic surfaces, since the latter all have Picard number 2.

Then the point is that the various polydiagonals are organized as a lattice, which makes it possible to take care of the unwanted contributions by a method of inclusion-exclusion, which finally leads to a formula with Bell polynomials.

**(6.2) Example.** Let me describe the situation when one is interested in 2-nodal curves, assuming everything happens as expected. This is taken from [31, Example 2.5], which is written in the context of curves in a linear system on a fixed surface.

The critical locus $X$ has dimension $-1$ over $Y$ (as it has codimension 3 in $F$, and $F$ has dimension 2 over $Y$), hence the fibre product $X \times_Y X$ is expected to have dimension $-2$ over $Y$. There are three loci of dimension at least $-2$ over $Y$ which appear as the support of a primary component of the scheme $X \times_Y X$, namely:
— the locus of ordered pairs of nodes, which projects onto the locus of binodal curves in $Y$;
— the diagonal in $X \times_Y X$, which projects onto the locus of 1-nodal curves;
— the locus of cusps.
The locus of ordered pairs of nodes has the expected dimension $-2$ over $Y$; the diagonal has dimension $-1$ over $Y$, hence it gives an excess contribution; the locus of cusps is contained in the diagonal, but has the expected dimension $-2$ over $Y$ and gives a proper contribution, with the multiplicity 2. We shall compute the contributions of the two latter loci to the intersection product $p_1^{-1}X \cdot p_2^{-1}X$; this will give us the contribution of ordered pairs of nodes, in a form that relates to Formula (1.6.1), as we shall see.

Let us first work out Formula (1.6.1) explicitly; the expressions are in terms, as before, of
$$v = c_1\bigl(\mathcal{O}_F(D)\bigr) \quad \text{and} \quad w_j = c_j\bigl(\Omega^1_{F/Y}\bigr) \text{ for } j = 1, 2.$$
The class of $X$ in $F$ is
$$b_1 = v^3 + w_1 v^2 + w_2 v.$$
From the algorithm in Paragraph (1.4), we find
$$Q_{b_1,2} = -7v - 6w_1,$$
and
$$\begin{aligned} b_2 &= (-7v - 6w_1) \cdot (v^3 + w_1 v^2 + w_2 v) \\ &= -7\,v^4 - 13\,v^3 w_1 - 6\,v^2 w_1^2 - 7\,v^2 w_2 \end{aligned}$$
(there is a term in $w_1 w_2$ which vanishes for degree/dimension reasons). The upshot is that binodal curves are enumerated by the following class in $Y$: $\frac{1}{2}\bigl((\pi_* b_1)^2 + \pi_* b_2\bigr)$.

The class $(\pi_* b_1)^2$ is the push-forward (or push-down) to $Y$ of the intersection product $p_1^{-1}X \cdot p_2^{-1}X$. We shall now see that the class $\pi_* b_2$ is the correction term relative to the unwanted contributions. Let us first work out the excess contribution of the diagonal, as this will also serve to introduce some notation. We apply [13, Proposition 9.1.1]: both $p_1^{-1}X$ and $p_2^{-1}X$ have codimension 3 in $F \times_Y F$, hence the proper (excess) contribution of the diagonal $\Delta_X \subseteq X \times_Y X$ to the intersection product is a class of dimension $n - 2$, where $n = \dim(Y)$, namely
$$\left\{ c\bigl(p_1^* N_{X/F}\big|_{\Delta_X}\bigr) \cdot c\bigl(p_2^* N_{X/F}\big|_{\Delta_X}\bigr) \cdot c\bigl(N_{\Delta_X/F^2}\bigr)^{-1} \cap [\Delta_X] \right\}_{n-2},$$
where $F^2$ denotes the fibre product $F \times_Y F$; also, $c$ denotes the total Chern class and $\{\_\}_k$ denotes the $k$-dimensional part of the class between braces; as usual, the inverse of a class is defined by taking the inverse in the sense of formal power series. Now, the embedding $\Delta_X \hookrightarrow F^2$ is the composition of the two regular embeddings $\Delta_X \hookrightarrow \Delta_F \hookrightarrow F^2$ (where $\Delta_F \subseteq F \times_Y F$ denotes the diagonal), with normal bundles $N_{X/F} \cong \mathcal{P}^1_{F/Y}(D)$ and $T_{F/Y}$ respectively; thus the class in the display above equals

$$\left\{ c\big(\mathcal{P}^1_{F/Y}(D)\big) \cdot c\big(T_{F/Y}\big)^{-1} \cap [X] \right\}_{n-2},$$

a class in the Chow group $\mathrm{CH}^1(X)$. By the exact sequence (2.3.1), we find this is

(6.2.1) $$\left\{ c\big(\mathcal{O}_F(D) \otimes \Omega^1_{F/Y}\big) \cdot c\big(\mathcal{O}_F(D)\big) \cdot c\big(T_{F/Y}\big)^{-1} \cap [X] \right\}_{n-2}.$$

**(6.3) Notation.** Following [18, Section 8], we let

$$\begin{aligned}\alpha &= 1 + \alpha_1 + \alpha_2 + \cdots \\ &= c\big(\mathcal{O}_F(D) \otimes \Omega^1_{F/Y}\big) \cdot c\big(T_{F/Y}\big)^{-1},\end{aligned}$$

see Paragraph (7.2). Thus, with the notation $v, w_1, w_2$ recalled above,

$$\begin{aligned}1 + \alpha_1 + \alpha_2 + \cdots &= \frac{(1+v)^2 + (1+v)w_1 + w_2}{1 - w_1 + w_2} \\ &= \big((1+v)^2 + (1+v)w_1 + w_2\big)\big(1 + (w_1 - w_2) + (w_1 - w_2)^2\big),\end{aligned}$$

hence

$$\begin{aligned}\alpha_1 &= 2(v + w_1); \\ \alpha_2 &= (v + w_1)(v + 2w_1); \\ \alpha_3 &= v(vw_1 + 3w_1^2 - 2w_2); \\ \alpha_4 &= v^2(w_1^2 - w_2);\end{aligned}$$

and $\alpha_k = 0$ for all $k > 4$ (note that several terms vanish for degree reasons in these computations).

**(6.4) Example (6.2) continued.** We now resume the computation of the proper (excess) contribution of $\Delta_X$ to $p_1^{-1}X \cdot p_2^{-1}X$. Since $X$ itself has pure dimension $n-1$, we need to consider the degree 1 homogeneous part of the class capped with $[X]$ in (6.2.1); with the notation above, it is

$$\alpha_1 + v = 3v + 2w_1,$$

hence the proper contribution of the diagonal is

$$(3v + 2w_1) \cdot (v^3 + w_1 v^2 + w_2 v).$$

This contribution does not include the contribution of the cuspidal locus, which is the support of a different component of the normal cone $\mathrm{NC}_{X \times_Y X}(F \times_Y F)$. To compute the latter, we rely on a formula of Kazarian's (in the case of a linear system on a fixed surface, one can use [XII, **] instead); these are presented in the following Section 7. Considering the class $R_{A_2}$ at [18, p. 711] and Theorem 8.4 in that same article, one finds that the (codimension 2) class on $Y$ enumerating 1-cuspidal curves is the push-down to $Y$ of

$$\alpha_1 \cdot (v^3 + w_1 v^2 + w_2 v) = 2(v + w_1) \cdot (v^3 + w_1 v^2 + w_2 v).$$

Since we agreed that the contribution to $p_1^{-1}X \cdot p_2^{-1}X$ of the 1-cuspidal locus has to be counted with the multiplicity 2, we find that the total of the contributions to be discarded is the pushdown to $Y$ of

$$\big((\alpha_1 + v) + 2\alpha_1\big) \cdot (v^3 + w_1 v^2 + w_2 v) = (7v + 6w_1) \cdot (v^3 + w_1 v^2 + w_2 v) :$$

this is $-b_2$ from the first part of the example. The upshot is that the contribution to $p_1^{-1}X \cdot p_2^{-1}X$ of the locus of 2-nodal curves is $(\pi_* b_1)^2 + \pi_* b_2$, which has to be divided by two as each binodal curve comes with two possible orderings of its nodes.

For $\delta > 2$ we don't quite have (yet?) as accurate a description of all contributions supported on some diagonal, as we shall see in Section 6.3. Still, these various contributions organize within the lattice of polydiagonals, and this is sufficient to understand the structure of node polynomials in terms of the Bell polynomials, as we will now see.

## 6.2 – Multiplicative structure of the lattice of polydiagonals

In this section we will state and prove (partially) Qviller's intersection theoretic definition of the classes $a_i = \pi_* b_i$ (in the notation of Theorem (1.6)), thus explaining the appearance of the Bell polynomials in Formula (1.6.1).

**(6.5) Partitions and polydiagonals.** We shall consider the various diagonals in the fibre product
$$X^r = X \times_Y \cdots \times_Y X.$$
They are in bijection with the (non-singleton) partitions $\varpi$ of $[\![1,r]\!] = \{1,\ldots,r\}$:[9] associated to $\varpi = \coprod_{1 \leqslant i \leqslant k} J_i$ is the diagonal $\Delta_\varpi \subseteq X^r$ defined as follows:

$$\Delta_\varpi = \big\{(x_1,\ldots,x_r) : \forall i = 1,\ldots,j,\ \forall s,t \in J_i,\ x_s = x_t\big\}$$

(if $\varpi$ is the singleton partition, this gives $X^r$ itself, which is not considered a diagonal). We denote by $\Pi_r$ the set of all partitions of $[\![1,r]\!]$, and by $\Pi_r^\circ$ the set of all non-singleton partitions of $[\![1,r]\!]$.

**(6.6) Definition.** *Let $r \geqslant 1$. For all non-singleton partition $\varpi \in \Pi_r^\circ$, we let $B_\varpi \in \mathrm{CH}_{n-r}(\Delta_\pi)$ be the part of the intersection product*

$$p_1^* X \cdot \ldots \cdot p_r^* X$$

*(where $p_1,\ldots,p_r$ are the various projections $X^r \to X$) suported in the diagonal $\Delta_\varpi$, in the sense of [13, Definition 6.1.2]: thus $B_\varpi$ is the sum of the contributions of the various distinguished subvarieties contained in $\Delta_\varpi$.*

For instance, in the situation of Example (6.2) above, for the single block partition $\varpi = [\![1,2]\!]$, $B_\varpi$ is the sum of both the excess contribution of the diagonal itself and the locus of 1-cuspidal curves.

**(6.7) Definition.** *Moreover, we let*

$$a_i(F/Y, D) = (-1)^{i-1}(i-1)!\,\pi_* B_{[\![1,r]\!]} \in \mathrm{CH}_{n-r}(Y),$$

*where $[\![1,r]\!]$ denotes the single block partition. Here we use by abuse of notation the same symbol for both $\pi : F \to Y$ and its restriction to the critical locus $X \subseteq F$.*

---

[9] in case you're wondering, '$\varpi$' is 'variant $\pi$'.

Note that $\Delta_{[\![1,r]\!]}$ is the small diagonal

$$\{(x_1,\ldots,x_r):\ x_1 = \cdots = x_r\} \subseteq X^r\,;$$

thus, it may be identified with $X$, and it makes sense in the above definition to consider the push-down of $B_{[\![1,r]\!]} \in \mathrm{CH}(\Delta_{[\![1,r]\!]})$ to $Y$ by $\pi$.

**(6.8) Theorem** (Qviller). *Let $F/Y$ be a family of surfaces and $D$ be a relative effective divisor. Let $\delta \geqslant 1$ and assume that the pair $(F/Y, D)$ is sufficiently generic for $\delta$-nodal curves (see below). Then, the class of the locus of $\delta$-nodal curves in $Y$ is*

$$(6.8.1) \qquad \frac{1}{\delta!}\, P_\delta\big(a_1(F/Y,D),\, \ldots,\, a_\delta(F/Y,D)\big)$$

*where $P_\delta$ is the $\delta$-th Bell polynomial as defined in Paragraph (1.1).*

The condition that the pair $(F/Y, D)$ be sufficiently generic for $\delta$-nodal curves means that the following conditions hold:
 (i) the locus $Y(\delta\mathsf{A}_1)$ of $\delta$-nodal curves has codimension $\delta$ in $Y$ at all its points;
 (ii) for all equisingularity classes $\mathsf{S}$ with at least $\delta$ proper singular points,[10] if $\mathsf{S} \neq \delta\mathsf{A}_1$, then the locus $Y(\mathsf{S})$ has codimension in $Y$ strictly larger than $\delta$ at all its points;
 (iii) the locus $Y(\infty)$ of non-reduced curves has codimension in $Y$ strictly larger than $\delta$ at all its points.

In the situation of a linear system $|L|$ on a fixed surface, i.e., $F = S \times |L|$ and $D$ is the universal divisor, cf. Paragraph (3.1) or Theorem (1.3), this sufficient genericity condition for $\delta$-nodal curves holds if $L$ is $r$-very ample by Kool–Shende–Thomas's proof of the Göttsche conjecture, see [XIV, Theorem (2.1)].

The following lemma captures the combinatorial structure which makes the Bell polynomials show up.

**(6.9) Lemma.** *Denote by $I_\delta$ the part of the intersection product $p_1^*X \cdot \ldots \cdot p_\delta^*X$ that is not supported on the polydiagonals. Then*

$$(6.9.1) \qquad I_\delta = \sum_{\varpi \in \Pi_\delta} m_\varpi B_\varpi$$

*where, for all partitions $\varpi$ of $[\![1,\delta]\!]$, the coefficient $m_\varpi$ is defined as follows: for all $i = 1,\ldots,\delta$, denote by $s_i(\varpi)$ the number of blocks of size $i$ in $\varpi$; then*

$$(6.9.2) \qquad m_\varpi = \prod_{i=1}^{\delta}\big((-1)^{i-1}(i-1)!\big)^{s_i(\pi)}.$$

*Proof.* By definition,

$$I_\delta = p_1^*X \cdot \ldots \cdot p_\delta^*X - \sum_{Z \subseteq \bigcup_{\varpi \in \Pi_\delta^\circ} \Delta_\varpi} (p_1^*X \cdot \ldots \cdot p_\delta^*X)_Z,$$

where the sum ranges over all distinguished varieties $Z$ (in the sense of [13, Definition 6.1.2]) of the intersection $p_1^*X \cdot \ldots \cdot p_\delta^*X$ contained in some diagonal of the fibre product $X \times_Y \cdots \times_Y X$, and for all $Z$, $(p_1^* X \cdot \ldots \cdot p_\delta^* X)_Z$ denotes the contribution of $Z$ to the intersection product,[11] i.e., the sum of those terms corresponding to an irreducible component of the normal cone with support equal to $Z$.

---

[10] in the terminology of Paragraph (4.18), this means that $\mathrm{Sing}(\mathsf{S})$ has at least $\delta$ distinct roots, equivalently the Enriques diagram of $\mathsf{S}$ has at least $\delta$ connected components.

By definition the distinguished subvarieties are irreducible, so we may rearrange the second summand in the above display in the form

$$(6.9.3) \qquad - \sum_{Z \subseteq \bigcup_{\varpi \in \Pi_\delta^\circ} \Delta_\varpi} (p_1^* X \cdot \ldots \cdot p_\delta^* X)_Z = \sum_{\varpi \in \Pi_\delta^\circ} \sum_{Z \subseteq \Delta_\varpi} m_\varpi . (p_1^* X \cdot \ldots \cdot p_\delta^* X)_Z,$$

where the $m_\varpi$'s are arranged so that for all $Z$, the relevant $m_\varpi$'s sum up to $-1$. Thus, the $m_\varpi$'s are defined by induction on $\Pi_\delta^\circ$ with respect to the partial order by refinement, as we will now see explicitly. Also, we will prove that they verify Equality (6.9.2).

For two partitions $\varpi$ and $\varpi'$ of $[\![1,\delta]\!]$, we say that $\varpi$ *refines* $\varpi'$, and we write $\varpi \prec \varpi'$,[12] if every block of $\varpi$ is contained in a block of $\varpi'$. We let $\hat{0}$ and $\hat{1}$ denote the singleton and the single-block partition of $[\![1,\delta]\!]$ respectively; thus $\hat{0}$ is a refinement of every partition, and every partition is a refinement of $\hat{1}$. Moreover, for all $\varpi$, we denote by $|\varpi|$ the number of its blocks; thus $|\hat{0}| = \delta$ and $|\hat{1}| = 1$.

Now to the definition of the $m_\varpi$'s. Let us first consider the largest diagonals, i.e., those corresponding to a partition $\varpi$ with only one non-trivial block of size $2$, equivalently those $\varpi$ with $|\varpi| = \delta - 1$. We set

$$(6.9.4) \qquad m_\varpi = -1 \quad \text{if} \quad |\varpi| = \delta - 1,$$

so that those distinguished varieties $Z$ contained in one only of the largest diagonals (and therefore in no smaller diagonal) occur with the correct coefficient $-1$ in the summation (6.9.3). Now assume $Z$ is contained in some smaller diagonals, and let $\Delta_\varpi$ be the smallest diagonal in which it is contained; then a diagonal $\Delta_{\varpi'}$ contains $Z$ if and only if $\varpi'$ is a refinement of $\varpi$. Thus we impose the recursive relation

$$(6.9.5) \qquad m_\varpi = -1 - \sum_{\substack{\varpi' \in \Pi_\delta^\circ \\ \varpi' \prec \varpi}} m_{\varpi'} \quad \text{for all} \quad \varpi \in \Pi_\delta^\circ,$$

so that each $Z$ contained in some diagonal occurs with total coefficient $-1$ in the summation (6.9.3). The two conditions (6.9.4) and (6.9.5) are sufficient to define $m_\varpi$ for all $\varpi \in \Pi_\delta^\circ$ by induction, and with such a definition the relation (6.9.3) holds.

Now let us set in addition $m_{\hat{0}} = 1$ for the singleton partition. As we have observed in Paragraph (6.5), one has $\Delta_{\hat{0}} = X^\delta$ (which strictly speaking is not a diagonal), hence

$$p_1^* X \cdot \ldots \cdot p_\delta^* X = \sum_{Z \subseteq \Delta_{\hat{0}}} m_{\hat{0}} . (p_1^* X \cdot \ldots \cdot p_\delta^* X)_Z,$$

and therefore

$$I_\delta = \sum_{\varpi \in \Pi_\delta} \sum_{Z \subseteq \Delta_\varpi} m_\varpi . (p_1^* X \cdot \ldots \cdot p_\delta^* X)_Z$$

(note that the summation ranges over $\Pi_\delta = \Pi_\delta^\circ \cup \{\hat{0}\}$). It is thus sufficient to prove Formula (6.9.2) for all $\varpi \in \Pi_\delta$ in order to complete the proof of the lemma.

---

[11] we avoid the notation $(p_1^* X \cdot \ldots \cdot p_\delta^* X)^Z$, which has a different meaning in [13, Definition 6.1.2].

[12] agreed, this is the same notation as for the proximity relation, but this should not create any trouble.

Observe that $m_\varpi$ is defined for all $\varpi \in \Pi_\delta$ by the two conditions

$$(6.9.6) \qquad m_{\hat{0}} = 1, \quad \text{and} \quad m_\varpi = -\sum_{\substack{\varpi' \in \Pi_\delta \\ \varpi' \prec \varpi}} m_{\varpi'} \quad \text{for all} \quad \varpi \in \Pi_\delta,$$

as they imply the two conditions (6.9.4) and (6.9.5). It follows that

$$m_\varpi = \mu(\hat{0}, \varpi)$$

where $\mu$ is the Möbius function of the partially ordered set $\Pi_\delta$, see [36, Section 3.7 and Display (3.15)].

To compute the relevant values of this Möbius function, consider the partially ordered sets $\Pi_r$ for all $r$, with Möbius function $\mu_r$, and with singleton and single-block partitions $\hat{0}_r$ and $\hat{1}_r$ respectively. The induction relations (6.9.6) imply that

$$\mu_r(\hat{0}_r, \hat{1}_r) = -(r-1)\,\mu_{r-1}(\hat{0}_{r-1}, \hat{1}_{r-1}),$$

hence

$$\mu_r(\hat{0}_r, \hat{1}_r) = (-1)^{r-1}(r-1)!\,.$$

Now, if $\varpi$ is a partition $J_1 \coprod \cdots \coprod J_k$ with each block $J_j$ of size $r_j$, one has

$$\mu(\hat{0}, \varpi) = \prod_{j=1}^{k} \mu_{r_j}(\hat{0}_{r_j}, \hat{1}_{r_j}) = \prod_{j=1}^{k} \bigl( (-1)^{r_j - 1}(r_j - 1)!\bigr),$$

hence finally

$$m_\varpi = \mu(\hat{0}, \varpi) = \prod_{i=1}^{r} \bigl(\mu_i(\hat{0}_i, \hat{1}_i)\bigr)^{s_i(\varpi)},$$

and Formula (6.9.2) holds, which proves the lemma. For more details on the computation of the Möbius function, please see [36, Example 3.10.4]. □

In order to prove the theorem, we will also need the following lemma. It explains why it is sufficient to understand the classes $B_{[1,r]}$ (or $B_{\hat{1}_r}$ depending on one's favourite notation), corresponding to the smallest diagonals $\Delta_{[1,r]} \subseteq X^r$ for all $r = 2, \ldots, \delta$, to understand the classes $B_\varpi$ of an arbitrary diagonal $\Delta_\varpi \subseteq X^\delta$. For example, the part of $p_1^* X \cdot \ldots \cdot p_6^* X$ suported in the diagonal

$$\Delta_{1|23|456} = \bigl\{(x_1, \ldots, x_6): \ x_2 = x_3 \quad \text{and} \quad x_4 = x_5 = x_6\bigr\}$$

may be computed from the parts of $p_1^* X \cdot p_2^* X$ and $p_1^* X \cdot p_2^* X \cdot p_3^* X$ in $X^2$ and $X^3$ supported in the small diagonals $X_{[1,2]} \subseteq X^2$ and $X_{[1,3]} \subseteq X^3$. This explains the relevance of Definition (6.7) and its role in Formula (6.8.1).

**(6.10) Lemma.** *For all partitions $\varpi \in \Pi_\delta$, the following equality of classes in $Y$ holds:*

$$\gamma_{\delta *} B_\varpi = \prod_{i=1}^{\delta} \bigl(\gamma_{i *} B_{[1,i]}\bigr)^{s_i(\varpi)},$$

*where $\gamma_i$ denotes the projection $F^i = F \times_Y \cdots \times_Y F \to Y$ for all $i$. Recall also that $s_i(\varpi)$ is the number of blocks of size $i$ in the partition $\varpi$.*

Let me omit the proof of this lemma, as I believe the statement is intuitive and the proof is essentially technical and doesn't add much to one's global understanding of the story. Please refer to [31, Proposition 3.10] for details. Now, let us prove the main result of this section.

*Proof of Theorem (6.8).* The genericity assumption ensures that the push-down to $Y$ of the class $I_\delta$ of Lemma (6.9) enumerates $\delta$-nodal curves, with an ordering of the nodes, in the family $D$. Indeed, curves with strictly less than $\delta$ proper singular points correspond in $X^\delta$ to loci that are contained in the polydiagonals, hence their contribution to the intersection product $p_1^* X \cdot \ldots \cdot p_\delta^* X$ is discarded in the definition of $I_\delta$. On the other hand, among curves with at least $\delta$ singular points, those with more singularities than $\delta$ nodes occur in codimension at least $\delta+1$ in $Y$, hence they don't contribute to $p_1^* X \cdot \ldots \cdot p_\delta^* X$. The upshot is that $\delta$-nodal curves in $Y$ are enumerated by the class $\frac{1}{\delta!} \gamma_{\delta *} I_\delta$, with $\gamma_\delta$ the projection $F^\delta \to Y$ as in Lemma (6.10).

Now, one has $I_\delta = \sum_{\varpi \in \Pi_\delta} m_\varpi B_\varpi$ as in Lemma (6.9). For all partitions $\varpi \in \Pi_\delta$, by Lemma (6.10) and the definitions of the coefficients $m_\varpi$ in Lemma (6.10) and the classes $a_i(F/Y, D)$ in Definition (6.7), the following equality of classes in $Y$ holds,

$$\gamma_{\delta *}(m_\varpi B_\varpi) = \prod_{i=1}^\delta a_i(F/Y, D)^{s_i(\varpi)}.$$

It follows that

$$\gamma_{\delta *} I_\delta = \sum_{s_1 + \cdots + \delta s_\delta = \delta} e_{s_1,\ldots,s_\delta} \left( \prod_{i=1}^\delta a_i(F/Y, D)^{s_i} \right)$$

where $e_{s_1,\ldots,s_\delta}$ is the number of partitions $\varpi \in \Pi_\delta$ with $s_1$ blocks of size 1, ..., $s_\delta$ blocks of size $\delta$. By Corollary (A.4) on Bell polynomials, the above equality reads

$$\gamma_{\delta *} I_\delta = P_\delta \big( a_1(F/Y, D), \ldots, a_\delta(F/Y, D) \big),$$

which proves the theorem. □

### 6.3 – Further considerations

Theorem (6.8) presents the issue, as acknowledged by Qviller himself, that although the classes $a_i(F/Y, D)$ are given an intersection theoretic definition, they are not quite directly computable. In particular, in view of the Göttsche conjecture, it would be desirable to have a direct proof that, in the context of a linear system on a fixed surface (i.e., $F = S \times |L|$ and $D$ is the universal divisor over $|L|$), the classes $a_i(F/Y, D)$ are linear combinations of the four Chern numbers of the polarized surface $(S, L)$, i.e., $d = L^2$, $k = L \cdot K_S$, $s = K_S^2$, $x = c_2(S)$. Qviller gives a number of hints and conjectures pointing in this direction in [31, Sections 4–5], that I will leave aside in the present text.

Also, the proof of Theorem (6.8) does not provide as precise a picture as that gotten in Example (6.2) when we computed the class of the locus of binodal curves. Here I will very briefly report on Qviller's considerations in this direction.

**(6.11) Excess contribution of the polydiagonals.** It is interesting to isolate, in the intersection product $p_1^* X \cdot \ldots \cdot p_\delta^* X$ in $F \times_Y \cdots \times_Y F$, the contributions of those distinguished varieties of the correct dimension $n - \delta$, where $n = \dim Y$. In the context of Example (6.2), this means the contribution of both the binodal locus and the cuspidal locus; note that the latter is contained in the diagonal.

This amounts to the computation of the (excess) contributions of all the diagonals $\Delta_\varpi$ in $X^\delta = X \times_Y \cdots \times_Y X$. To do this in a sensible way, one would try to express the result in terms of the contributions of the small diagonals $\Delta_{[1,i]}$ in $X^i$ for all $i = 2, \ldots, \delta$, as has been done for Theorem (6.8). One would thus build up a definition similar to Definition (6.7) for the $a_i(F/Y, D)$, in which the class $B_{[1,i]}$ (the sum of the contributions of the various distinguished subvarieties contained in $\Delta_{[1,i]}$) would be substituted by the mere contribution of the sole distinguished subvariety $\Delta_{[1,i]}$.

It is indeed possible to express the contribution of an arbitrary polydiagonal $\Delta_\varpi$ in these terms, as the appropriate result analogous to Lemma (6.10) holds, see [31, Theorem 4.3]. This doesn't work directly, however, to express the sum of the contributions of all the polydiagonals: the problem is that these excess contributions do not satisfy the principle of inclusion-exclusion, which is the guiding line of Lemma (6.9) above. The basic issue is that Segre classes do not satisfy the principle of inclusion-exclusion.

Let me try to make this intelligible on an example, taken from [31, Example 4.4]. Consider $\delta = 3$. Then there are four diagonals to be considered, the small diagonal $\Delta_{123}$ and the three large diagonals $\Delta_{12|3}$, $\Delta_{1|23}$, and $\Delta_{2|13}$. All three large diagonals contain the small diagonal. Thus, according to the principle of inclusion-exclusion, if one sums the contributions of the three large diagonals, then the contribution of the small diagonal is counted three times instead of one only, hence the correct total contribution should be the sum of the contributions of the three large diagonals minus twice the contribution of the small diagonal. It turns out however, that this does not give the correct answer.

This being observed, Qviller suggests to compute correction terms to the formula predicted by the inclusion-exclusion principle. To do this end, he uses the following multiple point formulae.

**(6.12) Multiple point formulae.** The set of $r$-fold points of a map of schemes $f : V \to W$ is the set $M_r \subseteq V$ of points $x$ such that there exists $r$ distinct points (possibly infinitely near each other) in $f^{-1}f(x)$. It is endowed with a natural scheme structure. The object we are interested in in the previous paragraph, namely the sum of the contributions to $p_1^* X \cdot \ldots \cdot p_\delta^* X$ in $F \times_Y \cdots \times_Y F$ of those distinguished varieties of the correct dimension $n - \delta$, may be understood as the scheme of $\delta$-fold points of the map $\pi : X \to Y$. For instance, in the situation of Example (6.2), the scheme of 2-fold points of $X \to Y$ consists of the locus of binodal curves plus the locus of cuspidal curves, the latter with multiplicity 2.

In certain cases there exist formulae which give the class of $M_r$ in the Chow ring of $V$ in terms of the Chern classes of $f$, where the total Chern class of $f$ is defined by

$$c(f) = \frac{f^* c(T_W)}{c(T_V)}.$$

For instance, Kleiman gives in [19] the following formulae, valid for the map $\pi : X \to Y$ we are considering (the $c_i$'s are the Chern classes of $\pi$):

$$\begin{aligned} m_2 &= \pi^* \pi_*[X] - c_1 \cap [X]\,; \\ m_3 &= \pi^* \pi_* m_2 - 2c_1 \cap m_2 + 2c_2 \cap m_1\,; \\ m_4 &= \pi^* \pi_* m_3 - 3c_1 \cap m_3 + 6c_2 \cap m_2 - 6(c_1 c_2 + 2c_3) \cap m_1. \end{aligned}$$

Qviller uses the latter formulae to derive the following expressions of the push-down to $Y$ of the contributions to $p_1^* X \cdot \ldots \cdot p_\delta^* X$ of distinguished varieties of the correct dimension for $\delta = 2, 3, 4$:

$$\begin{aligned} \pi_* m_2 &= P_2(d_1, d_2)\,; \\ \pi_* m_3 &= P_3(d_1, d_2, d_3 + d_3')\,; \\ \pi_* m_4 &= P_4(d_1, d_2, d_3 + d_3', d_4 + d_4')\,; \end{aligned}$$

here, $P_2, P_3, P_4$ are the Bell polynomials, for all $i$ the class $d_i$ is defined by the contribution of the small diagonal $\Delta_{[1,i]}$ to $p_1^* X \cdot \ldots \cdot p_i^* X$, with a coefficient $(-1)^{i-1}(i-1)!$ dictated by the principle of inclusion-exclusion, and the classes $d_i'$ are ad hoc correction terms compensating the failure of the principle of inclusion-exclusion for the polydiagonals.

The strategy in [19] to derive the above multiple point formulae is similar in spirit to the inductive strategy of Kleiman and Piene described in Section 2: to compute $m_3$, one identifies $M_3$ with the locus of double points of the map $M_2 \to Y$, and so forth.

## 7 – Kazarian's Thom polynomials

In this section I give some formulae by Kazarian [18] which generalize those of Theorems (1.3) and (1.6) to arbitrary multisingularities. The existence of such formulae had been previously proven by Li and Tzeng [28], see also Rennemo [32], and generalizes the predictions of the Göttsche conjecture.

**(7.1) Thom polynomials.** The original result proved by Thom (see [18, Section 1] and the references therein) is the following. Let $f : M \to N$ be a holomorphic map between non-singular complex varieties. Let $\mathsf{S}$ be a class of local singularities of maps: this is a type of singularity supported at one point, a typical example is the class of germs of maps with dimension of the kernel of the differential at least $r$; see [18, Section 1] for a precise definition. There exists a universal polynomial, nowadays called Thom polynomial, in the Chern classes of the map $f$, defined by

$$c(f) = 1 + c_1(f) + \cdots = \frac{f^* c(T_N)}{c(T_M)},$$

which gives the class $[M(\mathsf{S})]$ of the locus of points in $M$ at which $f$ has a singularity of type $\mathsf{S}$, if $f$ is a generic map.

This is generalized by Kazarian to multisingularities in [18]. In this context, a multisingularity is an arbitrary finite ordered set of local singularities $\mathsf{S} = (\mathsf{S}_1, \ldots, \mathsf{S}_r)$, and the corresponding locus

$$M(\mathsf{S}) \subseteq M \times_N \cdots \times_N M$$

is the set of $r$-tuples $(x_1, \ldots, x_r)$ such that the $x_i$'s are pairwise distinct and $f$ has a singularity of type $\mathsf{S}_i$ at $x_i$ for all $i = 1, \ldots, r$. Denote by $\tilde{m}_\mathsf{S}$ and $\tilde{n}_\mathsf{S}$ the push-forward of the class $[M(\mathsf{S})]$ to $M$ and $N$, by the first projection $p_1$ and by $fp_1$ respectively. They come with multiplicities determined by the repetitions of the various local singularity classes in $\mathsf{S}$: if $\mathsf{S}$ is some ordering of $k_1$ singularities of type $\mathsf{T}_1, \ldots, k_s$ singularities of type $\mathsf{T}_s$, with the $\mathsf{T}_j$'s pairwise distinct, set

$$\mathrm{aut}(\mathsf{S}) = k_1! \cdots k_s!;$$

then there are integral classes $m_\mathsf{S}$ and $n_\mathsf{S}$ such that

$$\tilde{m}_\mathsf{S} = \mathrm{aut}(\mathsf{S}^\flat).m_\mathsf{S} \quad \text{and} \quad \tilde{n}_\mathsf{S} = \mathrm{aut}(\mathsf{S}).n_\mathsf{S},$$

where $\mathsf{S}^\flat = (\mathsf{S}_2, \ldots, \mathsf{S}_r)$; with these definitions, one has $n_\mathsf{S} = f_* m_\mathsf{S}$.

Kazarian proves the existence of a universal polynomial with integral coefficients computing, for a generic proper map, the class of $\tilde{m}_\mathsf{S}$ in terms of the Chern class of the map $f$ and the so-called Landweber–Novikov classes: by definition the latter are the following classes in $N$,

$$s_I(f) = f_*(\, c_1(f)^{i_1} c_2(f)^{i_2} \cdots )$$

for all $I = (i_1, i_2, \ldots)$, but by abuse of notation we also consider them as classes on $M$, namely as the pull-back
$$f^* f_*\bigl( c_1(f)^{i_1} c_2(f)^{i_2} \cdots \bigr).$$

The article [18] concentrates on the case when $\dim(M) < \dim(N)$. In the application we are interested in, $f: M \to N$ will be the map $\pi : X \to Y$ where $X$ is the critical locus of a family of curves $D$ in a family of surfaces $\pi : F \to Y$.

Kazarian's proof of the existence of Thom polynomials for multisingularities is based on a topological argument using complex cobordism theory. More recently, Ohmoto [29] has proposed a proof in the context of algebraic geometry over an algebraically closed field of characteristic zero, using as a key ingredient intersection theory on relative Hilbert schemes of points and its relation with the Thom–Mather theory of singularities of mappings. Once the existence of Thom polynomials is proved, there are various possibilities to compute them. Kazarian has offered very many such computations, see [18] as well as the references therein. Here I will give a very short selection of those.

**(7.2) Legendrian singularities.** The particular Thom polynomials we are interested in concern a special class of maps: a holomorphic map $f : M \to N$ is a *Legendrian map* if it can be written as a composed map $M \to \mathbf{P}(\Omega_N^1) \to N$, where the second map is the projection and the first map is a Legendrian immersion, meaning that $\dim M = \dim N - 1$, and the natural contact form[13] on $\mathbf{P}(\Omega_N^1)$ pulls back to 0 on $M$. The Legendrian maps form a special class of maps whose typical singularities differ from the typical singularities of the generic maps. We will be interested in a particular instance of Legendrian maps that I will describe shortly. Before that, let me mention that Thom polynomials for Legendrian maps are given in terms of a particular set of classes, namely

(7.2.1) $\qquad v = c_1\bigl(\mathcal{O}_{\mathbf{P}(\Omega_N^1)}(1)\bigr) \quad \text{and} \quad \alpha_i = c_i\bigl(f^* T_N - T_M \otimes \mathcal{O}_{\mathbf{P}(\Omega_N^1)}(1)\bigr);$

the latter description means that the $\alpha_i$ are defined by the following formal identity giving the total Chern class $\alpha$,
$$\alpha = \frac{c\bigl(f^* T_N\bigr)}{c\bigl(T_M \otimes \mathcal{O}_{\mathbf{P}(\Omega_N^1)}(1)\bigr)}.$$

The following is a Legendrian map, and it should look familiar if you came here after reading previous parts of the text. Consider a smooth family $\pi : F \to Y$, and $D$ a smooth hypersurface of the total space $F$. Denote by $X \subseteq D$ the subset at which $D$ is tangent to the fibres of $\pi$, equivalently, at which the fiber of $D$ is singular. Under a suitable transversality condition, $X$ is smooth and the restriction $\pi : X \to Y$ is a Legendrian map; the map $X \to \mathbf{P}(\Omega_Y^1)$ is defined by mapping a point $x$ to the hyperplane $\pi_* T_{D,x}$ of $T_{Y,\pi(x)}$. In this particular example, the classes in Display (7.2.1) are

(7.2.2) $\qquad v = c_1\bigl(\mathcal{O}_F(D)\bigr) \quad \text{and} \quad \alpha = c\bigl(\mathcal{O}_F(D) \otimes \Omega_{F/Y}^1\bigr) \cdot c\bigl(T_{F/Y}\bigr)^{-1}.$

Kazarian's result on the existence of Thom polynomials in this context is the following.

**(7.3) Theorem** (Kazarian)**.** *For all sequences* $\mathsf{T} = (\mathsf{T}_1, \ldots, \mathsf{T}_r)$ *of classes of isolated hypersurface singularities, there exists a universal "residual" polynomial $R_\mathsf{T}$ in the classes $v$ and $\alpha_i$, independent of the order of the sequence* $\mathsf{T}$*, and such that for all sequences* $\mathsf{S} = (\mathsf{S}_1, \ldots, \mathsf{S}_r)$ *and*

---

[13]this is a 1-form $\eta$ on $\mathbf{P}(\Omega_N^1)$ such that $(\eta \wedge d\eta)$ is nowhere vanishing; it is deduced from the Liouville form on $\Omega_N^1$, which in local coordinates writes $\sum_{i=1}^n x_i dx_i$.

*all generic Legendrian map $f : M \to N$, the classes $m_{\mathsf{S}}$ and $n_{\mathsf{S}}$, defined as in Paragraph (7.1), are given by the following formulae:*

$$m_{\mathsf{S}} = \sum_{J_1 \sqcup \cdots \sqcup J_k = [\![1,r]\!]} R_{\mathsf{S}_{J_1}} \cdot f^* f_* R_{\mathsf{S}_{J_2}} \cdots f^* f_* R_{\mathsf{S}_{J_k}};$$

$$n_{\mathsf{S}} = \sum_{J_1 \sqcup \cdots \sqcup J_k = [\![1,r]\!]} f_* R_{\mathsf{S}_{J_1}} \cdot f_* R_{\mathsf{S}_{J_2}} \cdots f_* R_{\mathsf{S}_{J_k}},$$

*where the summation ranges over all possible partitions of the set $[\![1,r]\!]$ in a disjoint union $J_1 \sqcup \cdots \sqcup J_k$ of an arbitrary number of unordered subsets, such that $J_1$ is always the subset containing 1; for all subsets $J = \{j_1, \ldots, j_s\}$, $R_{\mathsf{S}_J}$ denotes the universal residual polynomial attached to the sequence $(\mathsf{S}_{j_1}, \ldots, \mathsf{S}_{j_s})$, in any order.*

**(7.4) On the structure of the formulae.** In the theorem, the term "residual" is used for the following reason. Consider a sequence $(\mathsf{S}_1, \mathsf{S}_2)$ consisting of two classes of singularities. In order to enumerate, in $N$ say, the fibres of $f$ with a singularity of type $\mathsf{S}_1$ and a singularity of type $\mathsf{S}_2$, the naive way to proceed is to intersect the two classes $f_* R_{\mathsf{S}_1}$ and $f_* R_{\mathsf{S}_2}$, which enumerate fibers with a singularity $\mathsf{S}_1$ and fibers with a singularity $\mathsf{S}_2$, respectively. The problem however is that there is a lot of unwanted contributions in this intersection product, as has been observed at length in the previous Section 6. The residual polynomial $R_{\mathsf{S}_1, \mathsf{S}_2}$ is designed exactly to compensate for these unwanted contributions.

In fact, if one considers a sequence $\mathsf{S} = (r\mathsf{S}_0)$ consisting of $r$ times the same singularity $\mathsf{S}_0$, then the formulae of Theorem (7.3) may be phrased in terms of the Bell polynomials, namely

$$n_{r\mathsf{S}} = P_r(f_* R_{\mathsf{S}_0}, f_* R_{2\mathsf{S}_0}, \ldots, f_* R_{r\mathsf{S}_0})$$

with $P_r$ the $r$-th Bell polynomial. This is observed by Kazarian at [18, p. 683], and the computations leading to this formulation are exactly the same as those in the proof of Theorem (6.8). In brief, for $\mathsf{S} = (r\mathsf{S}_0)$, for all subsets $J \subseteq [\![1,r]\!]$, $R_{\mathsf{S}_J} = R_{|J|.\mathsf{S}_0}$, hence

$$\sum_{J_1 \sqcup \cdots \sqcup J_k = [\![1,r]\!]} f_* R_{\mathsf{S}_{J_1}} \cdot f_* R_{\mathsf{S}_{J_2}} \cdots f_* R_{\mathsf{S}_{J_k}} = \sum_{s_1 + \cdots + r s_r = r} e_{s_1, \ldots, s_r} \left( \prod_{i=1}^r f_* R_{i\mathsf{S}_0} \right)$$

where $e_{s_1, \ldots, s_r}$ is the number of partitions of $[\![1,r]\!]$ with $s_1$ subsets of size 1, …, $s_r$ subsets of size $r$. By Corollary (A.4) on Bell polynomials, the right-hand-side term in the above display is

$$P_r(f_* R_{\mathsf{S}_0}, f_* R_{2\mathsf{S}_0}, \ldots, f_* R_{r\mathsf{S}_0}),$$

as asserted.

In particular, when $\mathsf{S}_0$ is the ordinary double point singularity $\mathrm{A}_1$, Kazarian's residual classes $R_{\mathrm{A}_1}, \ldots, R_{r\mathrm{A}_1}$ coincide with the classes $b_1, \ldots, b_r$ in Theorem (1.6). You may verify this directly with example in the tables below, see also the computations in Example (6.2). For instance, from the formula in Paragraph (7.5) and the expressions in Paragraph (6.3), one finds

$$\begin{aligned} R_{\mathrm{A}_1^2} &= -v - 3\alpha_1 \\ &= -7v - 6w_1, \end{aligned}$$

which is $Q_{b_1, 2}$ from Paragraph (1.4), as we have computed in Example (6.2).

**(7.5) Residual classes of codimension at most** 4 **of Legendrian maps.** Reproduced below is the set of explicit formulae of [18, Section 9] for residual classes of Legendrian maps in the context of Theorem (7.3) above. Let me point out that the use of the classes $v$ and $\alpha_i$ of Display (7.2.1) makes the formulae rather compact. For explicit formulae of the $\alpha_i$ in terms of our usual classes $v, w_1, w_2$ in case $F$ has dimension 2 over $Y$, see Paragraph (6.3). Also, note that we are considering singularities of codimension at most 4 for the critical locus $X$, hence singularitites of codimension at most 5 for the original family of hypersurfaces $D$. The multisingularity sequences are written in multiplicative notation.

$$R_{A_1} = 1$$
$$R_{A_2} = \alpha_1$$
$$R_{A_1^2} = -v - 3\alpha_1$$
$$R_{A_3} = v\alpha_1 + 3\alpha_2$$
$$R_{A_1 A_2} = -6(v\alpha_1 + 2\alpha_2)$$
$$R_{A_1^3} = 2(v^2 + 19v\alpha_1 + 30\alpha_2)$$
$$R_{A_4} = v^2\alpha_1 + 4v\alpha_2 + 3\alpha_1\alpha_2 + 6\alpha_3$$
$$R_{D_4} = -v\alpha_2 + \alpha_1\alpha_2 - 2\alpha_3$$
$$R_{A_1 A_3} = -4(2v^2\alpha_1 + 5v\alpha_2 + 6\alpha_1\alpha_2 + 3\alpha_3)$$
$$R_{A_2^2} = -3(3v^2\alpha_1 + 8v\alpha_2 + 7\alpha_1\alpha_2 + 6\alpha_3)$$
$$R_{A_1^2 A_2} = 24(3v^2\alpha_1 + 7v\alpha_2 + 6\alpha_1\alpha_2 + 3\alpha_3)$$
$$R_{A_1^4} = -6(v^3 + 111v^2\alpha_1 + 239v\alpha_2 + 171\alpha_1\alpha_2 + 78\alpha_3)$$
$$R_{A_5} = v^3\alpha_1 - 4v^2\alpha_2 + 16v\alpha_1\alpha_2 - 12v\alpha_3 + 27\alpha_1\alpha_3 + 6\alpha_4$$
$$R_{D_5} = -2(2v^2\alpha_2 - 2v\alpha_1\alpha_2 + 7v\alpha_3 - 3\alpha_1\alpha_3 + 6\alpha_4)$$
$$R_{A_1 A_4} = -10(v^3\alpha_1 - 4v^2\alpha_2 + 14v\alpha_1\alpha_2 - 16v\alpha_3 + 21\alpha_1\alpha_3 - 6\alpha_4)$$
$$R_{A_1 D_4} = 4(5v^2\alpha_2 - 5v\alpha_1\alpha_2 + 16v\alpha_3 - 6\alpha_1\alpha_3 + 12\alpha_4)$$
$$R_{A_2 A_3} = -6(2v^3\alpha_1 - 10v^2\alpha_2 + 28v\alpha_1\alpha_2 - 39v\alpha_3 + 39\alpha_1\alpha_3 - 18\alpha_4)$$
$$R_{A_1^2 A_3} = 2(56v^3\alpha_1 - 220v^2\alpha_2 + 684v\alpha_1\alpha_2 - 951v\alpha_3 + 891\alpha_1\alpha_3 - 522\alpha_4)$$
$$R_{A_1 A_2^2} = 18(7v^3\alpha_1 - 20v^2\alpha_2 + 74v\alpha_1\alpha_2 - 96v\alpha_3 + 95\alpha_1\alpha_3 - 50\alpha_4)$$
$$R_{A_1^3 A_2} = -48(28v^3\alpha_1 - 55v^2\alpha_2 + 250v\alpha_1\alpha_2 - 318v\alpha_3 + 300\alpha_1\alpha_3 - 180\alpha_4)$$
$$R_{A_1^5} = 24(v^4 + 671v^3\alpha_1 - 701v^2\alpha_2 + 4863v\alpha_1\alpha_2 - 5844v\alpha_3 + 5490\alpha_1\alpha_3 - 3420\alpha_4)$$

**(7.6) Curves in a linear system on a fixed surface.** Theorem (7.3) may be applied to enumerate curves in a complete linear system $|L|$ on a surface $S$, arguing as in Section 3.1. The result then takes the following form, see [18, Theorem 10.1], in terms of the classes

$$d = L^2 \; ; \quad k = L \cdot K_S \; ; \quad s = K_S^2 \; ; \quad x = c_2(S).$$

**(7.7) Theorem** (Kazarian). *For all sequences* $\mathsf{T} = (\mathsf{T}_1, \ldots, \mathsf{T}_r)$ *of classes of isolated planar curve singularities, there exists a universal integral linear combination* $a_\mathsf{T}$ *of the numbers* $d, k, s, x$, *independent of the order of the sequence* $\mathsf{T}$, *and such that for all sequences* $\mathsf{S} = (\mathsf{S}_1, \ldots, \mathsf{S}_r)$, *and for all sufficiently generic complete linear systems* $|L|$ *on a surface, the number* $N_\mathsf{S}$ *of curves in* $|L|$ *with singularities of type* $\mathsf{S}$, *and passing through* $\dim |L| - \operatorname{codim}(\mathsf{S})$ *general points on the surface, is given by the formula*

$$N_\mathsf{S} = \frac{1}{\operatorname{aut}(\mathsf{S})} \sum_{J_1 \sqcup \cdots \sqcup J_k = [1,r]} a_{\mathsf{S}_{J_1}} a_{\mathsf{S}_{J_2}} \cdots a_{\mathsf{S}_{J_k}},$$

*as in Theorem (7.3).*

The number aut($\mathsf{S}$) is defined in Paragraph (7.1) above, and the expected codimension of a sequence of singularities is the same as that defined in Paragraph (4.22); see [18, Definition 1.1] for a definition in the more general context considered by Kazarian.

The genericity condition is the following. Denote by $Y_{\mathrm{reg}}$ the locus in $|L|$ of curves $C$ that are reduced, and such that the deformations of $C$ in $|L|$ induce a system of semi-universal deformations of the multigerm of singularities of $C$. Then the condition is that the complement of $Y_{\mathrm{reg}}$ in $|L|$ has codimension strictly larger than codim($\mathsf{S}$). This condition may be verified in practice with the techniques described in Section 5.1.

Now let us give the explicit formulae for the linear forms $a_{\mathsf{S}}$ for singularities of codimension up to 5, reproducing the table in [18, Section 10]. They are readily computed from those given above in Paragraph (7.5), argueing as in Section 3.1, but I felt it would be nice to include them for ease of reference. For $\mathsf{S} = A_1^r$, $a_{A_1^r}$ coincides with $a_r$ given in Paragraph (1.2) (fortunately).

$$a_{A_1} = 3d + 2k + x$$
$$a_{A_2} = 12d + 12k + 2s + 2x$$
$$a_{A_1^2} = -42d - 39k - 6s - 7x$$
$$a_{A_3} = 50d + 64k + 17s + 5x$$
$$a_{A_1 A_2} = -240d - 288k - 72s - 24x$$
$$a_{A_1^3} = 1380d + 1576k + 376s + 138x$$
$$a_{A_4} = 180d + 280k + 100s$$
$$a_{D_4} = 15d + 20k + 5s + 5x$$
$$a_{A_1 A_3} = -1260d - 1820k - 596s - 60x$$
$$a_{A_2^2} = -1260d - 1800k - 588s - 48x$$
$$a_{A_1^2 A_2} = 9000d + 12360k + 3864s + 456x$$
$$a_{A_1^4} = -72360d - 95670k - 28842s - 3888x$$
$$a_{A_5} = 630d + 1140k + 498s - 60x$$
$$a_{D_5} = 84d + 132k + 44s + 20x$$
$$a_{A_1 A_4} = -5460d - 9240k - 3740s + 200x$$
$$a_{A_1 D_4} = -420d - 624k - 196s - 100x$$
$$a_{A_2 A_3} = -6300d - 10332k - 4044s + 60x$$
$$a_{A_1^2 A_3} = 52920d + 84180k + 31816s + 240x$$
$$a_{A_1 A_2^2} = 53676d + 84456k + 31716s + 72x$$
$$a_{A_1^3 A_2} = -505008d - 770112k - 279792s - 5616x$$
$$a_{A_1^5} = 5225472d + 7725168k + 2723400s + 84384x$$

**(7.8) Contacts of hyperplanes with hypersurfaces.** Arguing as in Section 3.2, Kazarian gives applications of his Theorem (7.3) to the enumeration of hyperplanes satisfying tangency conditions with a fixed hypersurface in projective space.

In this case the Legendrian map is as follows. We consider a smooth, degree $d$, hypersurface $V \subseteq \mathbf{P}^n$. Then, we let $Y = \check{\mathbf{P}}^n$ be the dual projective space, parametrizing hyperplanes in $\mathbf{P}^n$, and $F \subseteq \mathbf{P}^n \times \check{\mathbf{P}}^n$ be the point-hyperplane incidence correspondence, i.e.,

$$F = \{(x, H) \in \mathbf{P}^n \times \check{\mathbf{P}}^n : x \in H\};$$

in Grothendieck notation, $F$ identifies with the projective bundle $\mathbf{P}(T_{\mathbf{P}^n})$. The map $\pi : F \to Y$ is the second projection. Next, we consider the divisor $D = V \times \check{\mathbf{P}}^n$ in $F$, and the critical locus $X$ of the map $\pi : D \to Y$. In this case, $X$ is the conormal variety of $V$, i.e.,

$$X = \{(x, H) \in V \times \check{\mathbf{P}}^n : H \text{ is tangent to } V \text{ at } x\},$$

and it identifies with the projectivization $\mathbf{P}(N_V)$ of the normal bundle of $V$ in $\mathbf{P}^n$.[14] Let me quote the two following sets of formulae from [18, p. 716–717], as they have been considered in this volume in [XII, C].

---

[14]Actually, it is better to see it as $\mathbf{P}(N_V(-1))$, so that the map $V \to V^\vee \subseteq \check{\mathbf{P}}^n$ is given by a linear subsystem of $|\mathcal{O}_{\mathbf{P}(N_V(-1))}(1)|$.

*(7.8.1) Contacts of lines with a plane curve.* Given a generic, smooth, degree $d$, curve $C \subseteq \mathbf{P}^2$,

$$N_{A_2} = 3d(d-2) \qquad\qquad N_{A_1^2} = \frac{1}{2}d(d-3)(d-2)(d+3)$$

are the number of flex and bitangent lines to $C$ respectively. These numbers are contained in the Plücker formulae, see [XII, Section **].

*(7.8.2) Contacts of planes with a surface in three-space.* Given a generic, smooth, degree $d$, surface $S \subseteq \mathbf{P}^3$,

$$N_{A_3} = 2d\,(11d-24)(d-2) \qquad N_{A_1 A_2} = 4d(d-3)(d-2)(d^3+3d-16)$$
$$N_{A_1^3} = \frac{1}{6}d(d-2)(d^7 - 4d^6 + 7d^5 - 45d^4 + 114d^3 - 111d^2 + 548d - 960)$$

are the number of planes cutting out on $S$, respectively, a curve with a tacnode, a curve with one node and one cusp, and a curve with three nodes. These numbers have been computed by Salmon, see [XII, Section **]: $N_{A_3}$ is the number of intersection points of the flecnodal and parabolic curves on $S$, $2N_{A_3} + N_{A_1 A_2}$ is the number of intersection points of the couple-nodal and the parabolic curves, and $N_{A_1^3}$ is the number of tritangent planes.

**(7.9) Contacts of linear spaces with hypersurfaces.** Kazarian gives more generally some formulae for the number of linear spaces of arbitrary dimension with prescribed contact conditions with a hypersurface in projective space.

In this case, the parameter space $Y$ is the Grassmannian $\mathbf{G}(k,n)$ of $k$-planes in $\mathbf{P}^n$ with the universal family $\pi : F \to Y$, the divisor is cut out by $V \times Y$ on $F \subseteq \mathbf{P}^n \times Y$, and

$$X = \big\{(x, \Lambda) \in V \times \mathbf{G}(k,n) : \ \Lambda \subseteq T_{X,x}\big\}$$

is a Grassmann bundle over $V$. I quote the following set of formulae from [18, p. 719].

*(7.9.1) Contacts of lines with a surface in three-space.* The Grassmannian of lines in $\mathbf{P}^3$ has dimension 4, and the following formulae encompass those contact conditions that directly give finitely many curves, without any further Schubert-type condition. Let $S$ be a generic, smooth, degree $d$, surface in $\mathbf{P}^3$.

$$N_{A_4} = 5d\,(7d-12)(d-4) \qquad\qquad N_{A_1^2 A_2} = \frac{1}{2}d(d-6)(d-5)(d-4)(d^3+9d^2+20d-60)$$
$$N_{A_1 A_3} = 2d\,(d-5)(3d-5)(d-4)(d+6) \qquad N_{A_1^4} = \frac{1}{12}d(d-7)(d-6)(d-5)(d-4)(d^3+6d^2+7d-30)$$
$$N_{A_2^2} = \frac{1}{2}d\,(d-5)(d-4)(d^3+3d^2+29d-60)$$

On the other hand, we have seen in [XII, Paragraph **] and [C, Theorem **] that the number of bitangent and flex-tangent lines to $S$ passing through a general point of $\mathbf{P}^3$ are respectively

$$d(d-1)(d-2) \quad\text{and}\quad \frac{1}{2}d(d-1)(d-2)(d-3);$$

they correspond to the singularities $A_1^2$ and $A_2$ respectively. Also, we have seen in [C, Theorem **] that the number of lines in $\mathbf{P}^n$ ($n > 3$) passing through a general point, and with two contacts of orders 2 and $n-1$ respectively (singularities $A_1$ and $A_{n-2}$) with a degree $d$ hypersurface, is $\prod_{k=0}^{n}(d-k)$. The following more general formula has appeared in a recent preprint.

**(7.10) Theorem** (Fehér–Juhász, [11, Theorem 4.7]). *Let $V$ be a generic, smooth, degree $d$ hypersurface in $\mathbf{P}^n$, and consider a general point $v \in \mathbf{P}^n$. For all $e_1, \ldots, e_m$ such that* $\sum_{k=1}^{m} k e_k = n - 1$, *the number of lines passing through $v$ and with tangency to $V$ corresponding to the singularity $A_1^{e_1} \cdots A_m^{e_m}$ equals*

$$\frac{1}{\prod_{k=1}^{m} e_k!} d(d-1) \cdots (d - \sum_{k=1}^{m}(k+1) e_k + 1).$$

In the context of the theorem, consider the Grassmannian of lines in $\mathbf{P}^n$: it has dimension $2(n-1)$; passing through a point is $n-1$ conditions, and the tangency conditions amount to another $n-1$ conditions.

In the same article, Fehér and Juhász give many more formulae, all built upon their main result, that for arbitrary $e_1, \ldots, e_m$, the number of lines with tangency to $V$ corresponding to $A_1^{e_1} \cdots A_m^{e_m}$ and intersecting a general linear subspace of $\mathbf{P}^n$ of codimension $\sum_{k=1}^{m} k e_k + 1$ is polynomial of degree $\sum_{k=1}^{m}(k+1) e_k$ in $d$.

# A – Basics on Bell polynomials

**(A.1) Complete and partial Bell polynomials.** Consider indeterminates $X_i$, $i \in \mathbf{N}^*$, and $T$. The *complete (exponential) Bell polynomials* $P_r$, $r \in \mathbf{N}$, are defined by the formal identity

$$\sum_{r \geqslant 0} P_r \frac{T^r}{r!} = \exp\left(\sum_{q \geqslant 1} X_q \frac{T^q}{q!}\right)$$

in $\mathbf{Q}(X_1, \ldots,)[[T]]$. The *partial Bell polynomials* $B_{r,k}$, $r, k \in \mathbf{N}$, $r \geqslant k$, on the other hand are defined by the formal identities

$$\sum_{r \geqslant k} B_{r,k} \frac{T^r}{r!} = \frac{1}{k!} \left(\sum_{q \geqslant 1} X_q \frac{T^q}{q!}\right)^k.$$

It follows directly that for all $r \in \mathbf{N}$, one has

$$P_r = \sum_{k \leqslant r} B_{r,k}.$$

**(A.2) Lemma.** *The polynomial $B_{r,k}$ involves only the indeterminates $X_1, \ldots, X_{r-k+1}$, and*

$$B_{r,k} = \sum_{\substack{i_1 + \cdots + i_{r-k+1} = k \\ i_1 + \cdots + (r-k+1)i_{r-k+1} = r}} \frac{r!}{i_1! \cdots i_{r-k+1}!} \frac{1}{(1!)^{i_1} \cdots ((r-k+1)!)^{i_{r-k+1}}} X_1^{i_1} \cdots X_{r-k+1}^{i_{r-k+1}}.$$

*Proof.* By Newton's multinomial formula,

$$\frac{1}{k!}\left(\sum_{q \geqslant 1} X_q \frac{T^q}{q!}\right)^k = \sum_{i_1+i_2+\cdots=k} \frac{1}{i_1! i_2! \cdots} \left(X_1 \frac{T}{1!}\right)^{i_1} \left(X_2 \frac{T^2}{2!}\right)^{i_2} \cdots$$

where the sum ranges over all sequences $(i_j)_{j \in \mathbf{N}^*}$ of non-negative integers such that $\sum_{j \geqslant 1} i_j = k$. The sequences that contribute to the $T^r$ term are those verifying in addition $\sum_{j \geqslant 1} j.i_j = r$; this implies $i_j = 0$ for all $j > r - k + 1$. Then the formula for $B_{r,k}$ follows directly. □

**(A.3) Lemma.** *Let $r, k \in \mathbf{N}$ such that $k \leqslant r$. For all non-negative integers $i_1, \ldots, i_{r-k+1}$ such that $i_1 + \cdots + i_{r-k+1} = k$ and $i_1 + \cdots + (r-k+1)i_{r-k+1} = r$, the number*

$$\frac{r!}{i_1! \cdots i_{r-k+1}!} \frac{1}{(1!)^{i_1} \cdots ((r-k+1)!)^{i_{r-k+1}}}$$

*is the number of unordered partitions of $[\![1, r]\!]$ into $k$ unordered blocks, of which $i_1$ have size 1, ..., $i_{r-k+1}$ have size $r - k + 1$.*

*Proof.* First, the number

$$\frac{r!}{(1!)^{i_1} \cdots ((r-k+1)!)^{i_{r-k+1}}}$$

is the number of partitions of $[\![1, r]\!]$ into $i_1$ blocks of size 1, ..., $i_{r-k+1}$ blocks of size $r - k + 1$, where there is an ordering of the blocks, but each block does not have an inner ordering. Indeed, there is a bijection between permutations of $[\![1, r]\!]$ and ordered partitions of $[\![1, r]\!]$ into $i_1$ ordered blocks of size 1, ..., $i_{r-k+1}$ ordered blocks of size $r - k + 1$: for each permutation, the ordered partition is simply obtained by putting the first term in the first block and so forth until $i_1$ blocks of size 1 have been chosen, then the next two terms in that order in the first block of size 2, and so on. Then, forgetting the inner ordering of each block, one obtains a total of

$$\frac{r!}{(1!)^{i_1} \cdots ((r-k+1)!)^{i_{r-k+1}}}$$

ordered partitions in unordered blocks. Finally, forgetting the ordering of the blocks, one obtains the asserted number. $\square$

The two above lemmas (A.2) and (A.3) have the following direct corollaries.

**(A.4) Corollary.** *The complete Bell polynomial $P_r$ involves only the indeterminates $X_1, \ldots, X_r$, and for all $(i_1, \ldots, i_r)$, the coefficient of its $X_1^{i_1} \cdots X_r^{i_r}$ term is the number of partitions of $[\![1, r]\!]$ into $i_1$ blocks of size 1, ..., $i_r$ blocks of size $r$.*

**(A.5) Corollary.** *The partial Bell polynomial $B_{r,k}$ has non-negative integer coefficients. It is homogeneous of degree $k$ with respect to the standard grading on $\mathrm{Z}[X_1, \ldots, X_{r-k+1}]$. It is homogeneous of degree $r$ with respect to the grading defined by assigning degree $i$ to $X_i$ for all $i$.*

**(A.6) Corollary.** *The complete Bell polynomial $P_{r,k}$ has non-negative integer coefficients. It is homogeneous of degree $r$ with respect to the grading defined by assigning degree $i$ to $X_i$ for all $i$.*

Finally, let us prove the following indentities.

**(A.7) Proposition.** *Let $r \geqslant 0$ be an integer, and consider indeterminates $X_1, \ldots, X_{r+1}$ and $Y_1, \ldots, Y_{r+1}$. The two following identitities hold:*

(A.7.1) $$P_{r+1}(X_1, \ldots, X_{r+1}) = \sum_{s=0}^{r} \binom{r}{s} X_{r-s+1} P_s;$$

(A.7.2) $$P_r(X_1 + Y_1, \ldots, X_r + Y_r) = \sum_{s=0}^{r} \binom{r}{s} P_s(X_1, \ldots) P_{r-s}(Y_1, \ldots).$$

*Proof.* To prove the first identity, let us differentiate with respect to the indeterminate $T$ the formal identity,

$$\sum_{r \geqslant 0} P_r \frac{T^r}{r!} = \exp\left(\sum_{q \geqslant 1} X_q \frac{T^q}{q!}\right);$$

one finds

$$\sum_{r\geqslant 0} P_{r+1}\frac{T^r}{r!} = \bigg(\sum_{q\geqslant 0} X_{q+1}\frac{T^q}{q!}\bigg) \exp\bigg(\sum_{q\geqslant 1} X_q\frac{T^q}{q!}\bigg)$$
$$= \bigg(\sum_{q\geqslant 0} X_{q+1}\frac{T^q}{q!}\bigg)\bigg(\sum_{r\geqslant 0} P_r\frac{T^r}{r!}\bigg)$$
$$= \sum_{r\geqslant 0}\bigg(\sum_{p+q=r}\frac{1}{p!q!}X_{q+1}P_p\bigg)T^r,$$

and (A.7.1) follows. For the other one, one computes

$$\exp\bigg(\sum_{q\geqslant 1}(X_q+Y_q)\frac{T^q}{q!}\bigg) = \exp\bigg(\sum_{q\geqslant 1} X_q\frac{T^q}{q!}\bigg)\cdot\exp\bigg(\sum_{q\geqslant 1} Y_q\frac{T^q}{q!}\bigg)$$
$$= \bigg(\sum_{r\geqslant 0} P_r(X_1,\ldots)\frac{T^r}{r!}\bigg)\bigg(\sum_{r\geqslant 0} P_r(Y_1,\ldots)\frac{T^r}{r!}\bigg)$$
$$= \sum_{r\geqslant 0}\bigg(\sum_{p+q=r}\frac{1}{p!q!}P_p(X_1,\ldots)P_q(Y_1,\ldots)\bigg)T^r,$$

and (A.7.2) follows. $\square$

## References

**Chapters in this Volume**[a]

[II]   E. Sernesi, *Deformations of curves on surfaces.*

[XI]   Th. Dedieu, *Enumerative geometry of K3 surfaces.*

[XII]  Th. Dedieu, *Some classical formulæ for curves and surfaces.*

[XIV]  F. Bastianelli and Th. Dedieu, *Göttsche conjecture and Göttsche–Yau–Zaslow formula.*

[C]    L. Busé and Th. Dedieu, *Generalized weight properties of resultants and discriminants, and applications to projective enumerative geometry.*

**External References**

[1] A. B. Altman, A. Iarrobino and S. Kleiman, *Irreducibility of the compactified jacobian*, in *Real and Complex Singularities (Proc. Ninth Nordic Summer School/NAVF Sympos. Math., Oslo, 1976)*, Sijthoff and Noordhoff, Alphen aan den Rijn, 1977, 1–12.

[2] V. I. Arnold, *Local normal forms of functions*, Invent. Math. **35** (1976), 87–109.

[3] E. Artal Bartolo and P. Cassou-Noguès, *Milnor number of plane curve singularities in arbitrary characteristic*, arXiv:2409.13520.

[4] M. Beltrametti and A. J. Sommese, *Zero cycles and k-th order embeddings of smooth projective surfaces*, in *Problems in the theory of surfaces and their classification. Papers from the meeting held at the Scuola Normale Superiore, Cortona, Italy, October 10-15, 1988*, London: Academic Press; Rome: Istituto Nazionale di Alta Matematica Francesco Severi, 1991, Appendix by Lothar Göttsche: Identification of very ample line bundles on $S^{(r)}$, 33–48; appendix: 44–48.

---

[a]see https://www.math.univ-toulouse.fr/~tdedieu/#enumTV

[5] F. Block, *Computing node polynomials for plane curves*, Math. Res. Lett. **18** (2011), no. 4, 621–643.

[6] ______, *Relative node polynomials for plane curves*, J. Algebraic Combin. **36** (2012), no. 2, 279–308.

[7] E. Casas-Alvero, *Singularities of plane curves*, Lond. Math. Soc. Lect. Note Ser., vol. 276, Cambridge: Cambridge University Press, 2000.

[8] Th. Dedieu and E. Sernesi, *Equigeneric and equisingular families of curves on surfaces*, Publ. Mat., Barc. **61** (2017), no. 1, 175–212.

[9] D. Eisenbud and J. Harris, *3264 and all that*, Cambridge: Cambridge University Press, 2016.

[10] F. Enriques and O. Chisini, *Lezioni sulla teoria geometrica delle equazioni e delle funzioni abgebriche.*, Zanichelli, Bologna, 1934.

[11] L. M. Fehér and A. P. Juhász, *Plücker formulas using equivariant cohomology of coincident root strata*, arXiv:2312.06430.

[12] S. Fomin and G. Mikhalkin, *Labeled floor diagrams for plane curves*, J. Eur. Math. Soc. (JEMS) **12** (2010), no. 6, 1453–1496.

[13] W. Fulton, *Intersection theory.*, 2nd ed. ed., Ergeb. Math. Grenzgeb., 3. Folge, vol. 2, Berlin: Springer, 1998.

[14] G.-M. Greuel and C. Lossen, *Equianalytic and equisingular families of curves on surfaces*, Manuscr. Math. **91** (1996), no. 3, 323–342.

[15] G.-M. Greuel, C. Lossen and E. Shustin, *Introduction to singularities and deformations*, Springer Monogr. Math., Berlin: Springer, 2007.

[16] ______, *Singular algebraic curves*, Springer Monogr. Math., Berlin: Springer, 2018, with an appendix by O. Viro.

[17] S. Katz, *Rational curves on Calabi-Yau threefolds*, in *Essays on mirror manifolds*, Cambridge, MA: International Press, 1992, 168–180.

[18] M. È. Kazaryan, *Multisingularities, cobordisms, and enumerative geometry.*, Russ. Math. Surv. **58** (2003), no. 4, 665–724.

[19] S. Kleiman, *Multiple-point formulas. I: Iteration*, Acta Math. **147** (1981), 13–49.

[20] S. Kleiman and R. Piene, *Enumerating singular curves on surfaces*, in *Algebraic geometry: Hirzebruch 70. Proceedings of the algebraic geometry conference in honor of F. Hirzebruch's 70th birthday, Stefan Banach International Mathematical Center, Warszawa, Poland, May 11–16, 1998*, Providence, RI: American Mathematical Society, 1999, 209–238.

[21] ______, *Node polynomials for families: Methods and applications*, Math. Nachr. **271** (2004), 69–90.

[22] ______, *Enriques diagrams, arbitrarily near points, and Hilbert schemes*, Atti Accad. Naz. Lincei, Cl. Sci. Fis. Mat. Nat., IX. Ser., Rend. Lincei, Mat. Appl. **22** (2011), no. 4, 411–451, with Appendix B by I. Tyomkin.

[23] ______, *Node polynomials for curves on surfaces*, SIGMA, Symmetry Integrability Geom. Methods Appl. **18** (2022), paper 059, 23.

[24] S. Kleiman and V. V. Shende, *On the Göttsche threshold (with an appendix by Ilya Tyomkin).*, in *A celebration of algebraic geometry. A conference in honor of Joe Harris' 60th birthday, Harvard University, Cambridge, MA, USA, August 25–28, 2011*, Providence, RI: American Mathematical Society (AMS); Cambridge, MA: Clay Mathematics Institute, 2013, 429–449.

[25] A. L. Knutsen, *On kth-order embeddings of K3 surfaces and Enriques surfaces*, Manuscr. Math. **104** (2001), no. 2, 211–237.

[26] M. Kool, V. Shende and R. P. Thomas, *A short proof of the Göttsche conjecture*, Geom. Topol. **15** (2011), no. 1, 397–406.

[27] T. Laarakker, *The Kleiman-Piene conjecture and node polynomials for plane curves in $\mathbb{P}^3$*, Sel. Math., New Ser. **24** (2018), no. 5, 4917–4959.

[28] J. Li and Y.-J. Tzeng, *Universal polynomials for singular curves on surfaces*, Compos. Math. **150** (2014), no. 7, 1169–1182.

[29] T. Ohmoto, *Universal polynomials for multi-singularity loci of maps*, 2024, preprint arXiv:2406.12166.

[30] N. Qviller, *The di Francesco-Itzykson-Göttsche conjectures for node polynomials of $\mathbb{P}^2$*, Int. J. Math. **23** (2012), no. 4, 19, Id/No 1250049.

[31] ______, *Structure of node polynomials for curves on surfaces*, Math. Nachr. **287** (2014), no. 11-12, 1394–1420.

[32] J. V. Rennemo, *Universal polynomials for tautological integrals on Hilbert schemes*, Geom. Topol. **21** (2017), no. 1, 253–314.

[33] T. Sasajima and T. Ohmoto, *Classical formulae on projective surfaces and 3-folds with ordinary singularities, revisited*, Saitama Math. J. **31** (2017), 141–159.

[34] ______, *Thom polynomials in $\mathcal{A}$-classification I: counting singular projections of a surface*, in *Schubert varieties, equivariant cohomology and characteristic classes—IMPANGA 15*, EMS Ser. Congr. Rep., Eur. Math. Soc., Zürich, 2018, 237–259.

[35] E. I. Shustin, *On manifolds of singular algebraic curves*, Sel. Math. Sov. **10** (1983), no. 1, 27–37.

[36] R. P. Stanley, *Enumerative combinatorics. Vol. 1.*, 2nd ed., Camb. Stud. Adv. Math., vol. 49, Cambridge: Cambridge University Press, 2012.

[37] Y.-J. Tzeng, *A proof of the Göttsche–Yau–Zaslow formula*, J. Differential Geom. **90** (2012), 439–472.

[38] I. Vainsencher, *Enumeration of $n$-fold tangent hyperplanes to a surface*, J. Algebr. Geom. **4** (1995), no. 3, 503–526.